\documentclass[11pt]{article}
\usepackage{graphicx}
\usepackage{amssymb,amsfonts,amsthm}

\usepackage{amsmath,amsthm,amssymb}

\usepackage{amscd}
\usepackage{amsfonts}
\usepackage{color}
\newtheorem{theorem}{Theorem}[section]

\newtheorem{corollary}[theorem]{Corollary}

\newtheorem{lemma}[theorem]{Lemma}
\newtheorem{proposition}[theorem]{Proposition}

\theoremstyle{definition}

\newtheorem{example}[theorem]{Example}
\newtheorem{problem}[theorem]{Problem}

\newtheorem{remark}[theorem]{Remark}

\newtheorem{definition}[theorem]{Definition}

\def\al{\alpha}                       
\def\ga{\gamma}

\newcommand{\RAAGS}{{right angled Artin groups\ }}
\newcommand{\RAAG}{{right angled Artin group\ } }

\newcommand\iv{^{-1}}

\newcommand\SL{{\mathrm{SL}}}

\newcommand\opp{^\textsl{opp}}

\newcommand\Z{{\mathbf Z}}

\newcommand\Lk{\textsl{Lk}}
\newcommand\cont{\textsl{cont}}
\newcommand\effcont{\textsl{econt}}
\newcommand\rcont{\textsl{rcont}}

\newcommand\Remark{\paragraph{Remark.}}

\newcommand{\X}{{\mathcal X}}
\newcommand\Star{\mathrm{Star}}

\newcommand\<{\langle}\def\>{\rangle}

\begin{document}

\title{Surface subgroups of right-angled Artin groups}
\author{John Crisp, Michah Sageev, Mark Sapir\thanks{The research of the third author was supported in part by NSF and BSF (the US-Israeli) grants.}}

\maketitle
\begin{abstract}
We consider the question of which right-angled Artin groups contain
closed hyperbolic surface subgroups. It is known that a right-angled
Artin group $A(K)$ has such a subgroup if its defining graph $K$
contains an $n$-hole (i.e. an induced cycle of length $n$) with
$n\geq 5$. We construct another eight ``forbidden" graphs and show that every graph $K$ on $\le
8$ vertices either contains one of our examples, or contains a hole
of length $\ge 5$, or has the property that $A(K)$ does not contain
hyperbolic closed surface subgroups. We also provide several sufficient conditions for a \RAAG to contain no hyperbolic surface subgroups.

We prove that for one of these
``forbidden" subgraphs $P_2(6)$, the right angled Artin group
$A(P_2(6))$ is a subgroup of a (right angled Artin) diagram group. Thus we
show that a diagram group can contain a non-free hyperbolic subgroup
answering a question of Guba and Sapir. We also show that fundamental groups of non-orientable surfaces can be subgroups of diagram groups. Thus the first integral homology of a subgroup of a diagram group can have torsion (all homology groups of all diagram groups are free Abelian by a result of Guba and Sapir).
\end{abstract}

\maketitle

\tableofcontents
\section{Introduction}

By a graph in this paper, we always mean finite non-oriented graph
without multiple edges or loops (edges whose initial and terminal
vertices coincide). Let $K$ denote a graph with the vertex set $K^0$
and edge set $K^1$.
We define the \emph{right-angled Artin group} 
$A(K)$ associated to $K$ to be the group with presentation
\[
A(K) = \<\ K^0\ \mid\ ab=ba \text{ if } [a,b]\in K^1\ \> \,.
\]
Such groups are sometimes referred to as \emph{graph groups} or {\em
partially commutative groups} in the literature.

\medskip

Much is already known about right-angled Artin groups and their
subgroups. For example:

\begin{itemize}
\item $A(K)\cong A(L)$ only when $K \cong L$ \cite{Droms}.

\item Every right-angled Artin group is bi-automatic, and so has solvable word and
conjugacy problem \cite{Ch}, and  does not contain nilpotent non-Abelian
subgroups.

\item Every right-angled Artin group is linear (since it is commensurable to a right-angled Coxeter group)
\cite{DJ}. In fact it embeds into $\SL_n(\Z)$ for some $n$ \cite{HW}.

\item Non-uniform (and many uniform) hyperbolic lattices embed into right-angled Artin groups
\cite{HaglundWise}.

\item If $K$ contains a hole of length $\ge 5$ then
$A(K)$ contains a copy of a hyperbolic surface group. \cite{DSS}
\end{itemize}

In this paper, we study the following problem:

\begin{problem}\label{1}
For which graphs $K$  does the right-angled Artin group $A(K)$
contain a hyperbolic surface subgroup?
\end{problem}

By a \emph{hyperbolic surface} we mean closed, compact surface of negative Euler characteristic. Let $S_g$ denote the orientable hyperbolic surface of genus $g\geq 2$.

Homomorphisms of a surface group $\pi_1(S_g)$ into another group $G$  are closely related to solving quadratic equations in $G$, an important area of group theory. For example there exists a well known reformulation due to Stallings and Hempel (see \cite{Stal, Hem, LS, GK} of the Poincare conjecture in terms of homomorphisms of the fundamental group of a closed hyperbolic surface $S_g$ of genus $g$ into $F_g\times F_g$ where $F_g$ is the free group of rank $g$ (the idea goes back to Maskit and Papakyriakopoulous, and was explored by Jaco, Waldhausen, Olshanskii and others). Quadratic equations, and homomorphisms of $\pi_1(S_g)$ into free groups play the key role in Makanin-Razborov theory (see \cite{Raz}, and \cite{GK}). A description of all solutions of quadratic equations in $F_g$ and other hyperbolic groups is obtained by Lysenok and Grigorchuk-Lysenok \cite{L,L1,GL}.

As was pointed out by Olshanskii \cite{Ol89} homomorphisms of a surface group $\pi_1(S_g)$ into a group $G$ given by presentation $\langle X\mid R\rangle$ are in natural correspondence with van Kampen diagrams on $S_g$. In fact many results about solutions of quadratic equations in groups, and homomorphisms of surface groups are most conveniently described in terms of diagrams on surfaces. For example a result of Lysenok \cite{L1} can be formulated as follows: for every hyperbolic group $G$ there exist only finitely many ``minimal" diagrams over $G$ on a surface $S_g$ up to the action by the mapping class group of $S_g$ and ${\mathrm{Aut}}(G)$.

In the case of \RAAGS, it is more convenient sometimes to study the dual pictures, i.e. {\em dissection curve diagrams} on surfaces previously used by Crisp and Wiest in \cite{CrispWiest}.
Given a graph $K$ and a compact surface with boundary $(S,\partial S)$ , a $K$-dissection on $(S,\partial S)$ consists of a collection of essential simple closed curves and properly embedded arcs ( i.e: arcs which intersect the boundary in their endpoints). The curves and arcs are transversally oriented, and
labeled by vertices of $K$ and two curves or arcs intersect only if
their labels are adjacent in $K$. For any choice of basepoint $x\in S$, the
homomorphism $\phi_x\colon \pi_1(S,x)\to A(K)$ corresponding to a given
dissection diagram is very natural: given an element $a\in \pi_1(S,x)$ represented by the loop $\alpha$ at $x$, we go along $\alpha$ (starting at
the base point $x$) reading off the labels of the dissection curves
we cross. The resulting word represents the element $\phi_x(a)\in A(K)$.

Constructing a dissection diagram for an injective homomorphism of $\pi_1(S_g)$ into a \RAAG or proving that such a diagram does not exist is usually a difficult task. This might be expected in view of the relationship already mentioned between the Poincar\'e conjecture and the homomorphisms of $\pi_1(S_g)$ into the relatively straightforward \RAAG $F_g\times F_g$ together with their associated dissection diagrams (see also the ``simple closed curve in the kernel" conjecture in \cite{Stal}).

We provide both negative and positive results for Problem \ref{1}.
One result on the negative side is that whenever the graph $K$ is chordal, i.e. does not contain induced cycles of length $\geq 3$, the right-angled Artin group $A(K)$ has no subgroup isomorphic to a hyperbolic surface group. A proof of this is given in Section \ref{chordal} where in fact we prove the much stronger result that $A(K)$ has no one-ended hyperbolic subgroup at all if $K$ is chordal.
In addition to this result, we develop a set of "reduction moves" which allow us to
reduce the question about a given graph to the same question about a simpler graph.
These reduction moves are the contents of Sections \ref{seplemma}, and \ref{redmoves}.
While some of these reduction moves are a bit involved (see Proposition
\ref{sepsub} and Section \ref{8}), a prototype to keep in mind is the following: if $K$
contains a cut point, which is the intersection of two subgraphs
$K_1$ and $K_2$ and if $A(K)$ contains a hyperbolic surface
subgroup, then either $A(K_1)$ or $A(K_2)$ contains a hyperbolic
surface subgroup. To find out how the reduction moves work, one can first read the definitions and statements of Section \ref{ndsets}, \ref{doubling},  and then read Sections \ref{examples} and \ref{8} where many examples are given.

On the other hand, we prove that many right-angled Artin groups
which do not contain long holes and so are not covered by \cite{DSS} contain hyperbolic surface
subgroups. Note that if a graph $L$ contains $K$ as an induced
subgraph, then $A(K)<A(L)$. Thus whenever one shows that a graph
contains a surface subgroup, one has shown this for any graph
containing the original graph as an induced subgraph. Recently,  some results about embeddings of surface subgroups into right
angled Artin groups were obtained by Kim \cite{Kim1}, who also showed that the right-angled Artin group associated to the triangular prism contains a surface subgroup. See
Section \ref{kim} for a brief discussion these results.
In particular, in Section \ref{embedding} we find eight new minimal ``forbidden" graphs $P_1(6), P_2(6), P_1(7), P_2(7), P_1(8)-P_4(8)$ (we write the number of vertices in parentheses): the \RAAGS corresponding to these graphs contain hyperbolic surface subgroups but \RAAGS corresponding to proper subgraphs do not. Here $P_1(6)$ is the triangular prism (= the anti-hole of length 6), $P_2(6)$ is the prism with a diagonal (= the complement of a path of length 5). The embedding results follow from a general statement (Proposition \ref{thcuts}) which
allows one to check if a given $K$-dissection diagram corresponds to
an injective homomorphism $\pi_1(S)\to A(K)$.

As an application of our results, we prove

\begin{theorem}\label{8vertex}
For every graph $K$ with at most 8 vertices, either
$A(K)$ contains one of our ``forbidden" subgraphs or $A(K)$ does not
contain a hyperbolic surface subgroup and can be reduced to the empty graph using one of the reduction moves.
\end{theorem}

A computer assisted proof of Theorem \ref{8vertex} is in Section \ref{8}.

Although this theorem is of obvious limited strength, it shows that our methods are powerful enough to deal with large classes of \RAAGS (there are more than 400 eight vertex graphs that do not contain long holes and our forbidden subgraphs, and such that the corresponding \RAAG is not decomposable into a free or direct product). There are more than 50 of these graphs that require the full strength of our reduction moves to be completely reduced.

We do not know how close we are to a complete answer to Problem \ref{1}.
We do not know any graph $K$ such that  $A(K)$ does not contain
hyperbolic surface subgroups and $K$ cannot be reduced to a one-vertex
graph by our reduction moves. We also do not know if there exists a
$K$ such that $A(K)$ contains a hyperbolic surface subgroup but $A(K)$
does not contain right angled Artin groups corresponding to $n$-holes with $n\ge 5$ or to one of our ``forbidden"
subgraphs $P_1(6)-P_4(8)$.

Theorem \ref{8vertex} has some unexpected applications.
First we show that $A(P_2(6))$ is a subgroup of a diagram group. Thus we prove the following statement answering a question by Guba and Sapir. For the definition of diagram groups, and for the motivation see Section \ref{dg}.

\begin{theorem} A diagram group can contain a hyperbolic surface subgroup. In particular not every hyperbolic group faithfully representable by diagrams is free.
\end{theorem}

Another corollary deals with the following question.

\begin{problem}
Given two finite graphs $K, K'$, decide whether the \RAAG $A(K)$ embeds into the  \RAAG $A(K')$.
\end{problem}
The question is still wide open. The only known obstacles for embedding of $A(K)$ into $A(K')$ are the following:

\begin{itemize}
\item If $K$ contains a clique of size $n$ and $K'$ does not contain a clique of
this size then $A(K)$ cannot be a subgroup of $A(K')$.
\item If $A(K)$ contains a hyperbolic surface subgroup and $A(K')$ does not,
then $A(K)$ cannot be a subgroup of $A(K')$.
\end{itemize}

In particular, it is not clear whether for each of our ``forbidden" graph $P$, the \RAAG $A(P)$ contains the \RAAG $A(C_n)$ corresponding to a long hole $C_n, n\ge 5$. In fact, Kim \cite{Kim1} showed that $A(P_1(6))$ contains $A(C_5)$. On the other hand it is known \cite{GS2} that diagram groups cannot contain $A(C_n)$ for odd $n>3$. Hence we obtain that $A(P_2(6))$ cannot contain $A(C_n)$ for any odd $n>3$. It is not clear how to prove such a result directly. It is also not clear whether $A(P_2(6))$ contains $A(C_n)$ for even $n>5$ (it is not known whether a diagram group can contain $A(C_n)$ for even $n>5$). Kim also showed (using Remark \ref{doubleRemark} below) that $A(P_2(6))$ contains $A(P_2(7))$ as a subgroup. So in fact we currently have only 6 ``essential" forbidden graphs (and $C_n$, $n\ge 5$). 

One way to continue would be to describe graphs which do not contain
long holes and our exceptional forbidden subgraphs in some algebraic way
using splittings over ``simple" subgraphs, and then try to prove
that any such graph can be simplified by one of our reduction moves.

A step in that direction has been done (upon our request) by M.
Chudnovsky and P. Seymour. They proved that any graph $K$ which does
not contain $n$-holes with $n\ge 5$ and induced copies of $P_1(6),
P_2(6)$ admits a skew partition, i.e. it non-trivially splits as an amalgam over a
join of two non-empty subgraphs. Slightly modifying their proof we prove the following stronger ``if and only if" statement.

For every subset $W$ of $K^0$, let $C(W)$ be the set (and the subgraph spanned by this set) of common neighbors of $W$. We say that a set of vertices $L\subseteq K^0$ {\em separates} vertices $u,v$ if $u$ and $v$ are different connected components of $K\setminus L$.

\begin{theorem}\label{thCW} A graph $K$
does not contain holes of length $\ge 5$ and induced subgraphs $P_1(6), P_2(6)$ if and only if for every two vertices $u,v$ at distance 2 and every co-component $W$ of $C(\{u,v\})$, the set $W \cup (C(W)\setminus \{u,v\})$ separates $u$ from $v$.
\end{theorem}

This theorem and the result of Chudnovsky and Seymour quoted above (see Lemma \ref{thCW2} below) show that a graph that does not contain long holes and induced copies of $P_1(6), P_2(6)$ can be constructed from the 1-vertex graphs by applying the following operations:

\begin{itemize}
\item taking the disjoint union of two graphs;
\item taking the join of two graphs;
\item amalgamating two graphs along a common subgraph that is a join
of two proper subgraphs.
\end{itemize}

That allows one to deal with these graphs using induction on their ``complexity" (the number of steps in their construction from 1-vertex graphs) because these correspond to direct and free products of \RAAGS. Hyperbolic surface subgroups cannot appear after steps of the first two types (joins and disjoint unions). Some types of amalgams also behave well in this respect (see Lemma \ref{product} below), but in general the situation is not clear. Six of our eight ``forbidden" graphs ($P_1(7)-P_4(8)$) are amalgams of smaller graphs over complete bi-partite graphs. It might be worthwhile  to start studying amalgams over non-trivial joins with amalgams over complete bi-partite graphs.
One more potential way of solving Problem \ref{1} is to establish result similar to Theorem \ref{thCW} for the smaller class of graphs avoiding all our ``forbidden" graphs.

We end the introduction with a few open problems. Positive solution of the first two of them would greatly advance our understanding of Problem \ref{1}. A solution of the third may shed some light on another well known open problem in group theory.

\begin{problem} Suppose that a graph $K$ splits as an amalgam of two proper subgraphs $K_1$, $K_2$ over a clique $L$. Suppose further that the \RAAG $A(K)$ contains a hyperbolic surface subgroup. Is it true that $A(K_1)$ or $A(K_2)$ also contain a hyperbolic surface subgroup? The question is open even when $K_1$ or $K_2$ is a clique itself.
\end{problem}

\begin{problem} Suppose that a graph $K$ contains two vertices $a, b$ with the same links and $A(K)$ contains a hyperbolic surface subgroup. Does it imply that $A(K\setminus \{a\})$ contains a hyperbolic surface subgroup.
\end{problem}

The next problem is an analog of the well known ``simple curve in the kernel" problem for 3-manifolds.

\begin{problem}\label{Per} Suppose that $A(K)$ contains no hyperbolic surface subgroup. Is it true that for every homomorphism $\phi\colon\pi_1(S_g)\to A(K)$, $g\ge 2$, there exists an essential simple closed curve on $S_g$ in the kernel of $\phi$?
\end{problem}

Note that the result is true when $A(K)$ is a free group (i.e. $K$ has no edges) \cite{Stal}. For complete bipartite graphs, i.e. when $A(K)$ is a direct product of two non-trivial free group, the problem is equivalent to the Poincar\'e conjecture \cite{Stal} (and so currently the only way to solve Problem \ref{Per} in this case is by using Ricci flows on 3-manifolds \cite{Per1, Per2}).

\begin{problem} Suppose that $A(K)$ does not contain hyperbolic surface subgroups. Is it true that every hyperbolic subgroup of $A(K)$ is free?
\end{problem}

Note that some of our methods of proving that a \RAAG does not contain hyperbolic surface subgroups can be (with some effort) generalized to prove that the \RAAG does not contain non-free hyperbolic subgroups at all. But some of the methods we employ use very specific properties of surfaces. For example, consider the graph $K$ presented on the picture of Step 8 in the proof of Theorem \ref{8vertex} in Section \ref{8}.  We prove that $A(K)$ does not contain hyperbolic surface subgroups. Does $A(K)$ contain non-free hyperbolic subgroups? If the answer is negative, we get a new method of proving non-existence of hyperbolic subgroups in \RAAGS. If the answer is positive we would get a negative solution of the well known Gromov's problem: does every 1-ended hyperbolic group contain a hyperbolic surface subgroup. Note that the same graph is a candidate for a counterexample to Problem \ref{Per}.

\section{Preliminaries}

\subsection{Terminology related to graphs}

We are going to use standard graph theory terminology. Here we
collect some of the terms. Let $K$ be a graph with vertex set $K^0$
and edge set $K^1$.

\begin{itemize}
\item a {\em subgraph of $K$ induced by a set of vertices $V$} is the
graph with vertex set $V$ and edge set $(V\times V)\cap K^1$
\item {\em the complementary graph} $K\opp$ is the graph with vertex set
$K^0$ where two vertices are adjacent if and only if they are not
adjacent in $K$;
\item a {\em clique} is a set of pairwise adjacent vertices of $K$;
\item a {\em stable set} is a set of pairwise non-adjacent
vertices of $K$ (i.e. it is a clique in $K\opp$);
\item a {\em hole} is an induced subgraph that is a cycle;
\item an {\em anti-hole} is a hole in $K\opp$;
\item a {\em (connected) component} of $K$ is a maximal connected
subgraph of $G$;
\item a {\em an anti-component} of $K$ is a component of $K\opp$.
\item a vertex $v$ in a subset of vertices $V$ of $K$ is called {\em
central in} $V$ if $v$ is adjacent to every other vertex in $V$.
\end{itemize}

\subsection{{Terminology related to curves and surfaces}}

{Let $S$ be a compact surface with boundary. \begin{itemize}
\item By an {\it
essential arc} in $S$, we mean a map $\alpha: [0,1]\to S$, with
$\alpha(0), \alpha(1)\in\partial S$ which is not homotopic relative
to the boundary into $\partial S$.  (We do not require arcs to be
embedded).

\item By an {\it essential closed curve} we mean a closed curve which is non-trivial in $\pi_1(S)$ and which cannot be homotoped into any boundary component.

\item We say that a collection of closed curves and arcs is in {\em minimal position} if for any two curves $\alpha$ and $\beta$, the number of intersections of $\alpha$ and $\beta$ is minimal among all curves $\alpha', \beta'$ where $\alpha'$ is homotopic to $\alpha$ and $\beta'$ is homotopic to $\beta$ (we consider free homotopies for closed curves and homotopies relative to the end points for arcs).
\end{itemize}

\subsection{The dissection diagrams}
\label{dc}

Let $(S,\partial S)$ be a surface with (possibly empty) boundary. Let $G=\langle X\mid
R\rangle$ be a finitely presented group. Let $\Psi$ be a {\em van Kampen
diagram} over the presentation of $G$ drawn on $S$. That is a polyhedral decomposition of $S$ with a cellular map into the presentation complex of $G$.
In other words, the diagram $\Psi$
is a graph drawn on $S$ with edges labeled by letters from $X$,
such that each connected component of $S\setminus \Psi^1$ is a
polygon with boundary path labeled by a word from $R^{\pm 1}$ (see
more details in \cite{book, LS}).

Given a van Kampen diagram $\Psi$ on $S$, one can define a
homomorphism $\phi\colon \pi_1(S)\to G$ as follows. As a base-point,
pick a vertex $v$ of $\Psi$. Let $\gamma$ be any loop at $v$. Since
all cells in the tessellation $\Psi$ are polygons, $\gamma$ is
homotopic to a curve that is a composition of edges of $\Psi$. Then
$\phi(\gamma)$ is the word obtained by reading the labels of edges
of $\Psi$ along $\gamma$. Since the label of the boundary of every
cell in $\Psi$ is equal to 1 in $G$, the words corresponding to any
two homotopic loops $\gamma$, $\gamma'$ represent the same element
in $G$. Hence $\phi$ is indeed well-defined. The fact that $\phi$ is
a homomorphism is obvious.

Conversely, the standard argument involving $K(.,1)$-complexes gives
that every injective homomorphism $\phi\colon \pi_1(S)\to G$
corresponds in the above sense to a van Kampen diagram over $G$ on
$S$.

If $G=A(K)$ is a right angled Artin group, then every cell in a van
Kampen diagram is a square, and instead of a van Kampen diagram on
$S$, it is convenient to consider its dual picture: pick a point
inside every cell, connect the points in neighbor cells by an edge
labeled by the label of the common edge of the cells. The result is the so called
{\em $K$-dissection diagram} of the surface, that was introduced
by Crisp and Wiest in \cite{CrispWiest}. The edges of the dual picture having the same labels
form collections of pairwise disjoint simple closed orientation
preserving curves and arcs connecting points on the boundary of $S$.
This is because every cell in the van Kampen diagram has exactly two pairs of opposite edges having the same labels and opposite orientation.
Each of these curves has a natural transversal direction. Each curve
is labeled by a vertex of $K$, two curves intersect only if their
labels are adjacent in $K$.

If $\Delta$ is the $K$-dissection diagram corresponding to a van
Kampen diagram $\Psi$ on $S$, then the corresponding homomorphism
$\phi_v\colon \pi_1(S)\to A(K)$ takes any loop $\gamma$ based at $v$
to the word of labels of the dissection curves and arcs of $\Delta$
crossed by $\gamma$ (a letter in the word can occur with exponent
$1$ or $-1$ according to the direction of the dissection curve
crossed by $\gamma$).

A $K$-dissection diagram $\Delta$ is called faithful if the corresponding homomorphism $\phi$ is faithful. Clearly $\phi$ is faithful only if every connected component of $S\setminus\Delta$ is a disc. The converse statement is far from being true.

There are several partial algorithms allowing
to check whether a homomorphism $\phi$ corresponding to the
$K$-dissection diagram is injective (see Section \ref{embedding}). But the answer to
the next question is still unknown.

Our general technique for showing that non-Abelian surface subgroups do not exist in a right-angled Artin group $A(K)$ is the following: we show that if $A(K)$ contains a non-Abelian surface subgroup, then so does $A(K')$ for some simpler graph $K'$. The graph $K'$ is either a factor in a decomposition of $K$ into an amalgam, or a result of removing certain edges of $K$, or the result of doubling of certain subgraph of $K$. These reduction statements are based on the following simple idea: to find a curve in the kernel of a homomorphism $\phi$ associated with a dissection
diagram $\Delta$ we are allowed to (a) take the boundary of a subsurface spanned by certain dissection curves and (b) take commutators of intersecting curves.

\section{Chordal graphs}
\label{chordal}

In this section, we show that the \RAAG associated to any chordal graph admits no hyperbolic surface subgroups.
A simplicial graph is said to be \emph{chordal} if
every circuit of length greater than 3 admits a ``chord" -- i.e: an
edge of $K$ which is not an edge in the circuit but whose endpoints
both lie in the circuit.

We shall say that a finite simplicial graph $K$ is \emph{treelike} if $K$ can be built by a finite number of glueings along cliques
(complete subgraphs), starting with (a finite number of) cliques. More precisely, the class of treelike graphs is the smallest class of finite
connected graphs which contains all finite cliques and all graphs $K_1\cup_X K_2$ where $X$ is a clique and $K_1$, $K_2$ are treelike.
The following is a standard result in graph theory due to Dirac \cite{Dirac} (see, for example, \cite{BP}).

\begin{lemma}\label{treelikelemma}(Dirac \cite{Dirac})
A finite connected simplicial graph is treelike if and only if it is chordal.
\end{lemma}

\begin{proposition}\label{treelike}
Let $K$ be a chordal graph. Then $A(K)$ contains no on-ended hyperbolic subgroups.
\end{proposition}

\proof We first claim that, since $K$ is chordal, the group $A(K)$ is the fundamental group of a graph of groups in which each vertex group is of the form $G_v=A(K_v)$ for some clique $K_v$ in $K$ and such that, for each edge $e=(u,v)$, the edge group $G_e$ is just $A(K_u\cap K_v)$. This implies  that $A(K)$ acts on a simplicial tree $T$ with vertex stabilizers isomorphic to free Abelian groups, such that
the stabilizer of any edge $e$ of the tree is a retract in the vertex stabilizers of $e_-$ and $e_+$.

To prove the claim we use the fact that, by Lemma \ref{treelikelemma}, $K$ is treelike. Therefore, either $K$ is a clique, in which case the statement holds trivially, or $A(K)$ decomposes nontrivially as an amalgamated product $A(K_1)\star_{A(X)}A(K_2)$ where $X$ is a clique. By induction on the number of vertices in $K$, each of $A(K_1)$ and $A(K_2)$ admits a graph of groups decomposition as claimed. Since it is free abelian $A(X)$ must lie in one of the vertex groups of each of the decompositions of the $A(K_i)$. It follows that $A(K)$ admits a graph of groups decomposition as required.
Note that since there are no HNN-extensions required, the resulting decomposition is always, in fact, a tree of groups.

Suppose now that $A(K)$ contains a one-ended hyperbolic subgroup $G$. Since every Abelian subgroup of $G$ is infinite cyclic,
$G$ acts on the tree $T$ with cyclic vertex and edge stabilizers. Moreover since every edge stabilizer in $A(K)$ is a retract in the corresponding
vertex stabilizers, the stabilizer of any edge $e$ of the tree in $G$ must coincide with stabilizers of both $e_-$ and $e_+$.
That immediately implies that $G$ is cyclic, a contradiction. \endproof

\section{Reduction via the doubling argument}
\label{doubling}

In this section we present our first reduction move, the proof of which is obtained by considering the double of a graph along one of its cliques.
Suppose as before that $K$ is a graph and that $L$ is an induced subgraph. We will be interested in two types of new graphs that can
be built from this.

\begin{enumerate}
\item The {\it double of $K$ along $L$}  is obtained by taking two
identical copies of $K$ and identifying them along $L$. The double is denoted $K*_L K$. It is easy to see that

\[
A(K*_L K)= A(K) *_{A(L)} A(K).
\]

\item The {\it central HNN-extension of $K$ over $L$} is obtained by
taking the graph $K$ adding a single vertex and joining it to all the
vertices of $L$. This extension is denoted $K*_L$. Then  we have

\[
A(K*_L)=A(K)*_{A(L)} = \langle A(K), t \vert txt^{-1} = x , x\in L \rangle
\]
\end{enumerate}

\begin{remark}\label{doubleRemark}
Note that $A(K*_L)$ contains an isomorphic copy of $H(K,L)=A(K*_L K)$ (it is isomorphic to the subgroup of $A(K*_L)$ generated by $A(K)$ and $tA(K)t\iv$ by \cite{LS}).  So, in order to show that $A(K*_L)$ contains a hyperbolic surface subgroup, it is enough to show that $A(K*_L K)$ does so. Note also that $A(K*_L K)$ contains a hyperbolic surface subgroup if and only if $A(K*_L K)*\Z$ does.
This was used recently by Kim \cite{Kim2} to show that $A(P_2(7))< A(P_2(6))$.
\end{remark}

Consider the following 2-complex $D(K,L)$ with fundamental group $H(K,L)$.
Start with two copies $X_K$
and $X_{K'}$ of the square $2$-complex for $A(K)$. These complexes contain isometrically embedded copies of $X_L$. Consider the mapping cylinder of $X_L$, i.e. $X_L\times [0,1]$, and identify $X_L\times\{0\}$ with the copy of $X_L$ in $X_K$ and $X_L\times \{1\}$ with $X_{L'}\subseteq X_{K'}$. The $1$-skeleton of the mapping cylinder consist of the edges in $L, L'$ and the edges connecting $x\in L$ with their copies $x'\in L'$. We shall denote these edges by $t$ with indices. The two-cells are the 2-cells in $X_L, X_{L'}$ plus the squares with two opposite edges $e, e'$ which are copies of each other in $X_L$, $X_{L'}$ and two opposite $t$-edges connecting $e_-$ with $e'_-$ and $e_+$ with $e'_+$. It is easy to deduce from the van Kampen theorem that $\pi_1(D(K,L))$ is isomorphic to $H(K,L)$.

Suppose we have a dissection diagram $\Delta$ for $(S,\partial S)$
associated to a homomorphism $\phi\colon \pi_1(S) \to A(K)$, so that the boundary
components of $S$ have content in $L$. It is the dual picture of a van Kampen diagram $\Psi$ on $S$ over the presentation of $A(K)$. Take a copy $(S',\partial S')$ of $S$ together with a copy $\Psi'$ (over $A(K')$) of the van Kampen diagram $\Psi$.  Let us connect each pair of corresponding boundary components $\gamma, \gamma'$ of $S$ and $S'$ by an annulus $\gamma\times [0,1]$. Let us denote the resulting surface by $D(S)$. Since $\gamma$ is a concatenation of edges of $\Psi$ with labels from $L$, and $\gamma'$ is a similar concatenation of edges with labels from $L'$, we can tessellate each of the annuli by squares corresponding to the cells $t\iv x't'x\iv$ of $D(K,L)$ where $t, t'$ are two $t$-edges. Consider the van Kampen diagram $D(\Psi)$ over $D(K,L)$ on $D(S)$ consisting of the polyhedral decomposition just described and
the cellular map into $D(K,L)$ extending the maps of $\Psi$ and $\Psi'$ and mapping
the $t$-edges of the polyhedral decomposition of the mapping cylinder to the corresponding edges of $D(K,L)$.

\begin{figure}[htbp]
\begin{center}
\includegraphics[scale=1.2]{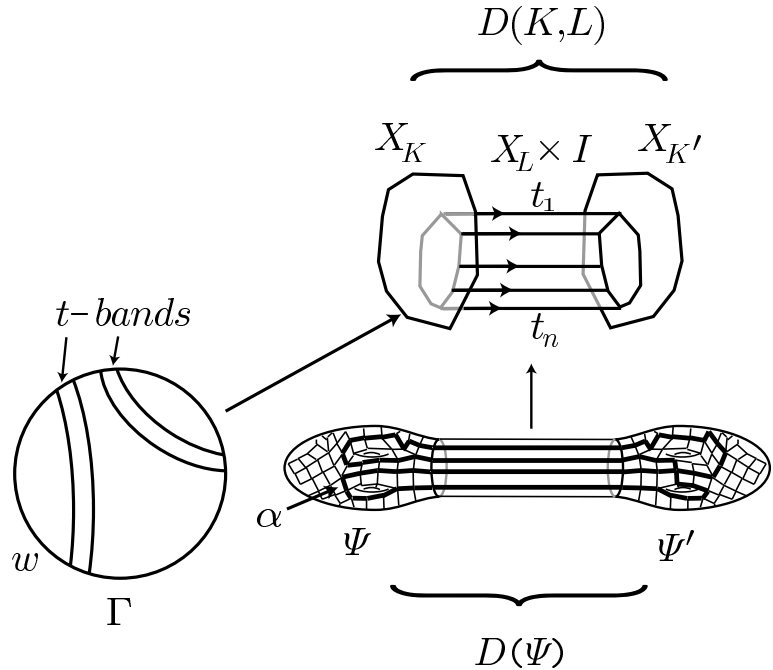}
\end{center}
\caption{Doubling a surface van Kampen diagram.}\label{figexdoubl}
\end{figure}

Our first lemma will address the question of when the diagram $D(\Psi)$ is faithful.

A van Kampen diagram $\Psi$ on $(S, \partial S)$ is called {\it essential} if it is faithful and for every path
$\alpha$ with endpoints on $\partial S$, we have $\rcont(\alpha) \setminus L
\not=\emptyset$ (here the reduced content is the one relative to
the endpoints of $\alpha$).

\begin{theorem}
\label{DoublingLemma}
Let $\Psi$ be an essential van Kampen diagram over $X_K$ on a surface
$(S,\partial S)$  with $\cont(\partial S)\subset L$.
Then $D(\Psi)$ is faithful.
\end{theorem}

\begin{proof}  Suppose that $D(\Psi)$ is not faithful. We
then have a non-trivial polygonal loop $\alpha$ in $D(\Psi)$ that is in the kernel
of the corresponding homomorphism $\phi\colon \pi_1(D(S))\to H(K,L)$.
Let $w$ be the word in the edges of $D(K,L)$ corresponding to $\alpha$.
Since $\alpha$ is $0$-homotopic, there exists a disc van Kampen diagram $\Gamma$ over $D(K,L)$ with boundary label $w$.

Since every 2-cell in the complex $D(K,L)$ involving $t$-letter, has exactly two opposite $t$-edges, we can consider $t$-bands (in another terminology, $t$-corridors) in $\Gamma$. The standard technique one can eliminate $t$-annuli, so we can assume that every (maximal) $t$-band in $\Gamma$ connects two edges on the boundary.

Note that $t$-bands do not intersect. Consider the innermost $t$-band $B$ in $\Gamma$ and the disc subdiagram bounded by a side $l$ of the band $B$ and a part $l'$ of $\partial\Gamma$ that does not contain $t$-edges. Without loss of generality we can assume that the label of the path $l$ is a word in $L$ (and not in $L'$). The path $l'$ corresponds to a subpath $\alpha'$ of $\alpha$. The image of $\alpha'$ does not contain $t$-edges, and so it is inside $X_K$ (the case when it is inside $X_{K'}$ can be easily excluded). Since the image of $l$ in $A(K)$ is inside $A(L)$, and the diagram $\Psi$ is essential, the reduced content of $\alpha'$ in $A(K)$ must be in $L$. Then we can homotop the subpath $\alpha'$ together with the $t$-edge preceding $\alpha'$ and the $t$-edge following $\alpha$ into $S'$ decreasing the number of $t$-edges in $\alpha$. We can conclude the proof by induction on the number of $t$-edges in $\alpha$.
\end{proof}

We now apply the above theorem to prove the following reduction statements.

\begin{corollary}\label{doub1}
Let $K$ be a graph $K=K_1\cup_L K_2$, so that $L$ is a clique.
Suppose that $A(K)$ contains a hyperbolic surface subgroup. Then either $A(K_1 *_L)$ or $A(K_2 *_L)$ also contains a hyperbolic surface subgroup.
\end{corollary}

\begin{proof} Consider a faithful $K$-dissection diagram $\Delta$ and the dual van Kampen diagram $\Psi$ on a surface $S$
associated to a $\pi_1$-injective map $f\colon \pi_1(S)\to A(K)$. Suppose that
there exist two curves $\alpha$ and $\beta$ in $\Delta$ with labels
in $K_1\setminus L$ and so that $\alpha$ and $\beta$ intersect. Using Lemma \ref{BasicLemma1}, we can find a subsurface $(S', \partial S')$ of $S$
with non-Abelian fundamental group and such that $\cont(S')$ is contained in either $K_1$ or $K_2$ and $\cont(\partial S')\subseteq L$. Without loss of generality, we may suppose that $\cont(S')\subseteq K_1$. If $S'=S$ ($\partial S' =\emptyset$) our diagram gives a hyperbolic
surface subgroup of $A(K_1)$ which is itself a subgroup of $A(K_1 *_L)$. We consider the case where $\partial S'$is nonempty.

Since our original map $\phi$ is injective,  so is
the restriction to $\pi_1(S')$. Moreover, by Remark \ref{doubleRemark}, it will suffice to show, by an application of  Theorem \ref{DoublingLemma},
that $A(K_1 *_L K_1)$ contains a hyperbolic surface subgroup. In order to apply Theorem \ref{DoublingLemma} to our situation we need only
show that the $U$-dissection diagram of $S'$ (the restriction of $\Delta$ on $S'$) is essential.

Let $\gamma$ be an essential arc in $S'$. We need to show that $\rcont(\gamma)\not\subset L$.
Suppose that there exists such an arc $\gamma$ with $\rcont[\gamma]
\subset L$. There exist boundary components $\delta$ and $\mu$ such
that $p=\gamma(0)\in\delta$ and $\gamma(1)\in\mu$. (Note that $\delta
$ and $\mu$ may denote the same boundary component.) Now we
consider two loops based at $p$: $\delta$ and $\nu=\gamma\mu\gamma^
{-1}$.  Since $\gamma$ is an essential arc and $S'$ is not an annulus,
$\delta$ and $\nu$ are not homotopic. Hence the subgroup generated by
$\delta$ and $\nu$ in the fundamental group based at the intersection point of these two curves, is free of rank 2. However the reduced content of the commutator $[\delta,\nu]$ is $\emptyset$ by Lemma \ref{comm}. Hence the image of $[\delta,\nu]$ is 1 in $A(K)$, a contradiction.
\end{proof}

Note that the reduction step given by Corollary \ref{doub1} may be useful applied in most cases where the graph $K$ contains
a separating clique $L$. The only exception to this is when $L$ separates just a single vertex off from the rest of $K$. In all other cases,
we have a separation $K=K_1 *_L K_2$ such that the graphs $K_i *_L$ each have fewer vertices than $K$.

Another useful formulation of doubling is the following corollary, whose proof proceeds exactly as above.

\begin{corollary}\label{doub2}
Let $K$ be a graph $K=K_1\cup_L K_2$, so that $L$ is a clique.
Suppose that the \RAAG obtained by doubling
$K_1$ along $L$ does not contain a hyperbolic surface subgroup but $A(K)$ contains a hyperbolic surface subgroup. Then the \RAAG
associated to the graph $K$ with all the edges in $K_1\setminus L$
removed also contains a hyperbolic surface subgroup.
\end{corollary}

\section{The separation lemma}
\label{seplemma}

We first observe that if $K$ is a disconnected graph, then $A(K)$ is
a free product of subgroups $A(K_i)$ where $K_i$ ranges over the
connected components. It follows, by the Kurosh subgroup Theorem
that if $A(K)$ contains a   hyperbolic surface subgroup then this
subgroup is conjugate into one of the free factors $A(K_i)$.

If, on the other hand, $K$ is a nontrivial join of two (or more)
graphs $K_1, K_2$ (i.e. every vertex of $K_1$ is adjacent to every
vertex of $K_2$), then $A(K)$ decomposes as a nontrivial direct
product $A(K_1)\times A(K_2)$.
In this case, any hyperbolic surface subgroup must project
faithfully to at least one factor. (In fact, if $\pi<G_1\times G_2$
is torsion free, then the presence of a nontrivial elements $g_1\in
\ker(\pi\to G_1)=G_2\cap\pi$ and $g_2\in \ker(\pi\to G_2)=G_1\cap
\pi$ imply that $\pi$ contains a subgroup $\Z^2=\<g_1,g_2\>$, a
contradiction when $\pi$ is a hyperbolic surface group.)

Thus Problem \ref{1} reduces easily to the case where $K$ is
connected and not a join of proper subgraphs.

If $K_1,..,K_n$, and $L$ are induced subgraphs  of $K$ such that
$K_i\cap K_j=L$, for all $i,j$,  and $K=\bigcup_{i=,..,n} K_i$,
then we say that $K$ is the result of \emph{gluing the subgraphs
$K_i$ along $L$} and write
\[
K= \bigcup\limits_L  K_i\,= K_1\cup_L...\cup_L K_n.
\]
In the case that every $K_i$ properly contains $L$, we say that $L$ is \emph{separating} or \emph{separates $K$}.

If $V$ is a set of vertices of a graph $K$, then the set $\cup_{v\in V} \Lk(v)\setminus V$ is denoted by $\Lk(V)$.
\emph{Note that this is non-standard usage of the terminology}.
Let $Y=\{Y_1,...,Y_s\}$ be a collection of subsets of $K^0$, $x\in K^0$. By $\Lk_Y(x)$ we
denote the set $\Lk(x)$ union with all $Y_i$ containing $x$. For every $Z\subseteq K^0$, $\Lk_Y(Z)=(\cup_{x\in Z} \Lk_Y(x))\setminus Z$.

The following lemma shows that a separation of $K$ induces a type of separation
of the dissection diagram.

\begin{lemma}[Basic Cutting Lemma]\label{BasicLemma1}
Suppose that $K=K_1\cup_L K_2\cup_L...\cup_L K_n$ and suppose that $(S,\partial S)$ is equipped with a $K$-dissection diagram $\Delta$.
Then in $S$, there exist
collections $B_1, ..., B_n$ of mutually disjoint non-null-homotopic
simple closed curves and arcs (with both ends on $\partial S$) such that
\begin{itemize}
\item[(i)] each $\ga\in B_i$ is isotopic to a composition of subcurves of
$(K_i\setminus L)$-curves;
\item [(ii)] for every $\alpha\in B_i$, $\cont(\alpha)\subset \Lk(K_i\setminus
 L)$;
\item [(iii)] if $S'$ is a connected component of $S\setminus \bigcup B_k$
then $\cont(S') \subseteq \Lk(K_i\setminus L)\cup (K_i\setminus L)$
for some $i$ or $\cont(S')\subseteq L$;
\item [(iv)] if $\Delta$ contains two intersecting $K_i\setminus L$-curves,
then there exists a component of $S\setminus \cup B_k$ with non-Abelian fundamental group and content in $K_i$.
\item[(v)] if a connected component $S'$ of $S\setminus \bigcup B_k$
contains $B_i$-curves and $B_j$-curves for $i\ne j$, then
$\cont(S')\subseteq L$;
\item [(vi)] if $\partial S$ is empty, then one of the connected components is
$S\setminus B$ has non-Abelian fundamental group;
\item [(vii)] every essential curve on $S$ that intersects a curve from $B_i$
also intersects a $K_i\setminus L$-curve from $\Delta$;

\end{itemize}
\end{lemma}

\begin{proof}
Let $\Gamma$ denote the union of all $K\setminus L$-dissection curves and arcs. This is a (not
necessarily connected) graph, where vertices are the intersection
points of the dissection curves and arcs and the end points of the arcs, and edges are parts of the dissection curves and arcs.

Consider the regular neighborhood $N(\Gamma)$. It is a (not
necessarily connected) subsurface of $S$. Let us attach every
component of $S\setminus N(\Gamma)$ that is a null-homotopic (relative to the boundary of $S$) disc to $N(\Gamma)$.
The resulting subsurface is denoted by $S'$. Let $B$ be the collection of all the boundary
components of $S'$. Note that since $K_i\setminus L$ does  not intersect $K_j\setminus L$ for $i\ne j$,
the content of each component of $S\setminus B$ is either in $L$ or in $K_i\setminus L$ for some $i$.
Each boundary component of $S'$ is a non-null-homotopic
simple closed curve whose content is in $\Lk(K_1\setminus
L)\subseteq L$ (this gives (i), (ii), (iii), (v), (vii)).

Property (iv) follows from the assumption that dissection curves are in minimal position with respect to each other and the fact that a surface with Abelian fundamental group cannot have two closed curves in minimal position that intersect.

Property (vi) follows from the fact that one cannot cut a surface without boundary into a collection of annuli by essential simple closed curves.
\end{proof}

\section{Reduction moves}
\label{redmoves}

Here we present several results allowing one to reduce the question
of whether a group $A(K)$ contains a hyperbolic surface subgroup to
the same question for $A(K')$ for simpler $K'$. We have mentioned
two of such statements already: if $A(K)$ is a free or direct
product of $A(K_1)$, $A(K_2)$ (i.e. if either $K$ is a disjoint
union of $K_1$ and $K_2$ or $K_1,K_2$ are full subgraphs of $K$,
$K^0=K_1^0\cup K_2^0$ and every vertex of $K_1$ is adjacent to every
vertex of $K_2$) then $A(K)$ contains a non-abelian surface group if
and only if one of $A(K_i)$ does ($i=1,2$).

We hope that we shall be able to find a \emph{complete} set of
reduction moves in the sense that, if $A(K)$ does not contain a
hyperbolic surface subgroup, then one could use these moves to
reduce $K$ to a 1-vertex graph.

\subsection{Nuclear and dense sets of vertices}
\label{ndsets}

\begin{remark}\label{pr}
A key fact about hyperbolic surface groups which we shall use in all
of our arguments is that the centralizers of non-identity elements
of hyperbolic surface groups are cyclic. In other words, if $\al$ and $\beta$ are
closed curves on the surface $S$, and $\ast$ is an intersection point of $\al$ and $\beta$, then the
elements $\al$ and $\beta$ of $\pi_1(S,\ast)$ commute only if
$\beta$ is homotopic (relative to the basepoint $\ast$) onto $\al$ ---
that is $\beta^k=\al^l$, for some integers $k,l\ne 0$.
\end{remark}

\begin{definition}
We say that two subgraphs $P,Q$ of a graph $K$ are {\em adjacent} if for
all vertices $p\in P$ and $q\in Q$ either $p=q$ or $p$ and $q$ are adjacent.
\end{definition}

\begin{lemma}[Edge reduction move]\label{edgered}
If $(a,b)$ is an edge of $K$ such that $\Lk(a)$ and $\Lk(b)$ are
adjacent subgraphs, then we write $K'$ for the graph obtained from
$K$ by forgetting the edge $(a,b)$ without deleting the vertices.
Then $A(K)$ has a hyperbolic surface subgroup only if $A(K')$ does.
\end{lemma}

\begin{proof}
Suppose that $S$ is a closed hyperbolic surface with a faithful
$K$-dissection diagram $\Delta$. Suppose that somewhere in the
dissection we can find an $a$-curve $\ga_a$ and a $b$-curve $\ga_b$
which intersect in an essential way. Taking any point
$\ast\in\ga_a\cap\ga_b$ as basepoint, we consider the homomorphism
$\phi$ corresponding to $\Delta$. Observe that $\phi(\ga_a)$ and
$\phi(\ga_b)$ commute, because $\phi(\ga_a)\in A(\Lk(a))$,
$\phi(\ga_b)\in A(\Lk(b))$ and the sets $\Lk(a)$ and $\Lk(b)$ are mutually adjacent.
 This is a contradiction, since essentially
intersecting simple closed curves in a hyperbolic surface should
generate a nonabelian free group (see Remark \ref{pr}). Thus, any
faithful dissection diagram is prohibited from admitting
intersections between $a$-curves and $b$-curves, and so induces an
injective map which factors through the right angled Artin group
$A(K')$.
\end{proof}

The above argument illustrates nicely our approach. However, to
obtain further general results it will be convenient to develop some
terminology.

\begin{definition}
Let $K$ be a graph, and $S$ a surface (possibly with boundary)
equipped with a $K$-dissection diagram. A dissection curve (arc)
labeled by $x$ is called an $x$-{\em curve (arc)}.

We define the \emph{content} of $S$ to be the set
\[
\cont(S)=\{ x\in K^0 \,:\, S \text{ has an $x$-curve or
$x$-arc}\,\}\,.
\]
If $S'$ is a subsurface of $S$ which is in general position with
respect to the dissection then $S'$ inherits a dissection (by simply
taking intersections of the dissecting curves and arcs with $S'$).
In this case we may define the content $\cont(S')$ of $S'$
accordingly.

If $\ga$ is a curve in $S$ which is transverse to the
dissection diagram then its content $\cont(\ga)$ is the set of $x\in K^0$ such that
$\ga$ crosses an $x$-curve or $x$-arc.
\end{definition}

\begin{definition}
Let $K$ be a graph, and $(S, \partial S)$ a surface (possibly with boundary)
equipped with a $K$-dissection diagram. If $\ga$ is any
curve in $S$ which is transverse to the dissection and $\ast$ is
a point on $\ga$ then the reduced content $\rcont_\ast(\ga)$ of
$\ga$ relative to $\ast$ is the smallest induced subgraph $L$ of $K^0$ such that
$\phi_\ast(\ga)$ is in $A(L)$. We define the \emph{effective content}
$\effcont(\ga)$ of $\ga$ to be the smallest subset $Z$ of $\cont(\gamma)$ such that

($\diamond$)$Z$ contains $\rcont_\ast(\ga)$ for every $\ast\in \ga$, and the set $\cont(\ga)\setminus Z$ is adjacent to $Z$.
\end{definition}

\begin{remark} Note that the intersection of subsets of $\cont(\gamma)$ satisfying ($\diamond$) also satisfy ($\diamond$). Hence every curve on $S$ has a well-defined effective content.
\end{remark}

\begin{center}
\begin{figure}[htbp]
\includegraphics[scale=1.4]{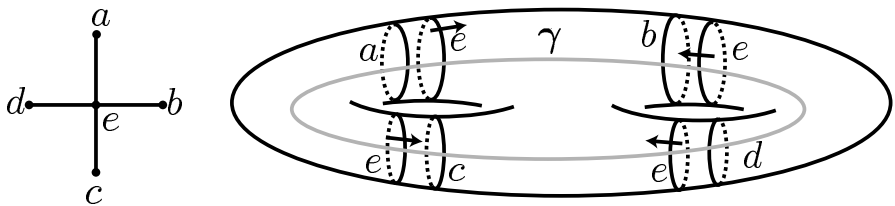}
\caption{The image of the curve $\gamma$ relative to any basepoint on $\gamma$ is inside the subgroup generated by $a,b,c,d$.}\label{figex0}
\end{figure}
\end{center}

\begin{example} On Figure \ref{figex0}, the graph $K$ consists of 5 vertices $a,b,c,d,e$ with $e$ connected to all other vertices and no other edges. The content of the curve $\ga$ is $\{a,b,c,d,e\}$ and the effective content is $\{a,b,c,d\}$.

\end{example}

\begin{remark} Note that the content of a curve or subsurface will always contain
the effective content of that curve or subsurface, but in general
may be strictly larger. Also, note that if $\ga$ is an $x$-curve of
the dissection then $\cont(\ga)\subset \Lk(x)$ where $\Lk(x)$
denotes the set of vertices of the graph $K$ which are adjacent to
the vertex $x$.
\end{remark}

If $V$ is a set of vertices of a graph $K$ then the set $\cap_{v\in
V} \Lk(v)$ is denoted by $C(V)$ (it is the generating set of the
centralizer of $V$ in $A(K)$).

Recall that if $K$ is a  graph we denote $K\opp$ the graph with the same vertex
set but with an edge between two vertices precisely when they are
non-adjacent in $K$.  Note that
if $P_1,...,P_k$ are connected components of $K\opp$, then $K$ is a
join of $P_1,...,P_k$, i.e. $K=P_1\star...\star P_k$.

\begin{definition}[Almost joins]\label{defaprod}
Let $L$ be a subgraph of $K$ and let $L=L_1\star...\star L_n$ be a decomposition of $L$ as a join of subgraphs. Suppose that $$K=K_1\cup_{L} K_2 ...\cup_{L} K_n \ {\rm and }\ \Lk(K_i-L)\subseteq L_i$$ Then we say that $K$ is {\em almost a join} of $K_1,...,K_n$
{\em over} $L_1,...,L_n$ (see Example \ref{exdense} and Figure \ref{exfigure31} below).
\end{definition}

\begin{lemma}[Separating product]\label{product}
Suppose that $K$ is almost a join of subgraphs $K_1,...,K_n$ over $L_1,...,L_n$ ($n\geq 2$).
Let $\Delta$ be a faithful $K$-dissection diagram on a
connected hyperbolic surface $(S,\partial S)$ such that
$\cont(\partial S)$ is in $L$. Then $\cont(S)\subseteq K_i$ for some $i$.
\end{lemma}

\begin{proof} Let us apply Lemma \ref{BasicLemma1} to the decomposition
$K=K_1\cup_L...\cup_L K_n$ to produce a collection of sets
$B_1,...,B_n$ of simple closed non-null-homotopic curves on $S$ such
that the content of each curve from $B_i$ is in $L_i$ (Lemma \ref{BasicLemma1} (ii)), and $B_i$ is
empty if $K_i\setminus L_i$ is empty. If $B_i$ is not empty for
more than one $i$, then one of the connected components $S'$ in
$S\setminus\bigcup B_j$ contains curves $\alpha$ and $\beta$ from
two different sets $B_i, B_j$ as boundary components. These curves
cannot be parallel because their contents (subsets of $L_i$ and
$L_j$ respectively) are disjoint. Hence $\pi_1(S')$ is not Abelian.
On the other hand $\cont(S')\subseteq L$, by Part (vi) of Lemma
\ref{BasicLemma1}. Pick a point $\ast$
on $\alpha$, and consider a curve $\gamma$ that starts at $\ast$,
goes to any point on $\beta$ along some curve $\delta\subset S'$,
then goes around $\beta$ and returns back to $\ast$ along $\delta$.
The image of $\gamma$ under the homomorphism $\phi\colon \pi_1(S)\to
A(K)$ induced by $\Delta$ is a word in $L$: $\phi(\gamma)=ubu\iv$
where $b=\phi(\beta)\in\<L_j\>$.  On the other hand, $c=\phi(\alpha)\in \<L_i\>$.
Since $L$ is a join of $L_1,...,L_n$,
$A(L)=A(L_1)\times...\times A(L_n)$. Therefore $ubu\iv$ is in $A(L_j)$. Hence $c$ and $ubu\iv$ commute.
Therefore
$\phi(\alpha)$ commutes with $\phi(\gamma)$, but $\alpha$ and
$\gamma$ generate a free non-Abelian subgroup in $\pi_1(S')$, a
contradiction (as we assumed that $\Delta$ is faithful). \end{proof}

We shall need the following notation. Let $U, V$ be two subsets of
$K^0$.  The decomposition of $U\opp$  into components corresponds to a canonical decomposition of $U$ as a join $U=U_1\star\ldots \star U_n$ and similarly for $V=V_1\star\ldots\star V_m$. Let $U'$ denote the union of all the $U_i$'s such that every vertex of $U_i$ is adjacent to $V$. We thus obtain a decomposition of $U$ as a join $U'\star U''$ so that the vertices of $U'$are adjacent  $V$. Similarly, we obtain a decomposition of $V=V'\star V''$ so that all the vertices of $V'$ are adjacent to $U$.
We then define $[U,V]=U''\cup V''$. The justification for this notation is the observation that that given any word $w\in A(U)$  and  any word $z\in A(V)$, we have that $[w,z]\in A([U,V])$.


An immediate application of the above observation is the following.

\begin{lemma}\label{comm} Let $\alpha$ and $\beta$ be two closed curves on
a surface $(S,\partial S)$ intersecting at a point $\ast$. Suppose that $(S,\partial S)$ is equipped with a $K$-dissection diagram
$\Delta$. Let $\gamma$ be the commutator curve $[\alpha,\beta]$ with base point $\ast$. Then the content of $\gamma$ is contained in the union of contents of $\alpha$ and $\beta$, and the reduced content of $\gamma$ (relative to $\ast$) is contained in $[\rcont_\ast(\alpha), \rcont_\ast(\beta)]$.
\end{lemma}

This lemma, in turn, immediately implies the following statement that justifies the somewhat unnatural definition of the effective content of a curve given above. This lemma will be used to construct curves in the kernels of homomorphisms $\phi_\ast$.

\begin{lemma}\label{commeff} Let $\Delta$ be a $K$-dissection diagram on $S$, $\alpha$ and $\beta$ be two intersecting closed curves on $S$. Let $\gamma$ be the commutator of $\alpha, \beta$ (based at an intersection point of these curves). Then $$\effcont([\alpha, \beta])\subseteq [\cont(\alpha), \cont(\beta)].$$
\end{lemma}

\begin{definition}[Nuclear subsets] \label{defnuclear} We shall say
that a subset $V$ of $K^0$ is {\em nuclear relative to} $Y\subseteq K^0$ if
there exists an ordering $x_1<...<x_m$ on the set $V$ such that for every $i=1,...,m$ one of the following conditions hold:

\begin{enumerate}
\item $x_i$ is in $C(\{x_1,...,x_{i-1}\})$ and $\Lk_Y(x_i)\setminus \{x_1,...,x_{i-1}\}$ is adjacent to $\{x_1,...,x_{i}\}$;
\item $x_i$ is not in $C(\{x_1,...,x_{i-1}\})$ and $\Lk_Y(x_i)$ is adjacent to $\{x_1,...,x_{i}\}$.
\end{enumerate}
\end{definition}

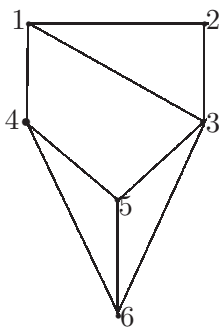
\begin{figure}[htbp]\begin{center}
\unitlength .5mm 
\linethickness{0.4pt}
\ifx\plotpoint\undefined\newsavebox{\plotpoint}\fi 
\begin{picture}(61.25,81.25)(0,0)
\multiput(36.25,33)(-.039968652,.03369906){638}{\line(-1,0){.039968652}}
\multiput(36,33.25)(.036833856,.03369906){638}{\line(1,0){.036833856}}
\put(35.75,33){\line(0,-1){30}}
\multiput(35.75,3)(.0337159254,.0738880918){697}{\line(0,1){.0738880918}}
\multiput(36,3)(-.033719346,.0705040872){734}{\line(0,1){.0705040872}}
\multiput(11.75,54)(.03125,3.28125){8}{\line(0,1){3.28125}}
\put(12,80.25){\line(1,0){48}}
\put(58.75,54.25){\line(0,1){26}}
\put(59,80){\circle*{1.5}}
\put(58.5,54){\circle*{1.5}}
\put(35.75,33.25){\circle*{1.5}}
\put(36,2.5){\circle*{1.5}}
\put(11.5,54){\circle*{2.12}}
\put(12.25,80.25){\circle*{1.12}}
\multiput(12.25,80.25)(.0599870298,-.0337224384){771}{\line(1,0){.0599870298}}
\put(9.5,81.25){\makebox(0,0)[cc]{1}}
\put(61.25,81.25){\makebox(0,0)[cc]{2}}
\put(61.25,53.75){\makebox(0,0)[cc]{3}}
\put(7.75,54.5){\makebox(0,0)[cc]{4}}
\put(38,31.75){\makebox(0,0)[cc]{5}}
\put(38.5,2.25){\makebox(0,0)[cc]{6}}
\end{picture}
\end{center}
\caption{Example of a nuclear subset}\label{figex30}
\end{figure}

\begin{example} In the graph on Figure \ref{figex30}, the set $\{3,4,5\}$ is nuclear relative to $\emptyset$. The ordering is $4<3<5$. The vertex $5$ is in $C(\{3,4\}$ and $\Lk(5)\setminus \{3,4\}=\{6\}$ is adjacent to $3,4,5$; the vertex $4$ is not adjacent to $3$, and $\Lk(4)=\{1,5,6\}$ is adjacent to $3,4$.
\end{example}

The following lemma is obvious.

\begin{lemma}\label{subsets} Let $Y$ be a subset of $K^{(0)}$. Then we have

\begin{enumerate}
\item Any one-vertex subset $\{x\}$ of $K^0$
is nuclear in $K$ relative to any subset of $\Lk(x)$.
\item If $V$ is nuclear in
$K$ relative to $Y$, then every subset of $V$ is nuclear in
$K$ relative to $Y$.
\end{enumerate}
\end{lemma}

\begin{definition}[Characteristic subgraphs] Suppose that
$K$ is an almost join of subgraphs $K_1,...,K_n$ over $L_1,...,L_n$.
For every $i=1,...,n$, if $K_i=L_i$, then we let $P_i=L_i$,
otherwise we let $P_i=K_i\cup \bigcup L_j$. The graphs $P_i$ are
called the {\em characteristic subgraphs of the almost join
decomposition of } $K$.
\end{definition}

\begin{definition}\label{defdense}
We say that a subset $X$ of $K^0$ is {\em dense} in $K$ relative to a subset $Y\subseteq K^0$ if for some decomposition of $K$ as an almost join
and every characteristic subgraph $P$ of that decomposition, $X\cap P$ is
nuclear in $P$ relative to $Y\cap P$.
\end{definition}

\begin{example}\label{exdense}
In the graph on Figure \ref{exfigure31}, the set $\{1,5\}$ is dense relative to
$\emptyset$ but is not nuclear relative to $\emptyset$. Indeed, the graph is an almost join of $K_1=\{1,2,3,4,6\}$ and $K_2=\{2,3,4,5,6\}$ over $L_1=\{2,3\}$ and $L_2=\{4,6\}$. The characteristic subsets are $K_1$ and $K_2$. The intersection of $N$ with each $K_i$ is a one-vertex subset which is dense in $K_i$ relative to $\emptyset$ by Lemma \ref{subsets}.

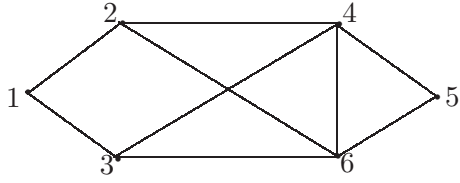
\begin{figure}[htbp]
\begin{center}
\unitlength .6mm 
\linethickness{0.4pt}
\ifx\plotpoint\undefined\newsavebox{\plotpoint}\fi 
\begin{picture}(118,37.25)(0,0)
\put(44,35){\line(1,0){48.5}}
\put(92.5,35){\line(0,-1){29.75}}
\put(92.5,5.25){\line(-1,0){49.5}}
\multiput(43,5.25)(.0565714286,.0337142857){875}{\line(1,0){.0565714286}}
\multiput(44,35.25)(.0549886621,-.0337301587){882}{\line(1,0){.0549886621}}
\multiput(24,19.5)(.044565217,.033695652){460}{\line(1,0){.044565217}}
\multiput(23.5,19.75)(.046511628,-.03372093){430}{\line(1,0){.046511628}}
\multiput(92.5,34.5)(.047109208,-.03372591){467}{\line(1,0){.047109208}}
\multiput(114.5,18.75)(-.055343511,-.033715013){393}{\line(-1,0){.055343511}}
\put(24,19.75){\circle*{1}}
\put(45,35){\circle*{1}}
\put(93,34.75){\circle*{1}}
\put(114.75,18.5){\circle*{1}}
\put(92.75,5.5){\circle*{1}}
\put(44,5){\circle*{1}}
\put(20.75,18.5){\makebox(0,0)[cc]{1}}
\put(42.25,37.25){\makebox(0,0)[cc]{2}}
\put(41.5,3.75){\makebox(0,0)[cc]{3}}
\put(95.25,37.25){\makebox(0,0)[cc]{4}}
\put(118,18.75){\makebox(0,0)[cc]{5}}
\put(94.75,4){\makebox(0,0)[cc]{6}}
\end{picture}

\end{center}
\caption{An example of a dense subset.} \label{exfigure31}
\end{figure}
\end{example}

The following lemma immediately follows from the definition.

\begin{lemma}\label{denrel} Let $X\subseteq K^0$ be a dense subset of $K$ relative to $Y$.
Let $N$ be any subset that is adjacent to $X$. Then $X$ is dense relative to $Y \cup N$.
\end{lemma}

\begin{lemma}\label{power}
Let $(S,\partial S)$ be a surface with non-abelian fundamental group
equipped with a faithful $K$-dissection diagram $\Delta$.
Then the effective content of an essential closed curve
$\gamma$ in $S$ cannot be dense in $K$ relative to $\cont(\partial S)$.
\end{lemma}

\begin{remark}
Since by Lemma \ref{subsets} every single-vertex set $\{x\}$ is nuclear in every graph
$K$ relative to any subset of $\Lk(x)$, Lemma \ref{subsets} implies,
in particular, every closed curve on $S$
whose effective content is a  single generator $x\in K^0$ and for which
$\cont(\partial(S))\subseteq \Lk(x)$ is null-homotopic.
\end{remark}

\proof[Proof of Lemma \ref{power}]
If $K$ is an almost join of $K_1,...,K_n$ over $L_1,...,L_n$ with some $K_i\ne L_i$ then by Lemma \ref{product}, we can
assume that for all but one $i=j$ $K_i=L_i$. Hence $K=P$ for some
characteristic subgraph $P$. If $K=L_1\cup ...\cup L_n$ then $A(K)$ is the direct product of $A(L_i)$. Since the $K$-dissection $\Delta$ is faithful, one of the $L_i$-subdissections $\Delta_i$ is faithful too. So we can assume again that $K=L_i$ for some $i$.

Thus we need to show that there is no non-null-homotopic curves on $S$ with effective content nuclear in $K$ relative $Y = \cont(\partial S)$.

By contradiction, suppose that such a curve $\gamma$ exists. We assume that $\gamma$ is in minimal position with respect to the curves in the dissection diagram.  If $X=\effcont(\gamma)$ is empty, the diagram $\Delta$ is not faithful (since $\gamma$ is non-null-homotopic, but its image
is $1$ in $A(K)$), a contradiction. So we can assume that
$X$ is not empty. Since that set is nuclear in
$K$ relative $Y$, there exists an ordering $x_1<x_2<...<x_m$ of elements of $X$ such that the conditions of the Definition \ref{defnuclear} hold. We can assume that $m=|X|$ is minimal possible for all such
$\gamma$ and that $m>0$.

Suppose first that $x_m\in C(\{x_1,...,x_{m-1}\})$ and $\Lk_Y(x_m)\setminus \{x_1,...,x_{m-1}\}$ is adjacent to $\{x_1,...,x_{m-1}\}$. Consider an $x_m$-curve $\beta$ intersecting $\gamma$ and let $\ast$ be the intersection point. If $\beta$ is an arc connecting two points on the boundary $\partial S$, then let $\delta$ be the closed curve composed of $\beta$ and the connected components of $\partial S$ intersecting $\delta$. This curve is not parallel to any power of $\gamma$ since otherwise the surface $S$ would be an annulus. Note that since $\beta$ intersects $\partial S$,  $x_m\in Y=\cont(\partial S))$ hence $\cont(\delta)\subseteq \Lk(x_m)\cup Y=\Lk_Y(x_m)$. If $\beta$ is a closed curve, then let $\delta=\beta$. Note that in this case also $\delta$ is not parallel to $\gamma$ and $\cont(\delta)\subseteq \Lk_Y(x_m)$. Thus we found a curve $\delta$ with content $\subseteq \Lk_Y(x_m)$ intersecting $\gamma$ at $\ast$.
Then by Lemma \ref{comm}, $$\rcont_\ast([\delta,\gamma])\subseteq [\Lk_Y(x_m), X].$$

Suppose first that $x_m$ is not in $C(\{x_1,...,x_{m-1}\})$. Then by our assumption, $\Lk_Y(x)$ is adjacent to $\{x_1,...,x_m\}$. Therefore $\rcont_\ast([\delta,\gamma])=\emptyset$, so $[\delta,\gamma]=1$ and $\effcont([\delta,\gamma])=\emptyset$, contradicting the minimality of $m$.

Now suppose that $x_m$ is in $C(\{x_1,...,x_{m-1}\}$. Then by our assumption, $\Lk_Y(x)\setminus \{x_1,...,x_m\}$ is adjacent to $\{x_1,...,x_m\}$. Then by Lemma \ref{comm}, $\rcont_\ast([\delta,\gamma])$ is inside $\{x_1,...,x_{m-1}\}$. Now take another point $\diamond$ on $[\delta,\gamma]$. Then the curve $[\delta,\gamma]_\diamond$ based at $\diamond$ is equal to the curve $[\delta,\gamma]_\ast$ based at $\ast$ conjugated by a curve $\varepsilon$ connecting $\ast$ and $\diamond$. The point $\diamond$ belongs to either $\delta$ or $\gamma$, so $\varepsilon$ is either a part of $\delta$ or a part of $\gamma$. If $\varepsilon$ is a part of $\delta$, then its content is in $\Lk_Y(x_m)$ which is a join of $\Lk_Y(x_m)\setminus \{x_1,...,x_{m-1}\}$ and $\Lk_Y(x_m)\cap \{x_1,...,x_{m-1}\}$. Since $\Lk_Y(x_m)\setminus \{x_1,...,x_{m-1}\}$ is adjacent to $\{x_1,...,x_{m-1}\}$, $\rcont_\diamond([\delta,\gamma])$ is also inside $\{x_1,...,x_{m-1}\}$. If $\varepsilon$ is a part of $\gamma$, then $\cont(\varepsilon)$ is inside $\cont(\gamma)$. Hence by the definition of effective content, $\cont(\varepsilon)\setminus X$ is adjacent to $X$. Since $x_m$ is adjacent to $\{x_1,...,x_{m-1}\}$, $\cont(\varepsilon)\setminus \{x_1,...,x_{m-1}\}$ is adjacent to $\{x_1,...,x_{m-1}\}$. Therefore the reduced content of $[\delta,\gamma]$ relative to $\diamond$ is inside $X\setminus \{x_m\}$. Therefore $X\setminus \{x_m\}$ contains the reduced content of $[\delta,\gamma]$ relative to any point of that curve. By Lemma \ref{comm}, the content of $[\delta,\gamma]$ is inside $\Lk_Y(x_m)\cup\cont(\gamma)$. Since $\Lk_Y(x_m)\setminus \{x_1,...,x_{m-1}\}$ is adjacent to $X\setminus \{x_m\}$ and $\cont(\gamma)\setminus\{x_1,...,x_{m}\}$ is adjacent to $X\setminus \{x_m\}$, we have that $\cont([\delta,\gamma])\setminus \{x_1,...,x_{m-1}\}$ is adjacent to $X\setminus\{x_m\}$. Hence the effective content of $[\delta,\gamma]$ is contained in $X\setminus\{x_m\}$. This contradicts the minimality of $m$.
\endproof

Given a graph $K$, in order to show that every $K$-dissection diagram on a surface $K$ is not faithful, one needs (by Lemma \ref{power}) to prove existence of a closed curve with nuclear effective content.

\begin{lemma}\label{ways}
Here are several ways to obtain closed non-null-homotopic curves on a $K$-dissected closed surface $S$:

\begin{enumerate}
\item[(1)] An $x$-curve for $x\in K^0$; its content is in $\Lk(x)$.

\item[(2)] Let $K=K_1\cup_L K_2$ be a non-trivial decomposition of $K$ (so that
$K_i\ne L$). Then Lemma \ref{BasicLemma1} says that there are closed essential curves with content in $\Lk(K_1\setminus L)$. If $K_1\setminus L$-curves in $\Delta$ intersect, then one of these curves bounds a subsurface with content $(K_1\setminus L)\cup \Lk(K_1\setminus L)$.

\item[(3)] If $\alpha$ is a curve as in (1) or (2) with content $L_1$,
$\alpha$ is a boundary component of a subsurface $S'\subseteq S$
with content $L_2$ and non-Abelian fundamental group, then we can consider
a closed essential curve $\beta$ in $S'$ intersecting $\alpha$ and form a commutator $[\alpha,\beta]$ whose effective content is inside $[L_1,L_2]$ by Lemma \ref{commeff}.

\item[(4)] If $\alpha$ and $\beta$ are curves as in (1) or (2) that intersect,
then one can form a commutator of these curves; its effective content is a subset of the commutator of the contents of $\alpha$ and $\beta$ by Lemma \ref{commeff}.
\end{enumerate}

\end{lemma}

This suggests the reduction moves which reduces the question of whether $A(K)$ contains a hyperbolic surface subgroup to simpler graphs: let $L$ be the effective content of a curve $\gamma$ constructed as in (1)-(4) above. If $L$ is dense in $K$ relative to $\emptyset$, and $A(K)$ contains a hyperbolic surface subgroup, then any $K$-dissection diagram does not contain $x$-curves (in Case (1)) or avoids intersections mentioned in the formulation of the corresponding case ((2), (3) or (4)). If the curve we construct is inside a subsurface obtained using Lemma \ref{BasicLemma1} and some decomposition $K=K_1\cup_U...\cup_U K_n$, then instead of the condition that the effective content of $\ga$ is dense in $K$ relative to $\emptyset$ we can assume that the effective content of $\ga$ is dense in the corresponding subgraph of $K$ relative to the content of the boundary of the subsurface.

Here are some concrete reduction moves used in proving Theorem \ref{8vertex}.

\begin{proposition}\label{sepsub}\label{sepsubs} Let $\Delta$ be a faithful $K$-dissection diagram on $S$.

\begin{enumerate}
\item If there exists a non-trivial decomposition $K=K_1\cup_L K_2$ with $L$
dense in $K$, then $\Delta$ does not contain $L$-curves;
\item Suppose that $K=K_1\cup_L K_2=K_1'\cup_{L'} ...\cup_{L'} K_n'$ be
two non-trivial decompositions of $K$. Let $U=\Lk(K_1\setminus L)$, $U_i'=K_i'\setminus L'$ for $i=1,...,n$. Suppose that $[U,L']$ is dense in $K$ relative to $\emptyset$, $U\cap (U_i'\cup \Lk(U_i'))$ is dense in $U_i'\cup \Lk(U_i')$ relative to $\{\Lk(U_i')$ for every $i$, and $U\cap L'$ is dense in $L$ relative to $\Lk(U_1), ..., \Lk(U_n)$. Then $\Delta$ does not contain $U$-curves. \end{enumerate}
\end{proposition}

\begin{figure}[htbp]
\begin{center}
\unitlength .5mm 
\linethickness{0.4pt}
\ifx\plotpoint\undefined\newsavebox{\plotpoint}\fi 
\begin{picture}(62.75,73)(0,0)
\put(14.75,24.25){\framebox(32.25,31)[cc]{}}
\multiput(30.75,72.25)(-.0337078652,-.0362359551){890}{\line(0,-1){.0362359551}}
\multiput(.75,40)(.0337301587,-.0379818594){882}{\line(0,-1){.0379818594}}
\multiput(30.5,6.5)(.0337108954,.0361380798){927}{\line(0,1){.0361380798}}
\multiput(61.75,40)(-.0337108954,.034789644){927}{\line(0,1){.034789644}}
\put(30.75,72.25){\circle*{1.5}} \put(62,39.5){\circle*{1.5}}
\put(30.75,7){\circle*{1.5}} \put(1.5,39.5){\circle*{1.5}}
\put(15.25,55){\circle*{1.5}} \put(48,55){\circle*{1.5}}
\put(47.25,24){\circle*{1.5}} \put(14.75,24){\circle*{1.5}}
\put(50.75,57.25){\makebox(0,0)[cc]{$a$}}
\put(12.5,22){\makebox(0,0)[cc]{$c$}}

\end{picture}
\end{center}
\caption{ $\{a,c\}$ is a separating non-adjacent pair of
vertices.}\label{figex1}
\end{figure}
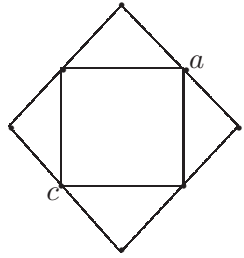

Here is an easy corollary from Proposition \ref{sepsub} which helps ruling out many graphs.

\begin{corollary}[Separating pair]\label{cor1}
Suppose that the connected graph $K$ contains no holes of
length greater than 4 and suppose that $(a,c)$ is a separating pair
of non-adjacent vertices in $K$, that is $K=\bigcup_{\{a,c\}} K_i$
for some collection of proper induced subgraphs $\{ K_1,..,K_n\}$
(as in Fig. \ref{figex1}).

Then $G(K)$ admits a hyperbolic surface subgroup if and only if
$G(K_i)$ does so, for some $i\in\{ 1,..,n\}$.
\end{corollary}

\begin{proof}
Let $X=\{a,c\}$. By a straightforward induction we may easily reduce
to the case of two components: $K=K_1\cup_X K_2$. By Proposition
\ref{sepsub}, part (1), we may suppose that $X$ is a minimal
separating subset, i.e: neither $a$ nor $c$ is a separating vertex
(since every 1-vertex set is nuclear). Let $x_i$ be a vertex in
$K_i\setminus\{a,c\}$, $i=1,2$. Since $K$ is connected, there is a
path in $K$ connecting $x_1$ and $x_2$. All these paths must go
through $a$ or $c$ since $\{ a,c\}$ is a separating pair. If all of
them contain $a$ (resp. $c$) then $a$ (resp. $c$) is a separating
vertex of $K$, which we have assumed is not the case. Therefore one
of these paths contains $a$ but not $c$ and another contains $c$ but
not $a$. This implies that $a$ and $c$ are connected by a path in
$K_1$ as well as by a path in $K_2$.

Note that, since $a$ and $c$ are non-adjacent, any pair of  induced
paths $\ga_1$ from $a$ to $c$ in $K_1$ and $\ga_2$ from $a$ to $c$
in $K_2$ combine to give an induced circuit in $K$. It follows,
since there are no induced circuits of length greater than $4$ in
$K$, and $a,c$ are not adjacent, that $\ga_1$ and $\ga_2$ are both
of length exactly 2. In fact, by this argument, any induced path
from $a$ to $c$ in $K$ is of length 2 and passes through
$Y=\Lk(a)\cap\Lk(c)$. In particular, $Y$ separates $a$ from $c$ in
$K$. Thus $K=K_1'\cup_{Y} K_2'$, where $X\cap K_1'=\{a\}, X\cap
K_2'=\{c\}$. Since one-element subsets are always nuclear, we can
complete the proof by applying Proposition \ref{sepsubs}.
\end{proof}

\subsection{Examples}
\label{examples}

We finish the section with several examples of complicated graphs that
can be completely reduced (to chordal graphs) by using Propositions
\ref{sepsub} and \ref{doub1}. For more examples see Section \ref{8}.


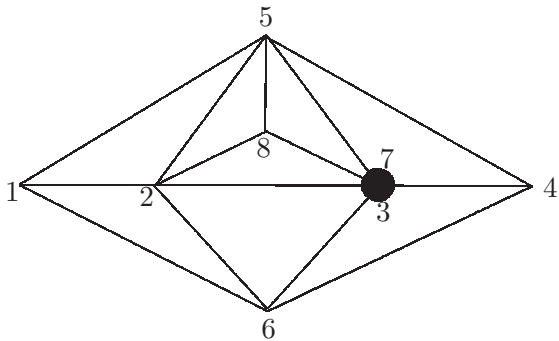
\begin{figure}[htbp]
\begin{center}
\unitlength .6 mm 
\linethickness{0.4pt}
\ifx\plotpoint\undefined\newsavebox{\plotpoint}\fi 
\begin{picture}(124.75,78)(0,0)
\multiput(7.25,40.75)(14.15625,-.03125){8}{\line(1,0){14.15625}}
\multiput(120.5,40.5)(-.0598377282,.0337221095){986}{\line(-1,0){.0598377282}}
\multiput(61.5,73.75)(-.0556690501,-.0337078652){979}{\line(-1,0){.0556690501}}
\multiput(7,40.75)(.0662650602,-.0337349398){830}{\line(1,0){.0662650602}}
\multiput(62,12.75)(.0713851762,.0337181045){823}{\line(1,0){.0713851762}}
\multiput(36.75,40.75)(.0337116155,-.0367156208){749}{\line(0,-1){.0367156208}}
\multiput(62,13.25)(.033719346,.0371253406){734}{\line(0,1){.0371253406}}
\multiput(86.75,40.5)(-.0337381916,.0452091768){741}{\line(0,1){.0452091768}}
\multiput(61.75,74)(-.0337001376,-.045735901){727}{\line(0,-1){.045735901}}
\multiput(37.25,40.75)(.069484241,.033667622){349}{\line(1,0){.069484241}}
\multiput(61.5,52.5)(.03125,2.6875){8}{\line(0,1){2.6875}}
\multiput(61.5,52.75)(.069522472,-.033707865){356}{\line(1,0){.069522472}}
\put(86.25,40.75){\line(0,1){0}}
\put(83.963,38.213){\rule{5.075\unitlength}{5.075\unitlength}}
\multiput(85.101,43.175)(0,-5.723){2}{\rule{2.798\unitlength}{.873\unitlength}}
\multiput(85.759,43.936)(0,-6.683){2}{\rule{1.482\unitlength}{.312\unitlength}}
\multiput(87.128,43.936)(-1.703,0){2}{\multiput(0,0)(0,-6.6){2}{\rule{.446\unitlength}{.228\unitlength}}}
\multiput(87.462,43.936)(-2.201,0){2}{\multiput(0,0)(0,-6.546){2}{\rule{.276\unitlength}{.174\unitlength}}}
\multiput(87.787,43.175)(-3.292,0){2}{\multiput(0,0)(0,-5.399){2}{\rule{.719\unitlength}{.549\unitlength}}}
\multiput(87.787,43.611)(-2.997,0){2}{\multiput(0,0)(0,-6.012){2}{\rule{.424\unitlength}{.289\unitlength}}}
\multiput(87.787,43.788)(-2.843,0){2}{\multiput(0,0)(0,-6.267){2}{\rule{.27\unitlength}{.19\unitlength}}}
\multiput(88.098,43.611)(-3.458,0){2}{\multiput(0,0)(0,-5.927){2}{\rule{.262\unitlength}{.205\unitlength}}}
\multiput(88.393,43.175)(-4.175,0){2}{\multiput(0,0)(0,-5.194){2}{\rule{.389\unitlength}{.344\unitlength}}}
\multiput(88.393,43.406)(-4.04,0){2}{\multiput(0,0)(0,-5.531){2}{\rule{.253\unitlength}{.218\unitlength}}}
\multiput(88.67,43.175)(-4.582,0){2}{\multiput(0,0)(0,-5.081){2}{\rule{.243\unitlength}{.231\unitlength}}}
\multiput(88.67,43.294)(-4.518,0){2}{\multiput(0,0)(0,-5.257){2}{\rule{.178\unitlength}{.17\unitlength}}}
\multiput(88.8,43.175)(-4.776,0){2}{\multiput(0,0)(0,-5.023){2}{\rule{.176\unitlength}{.173\unitlength}}}
\multiput(88.925,39.351)(-5.723,0){2}{\rule{.873\unitlength}{2.798\unitlength}}
\multiput(88.925,42.037)(-5.399,0){2}{\multiput(0,0)(0,-3.292){2}{\rule{.549\unitlength}{.719\unitlength}}}
\multiput(88.925,42.643)(-5.194,0){2}{\multiput(0,0)(0,-4.175){2}{\rule{.344\unitlength}{.389\unitlength}}}
\multiput(88.925,42.92)(-5.081,0){2}{\multiput(0,0)(0,-4.582){2}{\rule{.231\unitlength}{.243\unitlength}}}
\multiput(88.925,43.05)(-5.023,0){2}{\multiput(0,0)(0,-4.776){2}{\rule{.173\unitlength}{.176\unitlength}}}
\multiput(89.044,42.92)(-5.257,0){2}{\multiput(0,0)(0,-4.518){2}{\rule{.17\unitlength}{.178\unitlength}}}
\multiput(89.156,42.643)(-5.531,0){2}{\multiput(0,0)(0,-4.04){2}{\rule{.218\unitlength}{.253\unitlength}}}
\multiput(89.361,42.037)(-6.012,0){2}{\multiput(0,0)(0,-2.997){2}{\rule{.289\unitlength}{.424\unitlength}}}
\multiput(89.361,42.348)(-5.927,0){2}{\multiput(0,0)(0,-3.458){2}{\rule{.205\unitlength}{.262\unitlength}}}
\multiput(89.538,42.037)(-6.267,0){2}{\multiput(0,0)(0,-2.843){2}{\rule{.19\unitlength}{.27\unitlength}}}
\multiput(89.686,40.009)(-6.683,0){2}{\rule{.312\unitlength}{1.482\unitlength}}
\multiput(89.686,41.378)(-6.6,0){2}{\multiput(0,0)(0,-1.703){2}{\rule{.228\unitlength}{.446\unitlength}}}
\multiput(89.686,41.712)(-6.546,0){2}{\multiput(0,0)(0,-2.201){2}{\rule{.174\unitlength}{.276\unitlength}}}
\put(90.01,40.75){\line(0,1){.427}}
\put(89.98,41.18){\line(0,1){.211}}
\put(89.95,41.39){\line(0,1){.209}}
\put(89.91,41.6){\line(0,1){.206}}
\put(89.85,41.8){\line(0,1){.202}}
\put(89.78,42.01){\line(0,1){.197}}
\put(89.69,42.2){\line(0,1){.192}}
\put(89.6,42.39){\line(0,1){.186}}
\put(89.49,42.58){\line(0,1){.179}}
\multiput(89.38,42.76)(-.032,.0429){4}{\line(0,1){.0429}}
\multiput(89.25,42.93)(-.02762,.03269){5}{\line(0,1){.03269}}
\multiput(89.11,43.09)(-.02956,.03095){5}{\line(0,1){.03095}}
\multiput(88.96,43.25)(-.03139,.02909){5}{\line(-1,0){.03139}}
\multiput(88.81,43.4)(-.0331,.02712){5}{\line(-1,0){.0331}}
\multiput(88.64,43.53)(-.0434,.0313){4}{\line(-1,0){.0434}}
\put(88.47,43.66){\line(-1,0){.181}}
\put(88.29,43.77){\line(-1,0){.187}}
\put(88.1,43.87){\line(-1,0){.193}}
\put(87.91,43.97){\line(-1,0){.199}}
\put(87.71,44.05){\line(-1,0){.203}}
\put(87.5,44.11){\line(-1,0){.207}}
\put(87.3,44.17){\line(-1,0){.21}}
\put(87.09,44.21){\line(-1,0){.212}}
\put(86.87,44.24){\line(-1,0){.213}}
\put(86.66,44.26){\line(-1,0){.428}}
\put(86.23,44.25){\line(-1,0){.213}}
\put(86.02,44.23){\line(-1,0){.211}}
\put(85.81,44.19){\line(-1,0){.208}}
\put(85.6,44.14){\line(-1,0){.205}}
\put(85.4,44.08){\line(-1,0){.201}}
\put(85.2,44.01){\line(-1,0){.196}}
\put(85,43.92){\line(-1,0){.191}}
\put(84.81,43.82){\line(-1,0){.184}}
\put(84.62,43.72){\line(-1,0){.177}}
\multiput(84.45,43.6)(-.0424,-.0326){4}{\line(-1,0){.0424}}
\multiput(84.28,43.47)(-.03228,-.0281){5}{\line(-1,0){.03228}}
\multiput(84.12,43.32)(-.0305,-.03002){5}{\line(-1,0){.0305}}
\multiput(83.96,43.17)(-.02862,-.03182){5}{\line(0,-1){.03182}}
\multiput(83.82,43.02)(-.0333,-.0419){4}{\line(0,-1){.0419}}
\multiput(83.69,42.85)(-.0307,-.0438){4}{\line(0,-1){.0438}}
\put(83.56,42.67){\line(0,-1){.182}}
\put(83.45,42.49){\line(0,-1){.189}}
\put(83.35,42.3){\line(0,-1){.195}}
\put(83.26,42.11){\line(0,-1){.2}}
\put(83.19,41.91){\line(0,-1){.204}}
\put(83.12,41.7){\line(0,-1){.208}}
\put(83.07,41.5){\line(0,-1){.21}}
\put(83.03,41.28){\line(0,-1){.212}}
\put(83.01,41.07){\line(0,-1){.854}}
\put(83.03,40.22){\line(0,-1){.21}}
\put(83.07,40.01){\line(0,-1){.208}}
\put(83.12,39.8){\line(0,-1){.204}}
\put(83.19,39.6){\line(0,-1){.2}}
\put(83.26,39.4){\line(0,-1){.195}}
\put(83.35,39.2){\line(0,-1){.189}}
\put(83.45,39.01){\line(0,-1){.183}}
\multiput(83.56,38.83)(.0306,-.0439){4}{\line(0,-1){.0439}}
\multiput(83.69,38.65)(.0332,-.0419){4}{\line(0,-1){.0419}}
\multiput(83.82,38.49)(.02858,-.03185){5}{\line(0,-1){.03185}}
\multiput(83.96,38.33)(.03047,-.03005){5}{\line(1,0){.03047}}
\multiput(84.11,38.18)(.03225,-.02814){5}{\line(1,0){.03225}}
\multiput(84.27,38.04)(.0424,-.0326){4}{\line(1,0){.0424}}
\multiput(84.44,37.91)(.0443,-.03){4}{\line(1,0){.0443}}
\put(84.62,37.79){\line(1,0){.184}}
\put(84.81,37.68){\line(1,0){.19}} \put(85,37.58){\line(1,0){.196}}
\put(85.19,37.49){\line(1,0){.201}}
\put(85.39,37.42){\line(1,0){.205}}
\put(85.6,37.36){\line(1,0){.208}}
\put(85.81,37.31){\line(1,0){.211}}
\put(86.02,37.27){\line(1,0){.213}}
\put(86.23,37.25){\line(1,0){.641}}
\put(86.87,37.26){\line(1,0){.212}}
\put(87.08,37.29){\line(1,0){.21}}
\put(87.29,37.33){\line(1,0){.207}}
\put(87.5,37.39){\line(1,0){.203}}
\put(87.7,37.45){\line(1,0){.199}}
\put(87.9,37.53){\line(1,0){.193}}
\put(88.09,37.62){\line(1,0){.188}}
\put(88.28,37.73){\line(1,0){.181}}
\multiput(88.46,37.84)(.0434,.0313){4}{\line(1,0){.0434}}
\multiput(88.64,37.97)(.03313,.02709){5}{\line(1,0){.03313}}
\multiput(88.8,38.1)(.03142,.02906){5}{\line(1,0){.03142}}
\multiput(88.96,38.25)(.02959,.03092){5}{\line(0,1){.03092}}
\multiput(89.11,38.4)(.02765,.03266){5}{\line(0,1){.03266}}
\multiput(89.25,38.57)(.032,.0429){4}{\line(0,1){.0429}}
\put(89.37,38.74){\line(0,1){.179}}
\put(89.49,38.92){\line(0,1){.186}}
\put(89.6,39.1){\line(0,1){.192}}
\put(89.69,39.29){\line(0,1){.197}}
\put(89.78,39.49){\line(0,1){.202}}
\put(89.85,39.69){\line(0,1){.206}}
\put(89.9,39.9){\line(0,1){.209}}
\put(89.95,40.11){\line(0,1){.211}}
\put(89.98,40.32){\line(0,1){.431}}
\put(5.5,39.25){\makebox(0,0)[cc]{1}}
\put(35.25,38.25){\makebox(0,0)[cc]{2}}
\put(87.75,34.75){\makebox(0,0)[cc]{3}}
\put(88.5,46.5){\makebox(0,0)[cc]{7}}
\put(61.25,49){\makebox(0,0)[cc]{8}}
\put(124.75,40.25){\makebox(0,0)[cc]{4}}
\put(61.75,78){\makebox(0,0)[cc]{5}}
\put(62.25,9){\makebox(0,0)[cc]{6}}
\end{picture}
\end{center}

\caption{An illustration of Proposition \ref{sepsub}. The large dot
represents two non-adjacent vertices $3$ and $7$ having the same
links.}\label{figexk1}
\end{figure}

\begin{example} \label{ex3} Consider the graph $K$ on Figure
\ref{figexk1}. Let $X=\{2,3,5,6,7\}$, $K_1=\{1\}\cup X$,
$K_2=\{8\}\cup X$, $K_3=\{4\}\cup X$. Then $K=K_1\cup_X K_2\cup_X
K_3$. We are going to apply Proposition \ref{sepsub}, Part (2).

Assume first that a faithful dissection diagram $\Delta$ contains $1$-curves. Cutting $S$ along (disjoint) $1$-, $4$-, $8$-curves, we obtain at least one connected component $S_1$ with a non-Abelian fundamental group, content in $X$, and a boundary component having a content in $U=\Lk(1)$.
By Lemma \ref{ways}, part 3, there exists a closed essential curve $\gamma$ with effective content inside $[U, X]=\{5,6\}$. Let $V=C(\{5,6\})=\{1,2,3,4,7\}$. The set $V$ is a separator, $K=K_1'\cup_V K_2'$ where
$K_1'=\{5,8\}\cup V$, $K_2'=\{6\}\cup V$. Let $B$ be the collection of simple closed curves provided by Lemma \ref{BasicLemma1} for this decomposition of $K$. Then every connected component of $S\setminus B$ has content either in $V$ or in $K_1'$ or in $K_2'$. Moreover every curve from $B$ has content in $V$. Since $[U,V]=\emptyset$, the curve $\gamma$ cannot intersect a curve from $B$, so it is inside one of the connected components $S'$ of $S\setminus B$. Therefore $\effcont(\gamma)$ is a subset of $U\cap K_1'$ or $U\cap K_2'$ or $U\cap V$. Each of these sets contains at most one element and is nuclear (in $K$) by Lemma \ref{subsets}, a contradiction with Lemma \ref{power}.

Now suppose that $\Delta$ does not contain $1$-curves. Thus we need to consider
the $7$-vertex graph $K\setminus \{1\}$.

That graph is a join of $\{3,7\}$ and $\{2,4,5,7,8\}$. The group $A(\{3,7\})$ is free. Therefore the subdissection diagram $\Delta'$ of $\Delta$ consisting of $2,4, 5, 7,8$-curves is faithful. Finally note that the subgraph of $K$ spanned by the vertices $2,4,5,7,8$ is isomorphic to a subgraph of the graph on Figure \ref{figex1}.
\end{example}

\begin{figure}[htbp]
\begin{center}

\unitlength .8mm 
\linethickness{0.4pt}
\ifx\plotpoint\undefined\newsavebox{\plotpoint}\fi 
\begin{picture}(119.75,120.5)(0,30)
\multiput(58.75,45)(-.03372591,.037473233){467}{\line(0,1){.037473233}}
\put(43,62.5){\line(0,1){46}}
\multiput(43,108.5)(.038732394,.033702213){497}{\line(1,0){.038732394}}
\multiput(62.25,125.25)(.042553191,-.033687943){423}{\line(1,0){.042553191}}
\multiput(80.25,111)(.03125,-6.0625){8}{\line(0,-1){6.0625}}
\multiput(80.5,62.5)(-.040944123,-.03371869){519}{\line(-1,0){.040944123}}
\multiput(42.75,63)(1.2,-.033333){30}{\line(1,0){1.2}}
\multiput(78.75,62)(.05,.033333){30}{\line(1,0){.05}}
\multiput(80,62.75)(-.03372835,.0419325433){1097}{\line(0,1){.0419325433}}
\multiput(43,108.75)(.5522388,.0335821){67}{\line(1,0){.5522388}}
\multiput(80,111)(-.03372835,-.0435278031){1097}{\line(0,-1){.0435278031}}
\put(40.5,61.75){\makebox(0,0)[cc]{1}}
\put(59,42){\makebox(0,0)[cc]{2}}
\put(83.5,62.75){\makebox(0,0)[cc]{3}}
\put(39.25,108.75){\makebox(0,0)[cc]{6}}
\put(62.75,128.25){\makebox(0,0)[cc]{5}}
\put(83.25,110.75){\makebox(0,0)[cc]{4}}
\put(65.75,86){\makebox(0,0)[cc]{7}}
\multiput(80,111)(.04,.033684211){475}{\line(1,0){.04}}
\multiput(99,127)(.037234043,-.033687943){423}{\line(1,0){.037234043}}
\multiput(114.75,112.75)(-.6586538,-.0336538){52}{\line(-1,0){.6586538}}
\multiput(80.5,111)(.0337320574,-.0459330144){1045}{\line(0,-1){.0459330144}}
\multiput(115.75,63)(-.032609,2.163043){23}{\line(0,1){2.163043}}
\multiput(115,112.75)(-.0337243402,-.0488758553){1023}{\line(0,-1){.0488758553}}
\multiput(80.5,62.75)(2.35,.033333){15}{\line(1,0){2.35}}
\multiput(115.75,63.25)(-.03369906,-.035658307){638}{\line(0,-1){.035658307}}
\multiput(94.25,40.5)(-.03373494,.053614458){415}{\line(0,1){.053614458}}
\put(99.75,129.5){\makebox(0,0)[cc]{8}}
\put(118.75,114.25){\makebox(0,0)[cc]{9}}
\put(119.75,64.25){\makebox(0,0)[cc]{10}}
\put(95.75,37.75){\makebox(0,0)[cc]{11}}
\put(98.25,81){\makebox(0,0)[cc]{12}}
\end{picture}
\end{center}
\caption{Another graph illustrating Proposition
\ref{sepsubs}.}\label{figex2}
\end{figure}
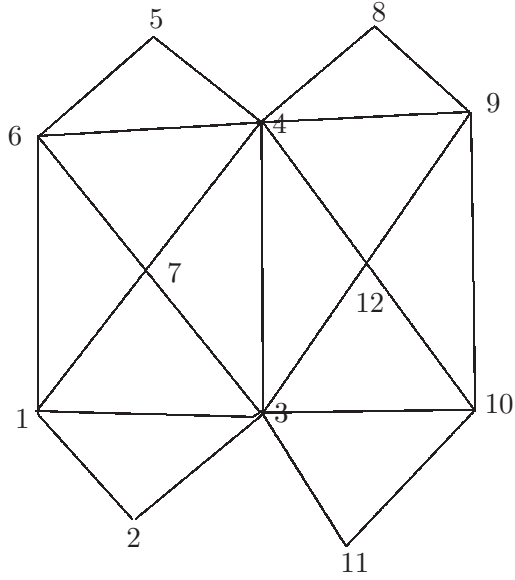

\begin{example} \label{ex2} Consider the graph $K$ on Figure \ref{figex2}.
We shall apply Proposition \ref{sepsubs}, part (2). Take
$X=\{ 6,7,3\}$. Then $K=K_1\cup_X K_2$ where
$K_1=\{1,2\}\cup X$, $K_2=K\setminus \{1,2\}$.
By Lemma \ref{BasicLemma1} for every faithful $K$-dissection diagram containing $K_1\setminus X$-curves and also $K_2\setminus X$-curves, there exists an essential curve $\ga$ with content in $X$.
Let $Y=C(X)=\{1,4,7\}$. Then $K=K_1'\cup_Y K_2'$ with $K_1'=\{5,6\}\cup Y$, $K_2'=(K\setminus K_1')\cup Y=K\setminus \{5,6\}$. Since $[Y,X]=\emptyset$ is dense in $K$ relative to $\emptyset$, we need to consider $X\cup K_1'$ and $X\cup K_2'$. The first set is $\{6,7\}$, the second set is $\{3,7\}$. The set $\{6,7\}$ is nuclear in $K_1'$ relative to $\emptyset$ (the ordering is $6<7$). By Lemma \ref{denrel},
$\{6,7\}$ is dense in $K_1'$ relative to $Y$ since $Y$ is adjacent to $\{6,7\}$. The set $\{3,7\}$ is nuclear in $K_2'$ relative to $\emptyset$ (the ordering is $3<7$), so it is dense in $K_2'$ relative to $Y$. This implies that no faithful $K$-dissection diagram can have $K_1\setminus X$-curves  and $K_2\setminus X$-curves. This allows us to reduce the graph $K$ (by removing some vertices). Continuing in this manner, one can reduce the graph to a 1-vertex graph.
\end{example}

\begin{example}
Graph on Figure \ref{exfigd} is an example illustrating Proposition \ref{doub1}.

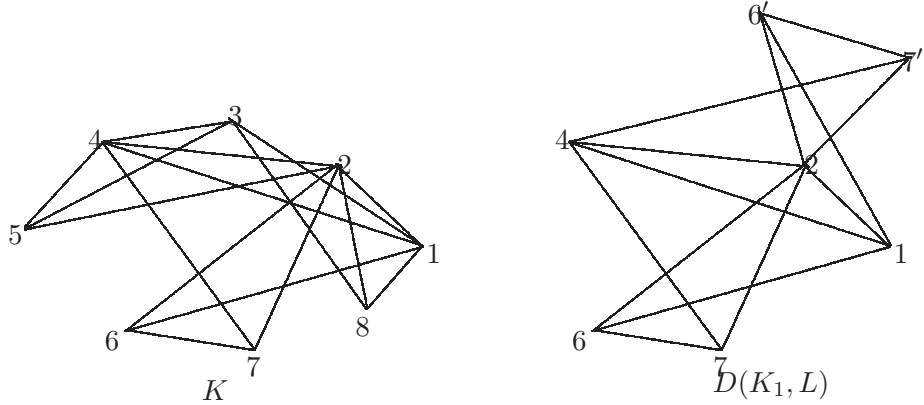
\begin{figure}[ht]
\begin{center}
\unitlength .6mm 
\linethickness{0.4pt}
\ifx\plotpoint\undefined\newsavebox{\plotpoint}\fi 
\begin{picture}(213.75,88.25)(0,0)
\put(107.5,34.25){\makebox(0,0)[cc]{1}}
\put(211,34.25){\makebox(0,0)[cc]{1}}
\put(87.75,54){\makebox(0,0)[cc]{2}}
\put(191.25,54){\makebox(0,0)[cc]{2}}
\put(63.5,64.75){\makebox(0,0)[cc]{3}}
\put(32.5,59.5){\makebox(0,0)[cc]{4}}
\put(136,59.5){\makebox(0,0)[cc]{4}}
\put(14.75,38.5){\makebox(0,0)[cc]{5}}
\put(36.25,15){\makebox(0,0)[cc]{6}}
\put(139.75,15){\makebox(0,0)[cc]{6}}
\put(67.5,9.25){\makebox(0,0)[cc]{7}}
\put(171,9.25){\makebox(0,0)[cc]{7}}
\put(91.75,18){\makebox(0,0)[cc]{8}}
\multiput(105.25,35.75)(-.051183432,.0337278107){845}{\line(-1,0){.051183432}}
\multiput(105,36)(-.1044721408,.0337243402){682}{\line(-1,0){.1044721408}}
\multiput(208.5,36)(-.1044721408,.0337243402){682}{\line(-1,0){.1044721408}}
\multiput(105.25,35.75)(-.119154676,-.033723022){556}{\line(-1,0){.119154676}}
\multiput(208.75,35.75)(-.119154676,-.033723022){556}{\line(-1,0){.119154676}}
\multiput(87.25,53.5)(-.32469512,.03353659){164}{\line(-1,0){.32469512}}
\multiput(190.75,53.5)(-.32469512,.03353659){164}{\line(-1,0){.32469512}}
\multiput(87,53.75)(-.16607565,-.033687943){423}{\line(-1,0){.16607565}}
\multiput(62.75,63.5)(-.21268657,-.03358209){134}{\line(-1,0){.21268657}}
\multiput(62.5,63.5)(-.0646306818,-.0337357955){704}{\line(-1,0){.0646306818}}
\multiput(34.25,59)(-.03371869,-.036608863){519}{\line(0,-1){.036608863}}
\multiput(86.5,54)(-.033723022,-.074190647){556}{\line(0,-1){.074190647}}
\multiput(190,54)(-.033723022,-.074190647){556}{\line(0,-1){.074190647}}
\multiput(39.25,17.25)(.21455224,-.03358209){134}{\line(1,0){.21455224}}
\multiput(142.75,17.25)(.21455224,-.03358209){134}{\line(1,0){.21455224}}
\multiput(86.25,54)(-.0428899083,-.0337155963){1090}{\line(-1,0){.0428899083}}
\multiput(189.75,54)(-.0428899083,-.0337155963){1090}{\line(-1,0){.0428899083}}
\multiput(104.75,35.75)(-.033653846,-.038461538){364}{\line(0,-1){.038461538}}
\put(53.75,4.75){\makebox(0,0)[cc]{}}
\put(157.25,4.75){\makebox(0,0)[cc]{}}
\multiput(104.75,36)(-.035580524,.033707865){534}{\line(-1,0){.035580524}}
\multiput(208.25,36)(-.035580524,.033707865){534}{\line(-1,0){.035580524}}
\multiput(86.5,53.5)(.03360215,-.16935484){186}{\line(0,-1){.16935484}}
\multiput(63,63.5)(.0337142857,-.0474285714){875}{\line(0,-1){.0474285714}}
\multiput(34,59)(.0337162837,-.045954046){1001}{\line(0,-1){.045954046}}
\multiput(137.5,59)(.0337162837,-.045954046){1001}{\line(0,-1){.045954046}}
\put(179.5,87.75){\makebox(0,0)[cc]{$6'$}}
\put(213.75,77.25){\makebox(0,0)[cc]{$7'$}}
\multiput(189.5,54)(-.033670034,.115319865){297}{\line(0,1){.115319865}}
\multiput(189.5,54)(.0337357955,.0340909091){704}{\line(0,1){.0340909091}}
\multiput(137.5,59)(.137067395,.033697632){549}{\line(1,0){.137067395}}
\multiput(208.75,36)(-.0337209302,.0601744186){860}{\line(0,1){.0601744186}}
\multiput(179.75,87.5)(.112794613,-.033670034){297}{\line(1,0){.112794613}}
\put(59,4){\makebox(0,0)[cc]{$K$}}
\put(182.25,5.25){\makebox(0,0)[cc]{$D(K_1,L)$}}
\end{picture}
\end{center}

\caption{An example illustrating Proposition \ref{doub1}.}\label{exfigd}
\end{figure}

Let $K$ be this graph. The clique $L=\{1,2,4\}$ separates the graph $K=K_1\cup_L K_2$ where $K_1=\{1,2,4,6,7\}$, $K_2=\{1,2,4,3,5,8\}$. The graph $D(K,L)$ obtained by doubling $K_1$ along $L$ is a join of $\{2\}$ and the 6-vertex graph spanned by $\{1,4,6,7,6',7'\}$. That graph can be reduced further by applying Corollary  \ref{cor1} because it contains several separating stable pairs of vertices ($\{1,7\}$, $\{6',4\}$, etc.). Thus $A(D(K_1,L))$ does not contain hyperbolic surface subgroups.

The graph $K'$ obtained from $K$ by removing the edge $\{6,7\}$  contains a stable set $X=\{5,6,7,8\}$ such that $K\setminus X$ is a join of $\{2,3\}$, $\{1\}$ and $\{4\}$. Suppose there exists a faithful $K'$-dissection diagram $\Delta$ on a hyperbolic closed surface $S$. Cutting $S$ along $5-,6-, 7-, 8$-curves, we obtain a surface $S'$ one of whose connected components $S_0$ has non-Abelian fundamental group. Taking two intersecting curves $\alpha, \beta$ in $S_0$, we obtain (using Lemma \ref{commeff}) an essential curve $\gamma=[\alpha,\beta]$ with effective content inside $\{2,3\}$ which is a nuclear set in $K$ relative to $\emptyset$ which contradicts Lemma \ref{power}. Thus by Proposition \ref{doub1}, $A(K)$ does not contain hyperbolic surface subgroups.
\end{example}

\section{Embedding results}
\label{embedding}

\subsection{Preliminaries}

We now consider methods to realize embeddings of surface groups into
right-angled Artin groups.

The following definition and lemma are motivated by the fact that
for any graph $K$ and and vertex $x\in K$, the right angled Artin
group $A(K)$ can be considered as an HNN extension of the right
angled subgroup $A(K\setminus \{x\})$ with free letter $a$.

\begin{definition}
Let $x$ denote a vertex of $K$. A subword $x^{\pm 1} u x^{\mp 1}$ of a word $w$
in the alphabet $K^0$ is called an $x$-pinch if $u$ commutes with $x$ in $A(K)$.
\end{definition}

The following lemma is an easy application of standard facts about HNN extensions.

\begin{lemma}\label{pinchlemma}
\begin{itemize}
\item[(1)] If $w$ is a word in $K^0$ representing an identity in $A(K)$, and $w$
contains a letter $x$, then $w$ contains an $x$-pinch.
\item[(2)] Suppose $u$ is a word in $K^0$ and $x\in K^0$ such that $ux=xu$ in
$A(K)$. Then for every letter $y$ in $u$ that is not adjacent to $x$ in $K$, the word $u$ contains a $y$-pinch.
\end{itemize}
\end{lemma}

\subsection{General statement}

In the next section, we shall introduce several finite graphs and
prove that the corresponding right angled Artin groups contain
hyperbolic  surface subgroups. In this Section, we introduce some
notation and a general statement used in the next Section. This
statement automatize proving that certain dissection diagram is
faithful.

Recall that for every graph $K=(V,E)$, $K\opp$ denotes the
complementary graph. Paths in $K\opp$ will be called {\em
anti-paths} in $K$, links in $K\opp$ will be called {\em anti-links}
in K, starts in $K\opp$ are {\em anti-stars} in $K$.

\begin{definition}
Consider the following data that can be assigned to every graph $K$.

\begin{itemize}
\item Linear order $\preceq$ on $V$;
\item Linear order $\preceq_v$ on the anti-star of every vertex $v\in V$.
\end{itemize}
In principle, $\preceq_v$ may not be the restriction of $\preceq$
onto the anti-link of $v$.

The data will be called a {\em load structure} on $K$. We say that a
graph $K$ is {\em loaded} if it is equipped with a load structure.
\end{definition}

Let $p=(v_1,v_2,...,v_k)$ be an anti-path in $K$.

By $\Theta(p)$ we denote the set of vertices of $K$ consisting of

\begin{itemize}
\item $\{v\mid v\succ v_1\}$;
\item for every $i=1,...,k-1$, all vertices $v$ that satisfy $v_{i+1}\preceq_{v_i} v$;
\item all vertices in the anti-star of $v_k$.
\end{itemize}

If we want to specify the graph $K$ and the surface $S$, we shall
write $\Theta(p, K, S)$ instead of $\Theta(p)$.

If $\Delta$ is a $K$-dissection diagram of a surface $S$,  and $Y$
is a set of vertices of $K$, then $S[Y]$ denotes the (possibly
disconnected) surface obtained by cutting $S$ along $v$-curves of
the dissection diagram for all $v\in Y$.

\begin{definition} Let $L$ be a load structure on a graph $K$.
Let $\Delta$ be a dissection diagram of a surface $S$ with boundary.
We say that $\Delta$ is {\em $L$-faithful} if the following
condition holds.

\begin{quote} (*) For every simple anti-path $p=(v_1,...,v_k)$,
every connected component of $S[\Theta(p)]$ is a polygon or a
polygonal annuli whose sides are subarcs of the dissection curves
such that no connected component has two sides labeled by $v_k$
having opposite orientation. In addition, if the component is an
annulus and a boundary component contains a $v_k$-arc, then that
boundary component should be a $v_k$-curve.
\end{quote}
\end{definition}

Let $K$ be a loaded graph. The (finite) set of all simple anti-paths
in $K$ will be denoted by $P_K$. We introduce the {\em lexicographic
order} on $P_K$: given $p=(v_1,...,v_k), p'=(v_1',...,v_m')$ we say
that $p$ is less than $p'$ if either $v_1\prec v_1'$ or $v_1=v_1'$
but $v_2\prec_{v_1} v_2'$ or $v_1=v_1', v_2=v_2'$ but
$v_3\prec_{v_2} v_3'$, etc., or, finally, if $m>k$ and $p$ is the
initial subpath of $p'$. Clearly this defines a linear order on
$P(K)$.

\begin{lemma} \label{cuts} Let $S$ be a surface (possibly with boundary),
$K$ be a loaded graph, with load structure $L$, $\Delta$ a
$K$-dissection diagram of $S$.

If the diagram $\Delta$ is $L$-faithful, then it is faithful.
\end{lemma}

\proof Let, by contradiction, $\gamma$ be a curve in the kernel. Let
$w$ be the word corresponding to $\gamma$. We can assume that
$\gamma$ is chosen in such a way in its homotopy class so that the
word $w$ is the shortest possible. Then $w=1$ in $A(K)$, so for
every letter (vertex) $v_1$ in $w$, $w$ contains a $v_1$-pinch
$w_1=v_1^{\pm 1}i(w_1)v_1^{\mp 1}$ (by Lemma \ref{pinchlemma}). By
the definition of a pinch, the word $i(w_1)$ must represent an
element in $A(K)$ that commutes with $v_1$ in $A(K)$. Hence if any
$v_2$ from the anti-link of $v_1$ occurs in $i(w_1)$, then $w_1$
must contain a $v_2$-pinch $w_2$. The word $i(w_2)$ may contain a
letter $v_3$ from the anti-link of $v_2$. Then $i(w_2)$ must contain
a $v_3$-pinch (by Part (2) or Lemma \ref{pinchlemma}), and so on.
The process stops when the word $i(w_k)$ does not contain vertices
from the anti-link of $v_k$. Since by the definition of a pinch,
$i(w_i)$ does not have occurrences of $v_i$, the anti-path $p=(v_1,
v_2,..., v_k)$ is simple. The finite set of all such anti-paths
corresponding to $w$ is denoted by $P_K(w)$.

Note that every word $w_i$ (more precisely, the occurrence of $w_i$
in $w$) corresponds to a subarc $\gamma[w_i]$ of $\gamma$. Thus we
get a sequence of nested subarcs
$$\gamma(w_k)\subset \gamma(w_{k-1})\subset ...\subset \gamma.$$

Note also that in the set $P_K(w)$, none of the anti-paths is an
initial anti-path of another. That is because if $(v_1,...,v_k)\in
P_K(w)$, then the $v_k$-pinch $w_k$ contains no vertices from the
anti-link of $v_k$.

Let $p=(v_1,...,v_k)$ be the maximal in the lexicographic order
anti-path from $P_K(w)$.

Note that because of the maximality condition, $\gamma[w_{1}]$ is
contained in the subsurface $S[\{v\mid v_1]$, $\gamma[w_2]$ is
contained in $S[\{v\mid v_1\preceq v \hbox{ or } v_{2}\preceq_{v_l}
v \}]$, and so on. Finally, $\gamma[w_k]$ is contained in
$S[\Theta(p)]$ (recall that $\Theta(p)=\{v\mid v_1\preceq v \hbox{
or } v_{2}\preceq_{v_1} v \hbox{ or }... \hbox{ or } v \hbox{ is in
the anti-star of } v_k\}].$ Moreover since $w_k$ is a $v_k$-pinch,
the curve $\gamma[w_k]$ must start and end on subarcs of
$v_k$-curves oriented in the opposite way, and should not be
homotopic to the $v_k$-subarc of the boundary (otherwise $w$ would
not be the shortest word corresponding to curves in the homotopy
class of $\gamma$). But this contradicts the definition of an
$L$-faithful $K$-dissection diagram.
\endproof

Let $\Delta$ be a $K$-dissection diagram on a surface $S$,
$X\subseteq K^0$. For every $X\subseteq K^0$ let $K[X]$ be the graph
induced by $K$ on the complement $K^0\setminus X$. If $K$ is loaded,
then we shall always assume that $K[X]$ inherits the load (i.e. the
orderings on $K[X]$ are restrictions of the orderings on $K$).

The {\em restriction} $\Delta[X]$ of $\Delta$ onto $S[X]$ is the
$K[X]$-dissection diagram on $S[X]$ consisting of the (essential)
intersections of the curves and arcs of $\Delta$ with $S[X]$.

\begin{proposition}\label{thcuts} Let $K$ be a loaded graph with load structure $L$,
$\X$ be a subset of $K^0$.
Let $\Delta$ be a $K$-dissection diagram on a surface $S$. Suppose
that
\begin{itemize}
\item[(1)] For every $x\in \X$, there is a load structure
on $K'=K[X\setminus \Star(x)]$, such that for every anti-path
$p=(x,v_2,...,v_k)$ in $K'$, the set $\Theta(p, K',
S[X\setminus\Star(x)])$ satisfies (*).

\item[(2)] The $K[\X]$-dissection diagram $\Delta[\X]$ on $S[\X]$ is
$L$-faithful.
\end{itemize}
Then the homomorphism $\pi_1(S)\to A(K)$ corresponding to $\Delta$
is faithful.
\end{proposition}

\proof Let $\gamma$ be an essential curve in the kernel. Let $w$ be
the word corresponding to $\gamma$. As before, we can assume that
$\gamma$ is chosen in such a way in its homotopy class so that the
word $w$ is the shortest possible.

Suppose first that $w$ does not contain letters from $\X$. Then
$\gamma$ is in $S[\X]$. Since the $K[\X]$-dissection diagram
$\Delta[\X]$ is faithful by (2), we get a contradiction with Lemma
\ref{cuts}.

Now suppose that $w$ has a letter from $\X$. Let us consider all
$x$-pinches in $w$ for $x\in \X$ and take an innermost pinch
$w'=x^{\pm 1}i(w')x^{\mp 1}$. Then $w'$ does not contain letters
from $\X\setminus\Star(x)$ (otherwise the pinch would not be
innermost). The sub-arc $\gamma'$ corresponding to this pinch is in
$S[\X\setminus\Star(x)]$ with terminal points on $x$-arcs of the
dissection diagram $\Delta[\X\setminus\Star(x)]$ oriented in the
opposite way.

Consider a load structure on $K'=K[\X\setminus\Star(x)]$ for which
the dissection diagram $\Delta[\X\setminus\Star(x)]$ satisfies
condition (1) of the theorem.  Let $p=(v_1,...,v_k)$ be the maximal
(in the lexicographic order) anti-path from $P_{K'}(w')$. Then
$v_1=x$ by the choice of the ordering $\succ$ and, as in the proof
of Lemma \ref{cuts},  we have a sequence of nested subarcs
$$\gamma(w_k)\subseteq ...\subseteq \gamma'$$
where $\gamma(w_i)$ corresponds to a $v_i$-pinch. By the maximality
of $p$, the subarc $\gamma(v_k)$ is in a connected component of the
subsurface $S(\Theta(p, K',S[X\setminus\Star(x)])$.

By (1), the subarc $\gamma(w_i)$ is homotopic to the $v_k$-subarc of
the boundary of the connected component, so $w$ can be shortened, a
contradiction.
\endproof

\begin{corollary} \label{corcuts} Suppose that $X\subseteq K^0$
is a stable set.
Suppose that there exists a load structure on $K[X]$ and a
$K$-dissection diagram $\Delta$ on $S$ such that for every anti-path
$(v_1,...,v_k)$ in $K$ where only $v_1$ can belong to $X$, the set
$\Theta(p)\setminus X$ satisfies (*) for the surface $S[X]$. Then
$\Delta$ is faithful.
\end{corollary}

\proof Let us prove that conditions (1) and (2) of Proposition
\ref{thcuts} hold.

For every $x\in X$, $\Star(X)\cap X=\{x\}$. Hence $K[X\setminus
\Star(x)]=K[X\setminus \{x\}$. Consider the load structure on $K[X]$
induced by the load structure on $K$, and extend it to
$K[X\setminus\{x\}]$ by setting $x\succ v$ for every $v\in K[X]$ and
$y\succ x$ for every $y\ne x\in X$. For every anti-path
$p=(x,v_2,...,v_k)$ then $\Theta=\Theta(p, K, S)$ contains $X$, so
the connected components of $S[\Theta]$ are the same as the
connected components of $S[X]$ cut by the $x$-curves with $x\in
K[X]$. Thus $\Theta$ satisfies (*) by the assumption of the
corollary. This gives (1).

Condition (2) follows directly from the conditions of the corollary.
\endproof

The significance of Corollary \ref{corcuts} is that it allows us to
deal with a subsurface $S[X]$ of $S$ which is in many cases much
simpler than $S$.

\subsection{Proofs for $n$-gones ($n\ge 5$), $P_1(6) - P_4(8)$}
\label{prev}

In all cases considered in this section, the surface $S$ is obtained
as a double of a planar surface $S_0$ obtained by identifying the
respected boundary components of $S_0$ and its copy $S_0'$. The
dissection diagram in each case is defined in $S_0$ and we consider
(almost) a copy of the dissection diagram on $S_0'$. As a result of
identification, the dissection arcs in $S_0$ become dissection
(closed) curves in $S=S_0\cup S_0'$. We shall use the following
convention of choosing the transverse directions on the dissection
curves. The directions on the arcs and closed curves in $S_0$ are
given on the pictures of dissection diagrams if needed. If the
direction is not given, it can be chosen arbitrarily. The directions
on the arcs in $S_0'$ are naturally  determined by the directions on
the corresponding arcs in $S_0'$. But the directions of the closed
dissection curves in $S_0'$ {\em are always chosen opposite to the
directions of the corresponding curves in $S_0$.}

We are going to apply Corollary \ref{corcuts}. Thus in each case, we
specify the set $X$ and the load structure on the graph.

\subsubsection{$n$-gones} The following two dissection diagrams serve
$n$-gons, $n\ge 5$. The first one is for even $n$, the second one -
for odd $n$.

\unitlength .5mm \linethickness{0.4pt}
\ifx\plotpoint\undefined\newsavebox{\plotpoint}\fi 
\begin{picture}(252,139.25)(0,0)
\put(142.5,109.75){\oval(219,50.5)[]}
\put(142.25,35.25){\oval(219,50.5)[]}
\put(55.75,110.5){\circle{23.71}} \put(50.75,36.75){\circle{23.71}}
\put(233,36.75){\circle{23.71}} \put(95.75,110.25){\circle{23.71}}
\put(84.5,36.5){\circle{23.71}} \put(199.25,36.5){\circle{23.71}}
\put(121,36.5){\circle{23.71}} \put(162.75,36.5){\circle{23.71}}
\put(227.5,110.75){\circle{23.71}}
\put(33,110.5){\line(1,0){10.75}}
\put(67.5,110.5){\line(1,0){16.5}}
\multiput(107.75,110.75)(3.5,-.0625){4}{\line(1,0){3.5}}
\put(203,111){\line(1,0){12.5}}
\multiput(239.5,110.5)(3.125,.0625){4}{\line(1,0){3.125}}
\put(94.25,139.25){\makebox(0,0)[cc]{1}}
\put(94,64.75){\makebox(0,0)[cc]{1}}
\put(37.75,115){\makebox(0,0)[cc]{2}}
\put(55.5,110){\makebox(0,0)[cc]{3}}
\put(50.5,36.25){\makebox(0,0)[cc]{3}}
\put(233.25,36.25){\makebox(0,0)[]{3}}
\put(75.5,113.5){\makebox(0,0)[cc]{4}}
\put(95.75,110.25){\makebox(0,0)[cc]{5}}
\put(84.5,36.5){\makebox(0,0)[cc]{5}}
\put(199.25,36.5){\makebox(0,0)[]{5}}
\put(113.75,114.25){\makebox(0,0)[cc]{6}}
\put(208.75,114.75){\makebox(0,0)[cc]{$_{n-2}$}}
\put(226.75,109.75){\makebox(0,0)[cc]{$_{n-1}$}}
\put(244.75,114){\makebox(0,0)[cc]{$_{n}$}}
\put(154.5,110.25){\makebox(0,0)[cc]{$...$}}
\put(140.25,76){\makebox(0,0)[cc]{$n$ is even.}}
\put(140,1.5){\makebox(0,0)[cc]{$n$ is odd.}}
\put(32.5,36.75){\line(1,0){6.25}}
\put(251.25,36.75){\line(-1,0){6.25}}
\put(62.75,37){\line(1,0){10}}
\put(221,37){\line(-1,0){10}}
\put(35.5,41){\makebox(0,0)[cc]{2}}
\put(248.25,41){\makebox(0,0)[]{2}}
\put(68,41){\makebox(0,0)[cc]{4}}
\put(215.75,41){\makebox(0,0)[]{4}}
\put(120.75,36.75){\makebox(0,0)[cc]{$_{n-2}$}}
\put(163,36.75){\makebox(0,0)[]{$_{n-2}$}}
\put(132.75,36.75){\line(1,0){18}}
\put(142.75,60.5){\line(0,-1){50.25}}
\put(102.25,37.25){\makebox(0,0)[cc]{$...$}}
\put(180.5,36.5){\makebox(0,0)[cc]{$...$}}
\put(146.25,54){\makebox(0,0)[cc]{$_{n}$}}
\put(137,40.25){\makebox(0,0)[cc]{$_{n-1}$}}
\end{picture}

\proof The set $X$ is empty. The order on the vertices is such that
$1\succ i$ for every $i\ne 1$, and the orders on the anti-links are
arbitrary. Let $p=(v_1,...,v_k)$ be any anti-path. For every $i$
then $\Theta(p)$ contains $\{1,2,...,n\}\setminus
\{i-1,i+1\}\cup\{1\}$ and satisfies (*).
\endproof

\subsubsection{The $6$-vertex graphs $P_1(6), P_2(6)$}
\label{p16p26}

\begin{center}
\unitlength .6mm 
\linethickness{0.4pt}
\ifx\plotpoint\undefined\newsavebox{\plotpoint}\fi 
\begin{picture}(255.75,69.25)(0,0)
\put(13.75,23.75){\line(1,0){18.25}}
\put(89,23.75){\line(1,0){18.25}}
\put(32,23.75){\line(1,0){18}}
\put(107.25,23.75){\line(1,0){18}}
\put(50,23.75){\line(1,0){17.75}}
\put(125.25,23.75){\line(1,0){17.75}}
\multiput(67.75,23.75)(-.05141844,.033687943){564}{\line(-1,0){.05141844}}
\multiput(143,23.75)(-.05141844,.033687943){564}{\line(-1,0){.05141844}}
\multiput(38.75,42.75)(-.044964029,-.033723022){556}{\line(-1,0){.044964029}}
\multiput(114,42.75)(-.044964029,-.033723022){556}{\line(-1,0){.044964029}}
\multiput(13.75,24)(.04417122,-.033697632){549}{\line(1,0){.04417122}}
\multiput(89,24)(.04417122,-.033697632){549}{\line(1,0){.04417122}}
\multiput(38,5.5)(.054644809,.033697632){549}{\line(1,0){.054644809}}
\multiput(113.25,5.5)(.054644809,.033697632){549}{\line(1,0){.054644809}}
\multiput(114.25,42.25)(.033682635,-.054640719){334}{\line(0,-1){.054640719}}
\put(11.25,22.5){\makebox(0,0)[cc]{1}}
\put(86.5,22.5){\makebox(0,0)[cc]{1}}
\put(30.5,26.25){\makebox(0,0)[cc]{2}}
\put(105.75,26.25){\makebox(0,0)[cc]{2}}
\put(52,26){\makebox(0,0)[cc]{3}}
\put(127.25,26){\makebox(0,0)[cc]{3}}
\put(69.75,24){\makebox(0,0)[cc]{4}}
\put(145,24){\makebox(0,0)[cc]{4}}
\put(38.5,46){\makebox(0,0)[cc]{5}}
\put(113.75,46){\makebox(0,0)[cc]{5}}
\put(40.5,3.5){\makebox(0,0)[cc]{6}}
\put(115.75,3.5){\makebox(0,0)[cc]{6}}
\put(212.38,36.38){\oval(86.75,58.75)[]}
\multiput(212,65.75)(-.03125,-7.34375){8}{\line(0,-1){7.34375}}
\put(188.5,35){\circle{11.24}} \put(234.75,35){\circle{11.24}}
\put(169,35.5){\line(1,0){13.5}}
\put(193.75,35.25){\line(1,0){35}}
\put(240,35){\line(1,0){15.75}}
\multiput(185.25,30.75)(-.040523691,-.033665835){401}{\line(-1,0){.040523691}}
\multiput(238,30.5)(.045984456,-.033678756){386}{\line(1,0){.045984456}}
\put(233.5,69.25){\makebox(0,0)[cc]{2}}
\put(214,55.75){\makebox(0,0)[cc]{1}}
\put(175.25,39){\makebox(0,0)[cc]{3}}
\put(188.5,35.25){\makebox(0,0)[cc]{4}}
\put(200.75,38.5){\makebox(0,0)[cc]{6}}
\put(234.75,35.5){\makebox(0,0)[cc]{4}}
\put(178,22.25){\makebox(0,0)[cc]{5}}
\put(244.25,23){\makebox(0,0)[cc]{5}}
\put(248.25,38.5){\makebox(0,0)[cc]{3}}
\multiput(106,24)(.03361345,.07983193){238}{\line(0,1){.07983193}}
\put(20.25,36.75){\makebox(0,0)[cc]{$P_1(6)$}}
\put(100.5,37.75){\makebox(0,0)[cc]{$P_2(6)$}}
\multiput(113.5,5.5)(.033667622,.052292264){349}{\line(0,1){.052292264}}
\multiput(38.75,42.75)(-.033653846,-.071153846){260}{\line(0,-1){.071153846}}
\multiput(49.25,24)(-.033639144,-.055045872){327}{\line(0,-1){.055045872}}
\end{picture}
\end{center}

\proof Let the set of pairwise non-adjacent vertices be $\{2,4\}$.
Orders on the graph and on the anti-links are arbitrary.

We need to show that the conditions of  Corollary \ref{cuts} hold.
Let $p=(v_1,...,v_k)$ be any anti-path in $P=P_1(6), P_2(6)$. Let
$\Theta=\Theta(p, K[X], S[X])$.

Suppose that $\Theta$ does not satisfy (*) for $S[X]=S_0$.

Consider six different possibilities for $v_k$.

{\bf 1}. Let $v_k=1$. Then $\Theta\supseteq \{1,3\}$ and satisfies
(*).

{\bf 2.} Let $v_k=2$ (and then $k=1$). Then $\Theta\supseteq \{6\}$
and satisfies (*).

{\bf 3.} Let $v_k=3$. Then $\Theta\supseteq \{1,3\}$ and satisfies
(*).

{\bf 4.} Let $v_k=4$ (and then $k=1$). Then $\Theta\supseteq \{1\}$
and satisfies (*).

{\bf 5.} Let $v_k=5$. Then $\Theta\supseteq \{5,6\}$ and satisfies
(*).

{\bf 6.} Let $v_k=6$. Then $\Theta\supseteq \{5,6\}$ and satisfies
(*).

Thus in all cases $\Theta$ satisfies (*), a contradiction.
\endproof

\subsubsection{The 7-vertex graph $P_1(7)$}

\begin{center}
\unitlength .5mm 
\linethickness{0.4pt}
\ifx\plotpoint\undefined\newsavebox{\plotpoint}\fi 
\begin{picture}(196.78,120.5)(0,0)
\multiput(12,25.5)(.067467652,.033733826){541}{\line(1,0){.067467652}}
\multiput(48.5,43.75)(.065543071,-.033707865){534}{\line(1,0){.065543071}}
\multiput(83.5,25.75)(-.060153584,-.033703072){586}{\line(-1,0){.060153584}}
\multiput(48.25,6)(-.062716263,.033737024){578}{\line(-1,0){.062716263}}
\put(48.5,43.75){\line(0,-1){37.5}}
\multiput(12,25.5)(8.90625,.03125){8}{\line(1,0){8.90625}}
\multiput(48.5,43.5)(-.034155598,-.033681214){527}{\line(-1,0){.034155598}}
\multiput(64.75,25.75)(-.033684211,-.040526316){475}{\line(0,-1){.040526316}}
\put(9.75,22.75){\makebox(0,0)[cc]{1}}
\put(30.25,23){\makebox(0,0)[cc]{2}}
\put(45.5,23.25){\makebox(0,0)[cc]{3}}
\put(60,23.25){\makebox(0,0)[cc]{4}}
\put(86,23.75){\makebox(0,0)[cc]{5}}
\put(49,47.5){\makebox(0,0)[cc]{6}}
\put(48.5,3){\makebox(0,0)[cc]{7}}
\put(196.78,60.5){\line(0,1){1.674}}
\put(196.75,62.17){\line(0,1){1.673}}
\multiput(196.67,63.85)(-.032,.4174){4}{\line(0,1){.4174}}
\multiput(196.55,65.52)(-.0298,.27746){6}{\line(0,1){.27746}}
\multiput(196.37,67.18)(-.0328,.23693){7}{\line(0,1){.23693}}
\multiput(196.14,68.84)(-.03113,.18341){9}{\line(0,1){.18341}}
\multiput(195.86,70.49)(-.03305,.16414){10}{\line(0,1){.16414}}
\multiput(195.53,72.13)(-.03171,.13587){12}{\line(0,1){.13587}}
\multiput(195.15,73.76)(-.033088,.12447){13}{\line(0,1){.12447}}
\multiput(194.72,75.38)(-.03196,.106947){15}{\line(0,1){.106947}}
\multiput(194.24,76.98)(-.033012,.0993){16}{\line(0,1){.0993}}
\multiput(193.71,78.57)(-.032028,.087329){18}{\line(0,1){.087329}}
\multiput(193.13,80.14)(-.032857,.081767){19}{\line(0,1){.081767}}
\multiput(192.51,81.7)(-.033574,.076688){20}{\line(0,1){.076688}}
\multiput(191.84,83.23)(-.032638,.068751){22}{\line(0,1){.068751}}
\multiput(191.12,84.74)(-.033214,.064777){23}{\line(0,1){.064777}}
\multiput(190.36,86.23)(-.033712,.061076){24}{\line(0,1){.061076}}
\multiput(189.55,87.7)(-.032828,.055401){26}{\line(0,1){.055401}}
\multiput(188.69,89.14)(-.033228,.052358){27}{\line(0,1){.052358}}
\multiput(187.8,90.55)(-.033569,.049485){28}{\line(0,1){.049485}}
\multiput(186.86,91.94)(-.032728,.045207){30}{\line(0,1){.045207}}
\multiput(185.87,93.3)(-.032994,.04276){31}{\line(0,1){.04276}}
\multiput(184.85,94.62)(-.033214,.040428){32}{\line(0,1){.040428}}
\multiput(183.79,95.92)(-.033391,.0382){33}{\line(0,1){.0382}}
\multiput(182.69,97.18)(-.033527,.0360686){34}{\line(0,1){.0360686}}
\multiput(181.55,98.4)(-.0336247,.0340263){35}{\line(0,1){.0340263}}
\multiput(180.37,99.59)(-.0346489,.0329827){35}{\line(-1,0){.0346489}}
\multiput(179.16,100.75)(-.036689,.0328469){34}{\line(-1,0){.036689}}
\multiput(177.91,101.86)(-.04003,.033692){32}{\line(-1,0){.04003}}
\multiput(176.63,102.94)(-.042365,.0335){31}{\line(-1,0){.042365}}
\multiput(175.31,103.98)(-.044815,.033262){30}{\line(-1,0){.044815}}
\multiput(173.97,104.98)(-.04739,.032976){29}{\line(-1,0){.04739}}
\multiput(172.6,105.94)(-.050104,.032638){28}{\line(-1,0){.050104}}
\multiput(171.19,106.85)(-.055007,.033483){26}{\line(-1,0){.055007}}
\multiput(169.76,107.72)(-.058245,.033058){25}{\line(-1,0){.058245}}
\multiput(168.31,108.55)(-.061696,.032565){24}{\line(-1,0){.061696}}
\multiput(166.83,109.33)(-.068359,.033452){22}{\line(-1,0){.068359}}
\multiput(165.32,110.06)(-.072652,.03284){21}{\line(-1,0){.072652}}
\multiput(163.8,110.75)(-.077302,.032134){20}{\line(-1,0){.077302}}
\multiput(162.25,111.4)(-.086943,.033063){18}{\line(-1,0){.086943}}
\multiput(160.69,111.99)(-.093084,.032178){17}{\line(-1,0){.093084}}
\multiput(159.1,112.54)(-.10656,.033227){15}{\line(-1,0){.10656}}
\multiput(157.5,113.04)(-.115206,.032094){14}{\line(-1,0){.115206}}
\multiput(155.89,113.49)(-.13549,.03332){12}{\line(-1,0){.13549}}
\multiput(154.27,113.89)(-.14885,.03181){11}{\line(-1,0){.14885}}
\multiput(152.63,114.24)(-.18303,.0333){9}{\line(-1,0){.18303}}
\multiput(150.98,114.54)(-.20696,.03116){8}{\line(-1,0){.20696}}
\multiput(149.33,114.78)(-.27708,.03309){6}{\line(-1,0){.27708}}
\multiput(147.66,114.98)(-.33356,.02953){5}{\line(-1,0){.33356}}
\put(146,115.13){\line(-1,0){1.672}}
\put(144.32,115.23){\line(-1,0){1.674}}
\put(142.65,115.27){\line(-1,0){1.674}}
\put(140.98,115.27){\line(-1,0){1.673}}
\put(139.3,115.21){\line(-1,0){1.671}}
\multiput(137.63,115.1)(-.33335,-.0318){5}{\line(-1,0){.33335}}
\multiput(135.97,114.94)(-.2373,-.02998){7}{\line(-1,0){.2373}}
\multiput(134.3,114.73)(-.20674,-.03257){8}{\line(-1,0){.20674}}
\multiput(132.65,114.47)(-.16452,-.0311){10}{\line(-1,0){.16452}}
\multiput(131,114.16)(-.14863,-.03283){11}{\line(-1,0){.14863}}
\multiput(129.37,113.8)(-.124854,-.031608){13}{\line(-1,0){.124854}}
\multiput(127.75,113.39)(-.114985,-.03288){14}{\line(-1,0){.114985}}
\multiput(126.14,112.93)(-.099685,-.031831){16}{\line(-1,0){.099685}}
\multiput(124.54,112.42)(-.092862,-.032812){17}{\line(-1,0){.092862}}
\multiput(122.96,111.86)(-.086715,-.033655){18}{\line(-1,0){.086715}}
\multiput(121.4,111.26)(-.077081,-.032661){20}{\line(-1,0){.077081}}
\multiput(119.86,110.6)(-.072426,-.033334){21}{\line(-1,0){.072426}}
\multiput(118.34,109.9)(-.065167,-.032443){23}{\line(-1,0){.065167}}
\multiput(116.84,109.16)(-.061472,-.032985){24}{\line(-1,0){.061472}}
\multiput(115.37,108.37)(-.058018,-.033454){25}{\line(-1,0){.058018}}
\multiput(113.92,107.53)(-.052748,-.032604){27}{\line(-1,0){.052748}}
\multiput(112.49,106.65)(-.04988,-.032979){28}{\line(-1,0){.04988}}
\multiput(111.09,105.73)(-.047164,-.033299){29}{\line(-1,0){.047164}}
\multiput(109.73,104.76)(-.044587,-.033567){30}{\line(-1,0){.044587}}
\multiput(108.39,103.75)(-.040819,-.032732){32}{\line(-1,0){.040819}}
\multiput(107.08,102.71)(-.038594,-.032935){33}{\line(-1,0){.038594}}
\multiput(105.81,101.62)(-.0364641,-.0330964){34}{\line(-1,0){.0364641}}
\multiput(104.57,100.49)(-.0344231,-.0332183){35}{\line(-1,0){.0344231}}
\multiput(103.36,99.33)(-.0333918,-.0342549){35}{\line(0,-1){.0342549}}
\multiput(102.2,98.13)(-.0332801,-.0362964){34}{\line(0,-1){.0362964}}
\multiput(101.06,96.9)(-.03313,-.038427){33}{\line(0,-1){.038427}}
\multiput(99.97,95.63)(-.032938,-.040653){32}{\line(0,-1){.040653}}
\multiput(98.92,94.33)(-.032702,-.042984){31}{\line(0,-1){.042984}}
\multiput(97.9,93)(-.033537,-.046996){29}{\line(0,-1){.046996}}
\multiput(96.93,91.63)(-.033231,-.049713){28}{\line(0,-1){.049713}}
\multiput(96,90.24)(-.03287,-.052583){27}{\line(0,-1){.052583}}
\multiput(95.11,88.82)(-.032449,-.055623){26}{\line(0,-1){.055623}}
\multiput(94.27,87.38)(-.033295,-.061305){24}{\line(0,-1){.061305}}
\multiput(93.47,85.9)(-.032771,-.065002){23}{\line(0,-1){.065002}}
\multiput(92.72,84.41)(-.0337,-.072257){21}{\line(0,-1){.072257}}
\multiput(92.01,82.89)(-.03305,-.076915){20}{\line(0,-1){.076915}}
\multiput(91.35,81.35)(-.032298,-.081989){19}{\line(0,-1){.081989}}
\multiput(90.73,79.8)(-.033281,-.092695){17}{\line(0,-1){.092695}}
\multiput(90.17,78.22)(-.032334,-.099523){16}{\line(0,-1){.099523}}
\multiput(89.65,76.63)(-.03346,-.114817){14}{\line(0,-1){.114817}}
\multiput(89.18,75.02)(-.032238,-.124693){13}{\line(0,-1){.124693}}
\multiput(88.76,73.4)(-.03358,-.14846){11}{\line(0,-1){.14846}}
\multiput(88.39,71.77)(-.03193,-.16436){10}{\line(0,-1){.16436}}
\multiput(88.07,70.12)(-.03361,-.20657){8}{\line(0,-1){.20657}}
\multiput(87.81,68.47)(-.03118,-.23715){7}{\line(0,-1){.23715}}
\multiput(87.59,66.81)(-.03348,-.33318){5}{\line(0,-1){.33318}}
\put(87.42,65.14){\line(0,-1){1.67}}
\put(87.3,63.47){\line(0,-1){1.673}}
\put(87.24,61.8){\line(0,-1){1.674}}
\put(87.22,60.13){\line(0,-1){1.674}}
\put(87.26,58.45){\line(0,-1){1.672}}
\multiput(87.35,56.78)(.02784,-.3337){5}{\line(0,-1){.3337}}
\multiput(87.49,55.11)(.03169,-.27725){6}{\line(0,-1){.27725}}
\multiput(87.68,53.45)(.03011,-.20711){8}{\line(0,-1){.20711}}
\multiput(87.92,51.79)(.03238,-.18319){9}{\line(0,-1){.18319}}
\multiput(88.21,50.14)(.03106,-.14901){11}{\line(0,-1){.14901}}
\multiput(88.55,48.5)(.03263,-.13566){12}{\line(0,-1){.13566}}
\multiput(88.94,46.88)(.031512,-.115367){14}{\line(0,-1){.115367}}
\multiput(89.39,45.26)(.032688,-.106727){15}{\line(0,-1){.106727}}
\multiput(89.88,43.66)(.033689,-.099073){16}{\line(0,-1){.099073}}
\multiput(90.42,42.07)(.032623,-.087108){18}{\line(0,-1){.087108}}
\multiput(91,40.51)(.033414,-.081541){19}{\line(0,-1){.081541}}
\multiput(91.64,38.96)(.032472,-.072816){21}{\line(0,-1){.072816}}
\multiput(92.32,37.43)(.033106,-.068527){22}{\line(0,-1){.068527}}
\multiput(93.05,35.92)(.033655,-.064549){23}{\line(0,-1){.064549}}
\multiput(93.82,34.44)(.032763,-.058411){25}{\line(0,-1){.058411}}
\multiput(94.64,32.98)(.033205,-.055175){26}{\line(0,-1){.055175}}
\multiput(95.5,31.54)(.033584,-.05213){27}{\line(0,-1){.05213}}
\multiput(96.41,30.13)(.032737,-.047556){29}{\line(0,-1){.047556}}
\multiput(97.36,28.75)(.033036,-.044983){30}{\line(0,-1){.044983}}
\multiput(98.35,27.41)(.033285,-.042534){31}{\line(0,-1){.042534}}
\multiput(99.38,26.09)(.033489,-.0402){32}{\line(0,-1){.0402}}
\multiput(100.45,24.8)(.033651,-.037971){33}{\line(0,-1){.037971}}
\multiput(101.57,23.55)(.0328073,-.034815){35}{\line(0,-1){.034815}}
\multiput(102.71,22.33)(.0329156,-.0328573){36}{\line(1,0){.0329156}}
\multiput(103.9,21.15)(.0358988,-.0337087){34}{\line(1,0){.0358988}}
\multiput(105.12,20)(.038031,-.033584){33}{\line(1,0){.038031}}
\multiput(106.37,18.89)(.040259,-.033418){32}{\line(1,0){.040259}}
\multiput(107.66,17.82)(.042593,-.03321){31}{\line(1,0){.042593}}
\multiput(108.98,16.79)(.045041,-.032956){30}{\line(1,0){.045041}}
\multiput(110.33,15.8)(.047614,-.032652){29}{\line(1,0){.047614}}
\multiput(111.71,14.86)(.052189,-.033492){27}{\line(1,0){.052189}}
\multiput(113.12,13.95)(.055234,-.033107){26}{\line(1,0){.055234}}
\multiput(114.56,13.09)(.058469,-.03266){25}{\line(1,0){.058469}}
\multiput(116.02,12.28)(.064608,-.033541){23}{\line(1,0){.064608}}
\multiput(117.51,11.5)(.068585,-.032985){22}{\line(1,0){.068585}}
\multiput(119.02,10.78)(.072874,-.032343){21}{\line(1,0){.072874}}
\multiput(120.55,10.1)(.0816,-.033269){19}{\line(1,0){.0816}}
\multiput(122.1,9.47)(.087166,-.032469){18}{\line(1,0){.087166}}
\multiput(123.67,8.88)(.099132,-.033513){16}{\line(1,0){.099132}}
\multiput(125.25,8.35)(.106784,-.032499){15}{\line(1,0){.106784}}
\multiput(126.85,7.86)(.124301,-.033716){13}{\line(1,0){.124301}}
\multiput(128.47,7.42)(.13571,-.03239){12}{\line(1,0){.13571}}
\multiput(130.1,7.03)(.14906,-.0308){11}{\line(1,0){.14906}}
\multiput(131.74,6.69)(.18325,-.03205){9}{\line(1,0){.18325}}
\multiput(133.39,6.4)(.20716,-.02974){8}{\line(1,0){.20716}}
\multiput(135.04,6.17)(.2773,-.0312){6}{\line(1,0){.2773}}
\multiput(136.71,5.98)(.33375,-.02725){5}{\line(1,0){.33375}}
\put(138.38,5.84){\line(1,0){1.672}}
\put(140.05,5.76){\line(1,0){1.674}}
\put(141.72,5.72){\line(1,0){1.674}}
\put(143.4,5.74){\line(1,0){1.673}}
\put(145.07,5.81){\line(1,0){1.67}}
\multiput(146.74,5.93)(.2776,.0284){6}{\line(1,0){.2776}}
\multiput(148.41,6.1)(.23709,.0316){7}{\line(1,0){.23709}}
\multiput(150.07,6.32)(.18357,.0302){9}{\line(1,0){.18357}}
\multiput(151.72,6.59)(.1643,.03222){10}{\line(1,0){.1643}}
\multiput(153.36,6.91)(.13603,.03102){12}{\line(1,0){.13603}}
\multiput(154.99,7.29)(.124635,.032459){13}{\line(1,0){.124635}}
\multiput(156.61,7.71)(.114758,.033663){14}{\line(1,0){.114758}}
\multiput(158.22,8.18)(.099466,.03251){16}{\line(1,0){.099466}}
\multiput(159.81,8.7)(.092636,.033445){17}{\line(1,0){.092636}}
\multiput(161.39,9.27)(.081932,.032444){19}{\line(1,0){.081932}}
\multiput(162.94,9.88)(.076857,.033186){20}{\line(1,0){.076857}}
\multiput(164.48,10.55)(.068915,.03229){22}{\line(1,0){.068915}}
\multiput(166,11.26)(.064944,.032886){23}{\line(1,0){.064944}}
\multiput(167.49,12.02)(.061246,.033404){24}{\line(1,0){.061246}}
\multiput(168.96,12.82)(.055566,.032548){26}{\line(1,0){.055566}}
\multiput(170.4,13.66)(.052525,.032963){27}{\line(1,0){.052525}}
\multiput(171.82,14.55)(.049654,.033319){28}{\line(1,0){.049654}}
\multiput(173.21,15.49)(.046936,.03362){29}{\line(1,0){.046936}}
\multiput(174.57,16.46)(.042926,.032778){31}{\line(1,0){.042926}}
\multiput(175.91,17.48)(.040595,.03301){32}{\line(1,0){.040595}}
\multiput(177.2,18.53)(.038368,.033198){33}{\line(1,0){.038368}}
\multiput(178.47,19.63)(.0362374,.0333444){34}{\line(1,0){.0362374}}
\multiput(179.7,20.76)(.0341957,.0334524){35}{\line(1,0){.0341957}}
\multiput(180.9,21.93)(.0331573,.0344819){35}{\line(0,1){.0344819}}
\multiput(182.06,23.14)(.0330317,.0365226){34}{\line(0,1){.0365226}}
\multiput(183.18,24.38)(.032867,.038652){33}{\line(0,1){.038652}}
\multiput(184.27,25.66)(.033713,.042196){31}{\line(0,1){.042196}}
\multiput(185.31,26.97)(.033488,.044647){30}{\line(0,1){.044647}}
\multiput(186.32,28.31)(.033215,.047223){29}{\line(0,1){.047223}}
\multiput(187.28,29.67)(.032891,.049938){28}{\line(0,1){.049938}}
\multiput(188.2,31.07)(.03251,.052806){27}{\line(0,1){.052806}}
\multiput(189.08,32.5)(.033352,.058077){25}{\line(0,1){.058077}}
\multiput(189.91,33.95)(.032876,.06153){24}{\line(0,1){.06153}}
\multiput(190.7,35.43)(.032327,.065224){23}{\line(0,1){.065224}}
\multiput(191.45,36.93)(.033206,.072485){21}{\line(0,1){.072485}}
\multiput(192.14,38.45)(.032524,.077139){20}{\line(0,1){.077139}}
\multiput(192.79,39.99)(.033501,.086774){18}{\line(0,1){.086774}}
\multiput(193.4,41.55)(.032648,.09292){17}{\line(0,1){.09292}}
\multiput(193.95,43.13)(.031654,.099742){16}{\line(0,1){.099742}}
\multiput(194.46,44.73)(.032676,.115043){14}{\line(0,1){.115043}}
\multiput(194.92,46.34)(.031386,.12491){13}{\line(0,1){.12491}}
\multiput(195.32,47.96)(.03256,.14869){11}{\line(0,1){.14869}}
\multiput(195.68,49.6)(.03081,.16457){10}{\line(0,1){.16457}}
\multiput(195.99,51.25)(.0322,.2068){8}{\line(0,1){.2068}}
\multiput(196.25,52.9)(.02956,.23735){7}{\line(0,1){.23735}}
\multiput(196.45,54.56)(.03121,.3334){5}{\line(0,1){.3334}}
\put(196.61,56.23){\line(0,1){1.671}}
\put(196.72,57.9){\line(0,1){2.6}}
\put(112.75,60){\circle{11.6}} \put(142,60){\circle{11.6}}
\put(171.75,59.75){\circle{11.6}}
\put(142.25,115.5){\line(0,-1){49.5}}
\put(142,54.25){\line(0,-1){48.5}}
\put(87,60.25){\line(1,0){20}}
\put(118.75,59.75){\line(1,0){17.25}}
\put(147.75,60){\line(1,0){18.25}}
\put(177.5,60){\line(1,0){19.25}}
\qbezier(114.25,65.75)(141.88,111)(170,65.25)
\qbezier(113.75,54.5)(142.88,5.63)(170.5,54.25)
\put(96.75,65.25){\vector(0,1){.07}}\put(96.75,56.75){\line(0,1){8.5}}
\put(127.5,56.75){\vector(0,-1){.07}}\put(127.5,63.75){\line(0,-1){7}}
\put(158.25,65.25){\vector(0,1){.07}}\put(158.25,58){\line(0,1){7.25}}
\put(185.25,56){\vector(0,-1){.07}}\put(185.25,63.75){\line(0,-1){7.75}}

\put(100.25,57.75){\makebox(0,0)[cc]{3}}
\put(112,59.5){\makebox(0,0)[cc]{2}}
\put(141.75,60.25){\makebox(0,0)[cc]{4}}
\put(171.75,59.25){\makebox(0,0)[cc]{2}}
\put(145.25,104.5){\makebox(0,0)[cc]{5}}
\put(122,81.75){\makebox(0,0)[cc]{6}}
\put(123.25,63.25){\makebox(0,0)[cc]{3}}
\put(154.75,63){\makebox(0,0)[cc]{3}}
\put(131.5,38.25){\makebox(0,0)[cc]{1}}
\put(22.5,40.25){\makebox(0,0)[cc]{$P_1(7)$}}
\put(189,63.5){\makebox(0,0)[cc]{3}}
\put(185.75,98.75){\makebox(0,0)[cc]{4}}
\put(144.5,14){\makebox(0,0)[cc]{7}}
\end{picture}
\end{center}

\proof Let $X=\{2,4\}$. The order on $P_1(7)$:
$2\succ4\succ3\succ1\succ5\succ6\succ7$. The partial orders on links
(we can extend these partial orders to linear orders arbitrarily):

\begin{itemize}
\item 2: $5\succ_27$,
\item 4: $1\succ_46$.
\end{itemize}

Let $p=(v_1,...,v_k)$ be an anti-path in $P_1(7)$ where only $v_1$
may belong to $X$. Suppose that $\Theta=\Theta(p)$ does not satisfy
(*) for $S[X]=S_0$. Consider seven cases for $v_k$.

{\bf 1.} Let $v_k=1$. Then $\Theta\supseteq \{1,3,5\}$ and satisfies
(*).

{\bf 2.} Let $v_k=2$ (and $k=1$). Then $\Theta\supseteq \{5,7\}$ and
satisfies (*).

{\bf 3.} Let $v_k=3$. Then $\Theta\supseteq \{1,3,5\}$ and satisfies
(*).

{\bf 4.} Let $v_k=4$ (and $k=1$). Then $\Theta\supseteq \{1,6\}$ and
satisfies (*).

{\bf 5.} Let $v_k=5$. Then $\Theta\supseteq \{1,3,5\}$ and satisfies
(*).

{\bf 6.} Let $v_k=6$. Then $\Theta\supseteq \{6,7\}$. Since
$\{3,6,7\}$ satisfies (*), and $3\succ6$, we have $k\ne 1$. Then
$v_{k-1}\in \{4,7\}$. Since $1\succ_46$ and $\{1,6,7\}$ satisfies
(*), we have $v_{k-1}=7$. Since $1\succ7$, $k-1\ne 1$, $v_{k-2}=2$.
But $5\succ_27$, hence $\Theta\supseteq \{5,6,7\}$ and satisfies
(*).

{\bf 7.} Let $v_k=7$. Then $\Theta\supseteq \{6,7\}$. Since
$5\succ7$ and $\{5,6,7\}$ satisfies (*), $k\ne 1$. Then $v_{k-1}\in
\{2,6\}$. Since $5\succ_27$, we have $v_{k-1}=6$.

Since $5\succ6$, $k-1\ne 1$, $v_{k-2}=4$. But $1\succ_46$, so $1\in
\Theta$ and $\Theta\supseteq \{1,6,7\}$ that satisfies (*), a
contradiction.
\endproof

\subsubsection{The 7-vertex graph $P_2(7)$ and the 8-vertex graph
$P_1(8)$} \label{p18}

\begin{center}
\unitlength .5mm 
\linethickness{0.4pt}
\ifx\plotpoint\undefined\newsavebox{\plotpoint}\fi 
\begin{picture}(285,137.5)(0,0)
\put(193,75.75){\circle{9.22}} \put(264,75.75){\circle{9.22}}
\put(228.5,75.75){\circle{9.22}}
\put(228.5,80.5){\line(0,1){52}}
\put(228.5,71){\line(0,-1){51}}
\put(171.75,76.25){\line(1,0){16.75}}
\put(197.5,76.25){\line(1,0){26}}
\multiput(232.75,76)(8.75,-.0833){3}{\line(1,0){8.75}}
\put(268.5,75.75){\line(1,0){16.5}}
\qbezier(193.75,80.5)(229.75,131.63)(264.75,80.25)
\qbezier(194,71)(229.25,19.88)(263.5,71.25)
\qbezier(180,76.25)(179.88,118)(228.25,118.75)
\qbezier(228.25,118.75)(275,119.13)(275.75,76)
\qbezier(275.75,76)(276.38,58.88)(267.5,59.25)
\qbezier(267.75,59.25)(249.13,59.38)(241,90)
\qbezier(241,89.5)(239.38,94.25)(228.25,93)
\qbezier(228.25,93)(216.25,94)(213.25,90)
\qbezier(213.25,90)(209.88,86)(192,58)
\qbezier(192,58)(180,45.75)(180,76.5)
\put(.5,48){\line(1,0){72.5}}
\put(84.75,48){\line(1,0){72.5}}
\multiput(73,48)(-.19458763,.11984536){194}{\line(-1,0){.19458763}}
\multiput(157.25,48)(-.19458763,.11984536){194}{\line(-1,0){.19458763}}
\multiput(35.25,71.25)(-.1796875,-.11979167){192}{\line(-1,0){.1796875}}
\multiput(119.5,71.25)(-.1796875,-.11979167){192}{\line(-1,0){.1796875}}
\put(35.25,70.75){\line(0,-1){42.5}}
\put(119.5,70.75){\line(0,-1){42.5}}
\multiput(35.25,28.25)(.23159509,.1196319){163}{\line(1,0){.23159509}}
\multiput(119.5,28.25)(.23159509,.1196319){163}{\line(1,0){.23159509}}
\multiput(1,48.25)(.20858896,-.1196319){163}{\line(1,0){.20858896}}
\multiput(85.25,48.25)(.20858896,-.1196319){163}{\line(1,0){.20858896}}
\multiput(18.75,48)(.11956522,.16485507){138}{\line(0,1){.16485507}}
\multiput(35.25,70.75)(.11931818,-.17045455){132}{\line(0,-1){.17045455}}
\multiput(119.5,70.75)(.11931818,-.17045455){132}{\line(0,-1){.17045455}}
\put(1.5,43.75){\makebox(0,0)[cc]{1}}
\put(85.75,43.75){\makebox(0,0)[cc]{1}}
\put(18.75,45){\makebox(0,0)[cc]{2}}
\put(103,45){\makebox(0,0)[cc]{2}}
\put(38.5,44.5){\makebox(0,0)[cc]{3}}
\put(122.75,44.5){\makebox(0,0)[cc]{3}}
\put(51,44.25){\makebox(0,0)[cc]{4}}
\put(135.25,44.25){\makebox(0,0)[cc]{4}}
\put(72,44){\makebox(0,0)[cc]{5}}
\put(156.25,44){\makebox(0,0)[cc]{5}}
\put(35.25,74.75){\makebox(0,0)[cc]{6}}
\put(119.5,74.75){\makebox(0,0)[cc]{6}}
\put(34.5,23){\makebox(0,0)[cc]{7}}
\put(118.75,23){\makebox(0,0)[cc]{7}}
\put(208,133.5){\makebox(0,0)[cc]{4}}
\put(192.5,75){\makebox(0,0)[cc]{2}}
\put(225.25,85.25){\makebox(0,0)[cc]{5}}
\put(228.5,75.5){\makebox(0,0)[cc]{4}}
\put(263.75,75.75){\makebox(0,0)[cc]{2}}
\put(175.25,71.5){\makebox(0,0)[cc]{3}}
\put(250.25,80){\makebox(0,0)[cc]{3}}
\put(216.25,79.25){\makebox(0,0)[cc]{3}}
\put(248.5,112.25){\makebox(0,0)[cc]{7}}
\put(225.25,110.25){\makebox(0,0)[cc]{$6$ or $0$}}
\put(225.25,35.25){\makebox(0,0)[cc]{6}}
\put(212.25,54.75){\makebox(0,0)[cc]{1}}
\multiput(119.5,71)(-.1194444,-.1388889){90}{\line(0,-1){.1388889}}
\multiput(108.75,58.5)(-.1184211,-.1798246){57}{\line(0,-1){.1798246}}
\multiput(108.75,58.25)(-.2767857,-.1190476){84}{\line(-1,0){.2767857}}
\multiput(108.5,58.25)(.1310976,-.1189024){82}{\line(1,0){.1310976}}
\multiput(108.25,58)(.5884146,-.1189024){82}{\line(1,0){.5884146}}
\put(108.75,55.25){\makebox(0,0)[cc]{0}}
\put(8.25,64){\makebox(0,0)[cc]{$P_2(7)$}}
\put(98,66.25){\makebox(0,0)[cc]{$P_1(8)$}}
\put(228.13,76.38){\oval(113.25,111.75)[]}
\put(280.25,79.75){\makebox(0,0)[cc]{3}}
\put(183.5,72.25){\vector(0,-1){.7}}\multiput(183.25,80.25)(.0833,-2.6667){3}{\line(0,-1){2.6667}}
\put(211,82.5){\vector(0,1){.7}}\multiput(211.5,73.75)(-.1,1.75){5}{\line(0,1){1.75}}
\put(239.75,70){\vector(0,-1){.7}}\put(239.75,81.75){\line(0,-1){11.75}}
\put(280.5,79.5){\vector(0,1){.7}}\put(280.5,73.25){\line(0,1){6.25}}
\end{picture}

The top arc intersecting both $2$-curves on the dissection diagram
is labeled by 6 for $P_2(7)$ and by $0$ for $P_1(8)$.
\end{center}

\proof Let $X=\{2,4\}$. The order on the graphs is:
$2\succ4\succ3\succ0\succ5\succ6\succ7$, partial orders on the
anti-links of vertices (the inequalities involving 0 - for $P_1(8)$
only):

\begin{itemize}
\item 2: $5\succ_27$, $5\succ_2 6$,
\item 4: $1\succ_40$, $1\succ_47$,
\end{itemize}

Assume that $\Theta=\Theta(p)$ does not satisfy (*) for $S[X]=S_0$
for some anti-path $p=(v_1,...,v_k)$.

{\bf 0 (for $P_1(8)$).} Let $v_k=0$. Then $\Theta\supseteq \{0,7\}$.
Since $3\succ0$ and $\{0,3,7\}$ satisfies (*), $k\ne 1$. Then
$v_{k-1}\in \{4,7\}$. But $1\succ_40$ and $\{0,1,7\}$ satisfies (*).
So $v_{k-1}=7$.

Since $3\succ7$, $k-1\ne 1$, $v_{k-2}\in \{2,4,6\}$. But
$5\succ_27$, so $v_{k-2}\ne 2$, $1\succ_47$ and $\{0,1,7\}$
satisfies (*), so $v_{k-2}\ne 4$. Thus $v_{k-2}=6$.

Since $3\succ6$, $k-2\ne 1$. Then $v_{k-3}=2$. But $5\succ_26$, so
$5\in \Theta$, a contradiction.

{\bf 1.} Let $v_k=1$. Then $\Theta\supseteq \{1,3,5\}$ and satisfies
(*).

{\bf 2.} Let $v_k=2$ (and $k=1$). Then $\Theta\supseteq \{5,7\}$ and
satisfies (*).

{\bf 3.} Let $v_k=3$. Then $\Theta\supseteq \{1,3, 5\}$ and
satisfies (*).

{\bf 4.} Let $v_k=4$ (and $k=1$). Then $\Theta\supseteq \{1,7\}$ and
satisfies (*).

{\bf 5.} Let $v_k=5$. Then $\Theta\supseteq \{1,3,5\}$ and satisfies
(*).

{\bf 6.} Let $v_k=6$. Then $\Theta\supseteq \{6,7\}$. Since
$\{3,6,7\}$ satisfies (*), and $3\succ6$, we have $k\ne 1$. Then
$v_{k-1}=7$.

Since $3\succ7$, $k-1\ne 1$. Then $v_{k-2}\in \{0,2,4\}$ ($0$ - in
the case of $P_1(8)$). Since $5\succ_27$, $1\succ_47$, and sets
$\{5,6,7\}$, $\{1,6,7\}$ satisfy (*), $v_{k-2}=0$ and the graph is
$P_1(8)$. Since $3\succ0$, $k-2\ne 1$. Then $v_{k-3}=4$. But
$1\succ_40$, $\{1,6,7\}$ satisfies (*), a contradiction.

{\bf 7.} Let $v_k=7$. Then $\Theta\supseteq \{6,7\}$ (and
$\Theta\supseteq \{0,6,7\}$ in the case of $P_1(8)$) and satisfies
(*).
\endproof

\subsubsection{The 8-vertex graph $P_2(8)$}

\begin{center}
\begin{minipage}[t]{2 in}
\unitlength .5mm 
\linethickness{0.4pt}
\ifx\plotpoint\undefined\newsavebox{\plotpoint}\fi 
\begin{picture}(110.25,90)(15,10)
\multiput(14.5,46.5)(11.65625,-.03125){8}{\line(1,0){11.65625}}
\qbezier(14.75,46.75)(58,125.5)(107.25,46.25)
\multiput(34.75,46.5)(.03920742,.033726813){593}{\line(1,0){.03920742}}
\multiput(58,66.5)(-.03125,-2.4375){8}{\line(0,-1){2.4375}}
\multiput(57.75,57.25)(-.072916667,-.033653846){312}{\line(-1,0){.072916667}}
\multiput(35,46.5)(.033333,-1.75){15}{\line(0,-1){1.75}}
\multiput(35.5,20.25)(.0337331334,.0397301349){667}{\line(0,1){.0397301349}}
\multiput(14.75,46.75)(.033717105,-.042763158){608}{\line(0,-1){.042763158}}
\multiput(35.25,34.25)(.0625,.033653846){364}{\line(1,0){.0625}}
\multiput(57.75,66.5)(.083052277,-.033726813){593}{\line(1,0){.083052277}}
\put(16, 60){\makebox(0,0)[cc]{$P_2(8)$}}

 \put(12,45.25){\makebox(0,0)[cc]{0}}
\put(32.25,43.5){\makebox(0,0)[cc]{7}}
\put(59,44){\makebox(0,0)[cc]{6}}
\put(110.25,46){\makebox(0,0)[cc]{1}}
\put(57.25,70.25){\makebox(0,0)[cc]{3}}
\put(60,57.5){\makebox(0,0)[cc]{5}}
\put(31.5,35.5){\makebox(0,0)[cc]{4}}
\put(39.5,21.5){\makebox(0,0)[cc]{2}}
\end{picture}
\end{minipage}
\begin{minipage}[t]{4 in}
\unitlength .6mm 
\linethickness{0.4pt}
\ifx\plotpoint\undefined\newsavebox{\plotpoint}\fi 
\begin{picture}(187.75,50.5)(20,10)
\put(160.13,43){\line(0,1){.418}}
\put(160.12,43.42){\line(0,1){.417}}
\put(160.08,43.84){\line(0,1){.414}}
\put(160.02,44.25){\line(0,1){.41}}
\put(159.94,44.66){\line(0,1){.405}}
\multiput(159.84,45.07)(-.0321,.0996){4}{\line(0,1){.0996}}
\multiput(159.71,45.46)(-.03018,.07809){5}{\line(0,1){.07809}}
\multiput(159.56,45.85)(-.0288,.06355){6}{\line(0,1){.06355}}
\multiput(159.38,46.24)(-.03235,.06181){6}{\line(0,1){.06181}}
\multiput(159.19,46.61)(-.03069,.05132){7}{\line(0,1){.05132}}
\multiput(158.97,46.97)(-.03355,.0495){7}{\line(0,1){.0495}}
\multiput(158.74,47.31)(-.03176,.04158){8}{\line(0,1){.04158}}
\multiput(158.49,47.64)(-.03029,.0353){9}{\line(0,1){.0353}}
\multiput(158.21,47.96)(-.03224,.03353){9}{\line(0,1){.03353}}
\multiput(157.92,48.26)(-.03409,.03164){9}{\line(-1,0){.03409}}
\multiput(157.62,48.55)(-.0403,.03337){8}{\line(-1,0){.0403}}
\multiput(157.29,48.82)(-.04213,.03103){8}{\line(-1,0){.04213}}
\multiput(156.96,49.06)(-.05008,.03267){7}{\line(-1,0){.05008}}
\multiput(156.61,49.29)(-.05186,.02978){7}{\line(-1,0){.05186}}
\multiput(156.24,49.5)(-.06237,.03126){6}{\line(-1,0){.06237}}
\multiput(155.87,49.69)(-.07685,.03321){5}{\line(-1,0){.07685}}
\multiput(155.48,49.85)(-.07861,.0288){5}{\line(-1,0){.07861}}
\multiput(155.09,50)(-.1001,.0304){4}{\line(-1,0){.1001}}
\put(154.69,50.12){\line(-1,0){.407}}
\put(154.28,50.22){\line(-1,0){.412}}
\put(153.87,50.29){\line(-1,0){.415}}
\put(153.46,50.35){\line(-1,0){.418}}
\put(153.04,50.37){\line(-1,0){.419}}
\put(152.62,50.38){\line(-1,0){.418}}
\put(152.2,50.36){\line(-1,0){.416}}
\put(151.79,50.32){\line(-1,0){.413}}
\put(151.37,50.25){\line(-1,0){.409}}
\put(150.96,50.16){\line(-1,0){.403}}
\multiput(150.56,50.05)(-.07921,-.0271){5}{\line(-1,0){.07921}}
\multiput(150.16,49.91)(-.07755,-.03154){5}{\line(-1,0){.07755}}
\multiput(149.78,49.75)(-.06303,-.02991){6}{\line(-1,0){.06303}}
\multiput(149.4,49.58)(-.06123,-.03343){6}{\line(-1,0){.06123}}
\multiput(149.03,49.37)(-.05078,-.03159){7}{\line(-1,0){.05078}}
\multiput(148.68,49.15)(-.04279,-.03011){8}{\line(-1,0){.04279}}
\multiput(148.33,48.91)(-.04101,-.03249){8}{\line(-1,0){.04101}}
\multiput(148.01,48.65)(-.03476,-.0309){9}{\line(-1,0){.03476}}
\multiput(147.69,48.37)(-.03295,-.03282){9}{\line(-1,0){.03295}}
\multiput(147.4,48.08)(-.03104,-.03464){9}{\line(0,-1){.03464}}
\multiput(147.12,47.77)(-.03265,-.04089){8}{\line(0,-1){.04089}}
\multiput(146.86,47.44)(-.03028,-.04267){8}{\line(0,-1){.04267}}
\multiput(146.61,47.1)(-.03179,-.05065){7}{\line(0,-1){.05065}}
\multiput(146.39,46.74)(-.03368,-.0611){6}{\line(0,-1){.0611}}
\multiput(146.19,46.38)(-.03016,-.06291){6}{\line(0,-1){.06291}}
\multiput(146.01,46)(-.03185,-.07742){5}{\line(0,-1){.07742}}
\multiput(145.85,45.61)(-.02741,-.0791){5}{\line(0,-1){.0791}}
\put(145.71,45.22){\line(0,-1){.403}}
\put(145.6,44.81){\line(0,-1){.408}}
\put(145.51,44.41){\line(0,-1){.413}}
\put(145.44,43.99){\line(0,-1){.416}}
\put(145.39,43.58){\line(0,-1){1.67}}
\put(145.45,41.91){\line(0,-1){.412}}
\put(145.53,41.5){\line(0,-1){.407}}
\put(145.62,41.09){\line(0,-1){.401}}
\multiput(145.74,40.69)(.02848,-.07872){5}{\line(0,-1){.07872}}
\multiput(145.88,40.29)(.0329,-.07698){5}{\line(0,-1){.07698}}
\multiput(146.05,39.91)(.03101,-.06249){6}{\line(0,-1){.06249}}
\multiput(146.23,39.53)(.02958,-.05197){7}{\line(0,-1){.05197}}
\multiput(146.44,39.17)(.03248,-.05021){7}{\line(0,-1){.05021}}
\multiput(146.67,38.82)(.03086,-.04225){8}{\line(0,-1){.04225}}
\multiput(146.92,38.48)(.03321,-.04044){8}{\line(0,-1){.04044}}
\multiput(147.18,38.16)(.03151,-.03421){9}{\line(0,-1){.03421}}
\multiput(147.47,37.85)(.0334,-.03237){9}{\line(1,0){.0334}}
\multiput(147.77,37.56)(.03518,-.03043){9}{\line(1,0){.03518}}
\multiput(148.08,37.28)(.04145,-.03193){8}{\line(1,0){.04145}}
\multiput(148.41,37.03)(.0432,-.02953){8}{\line(1,0){.0432}}
\multiput(148.76,36.79)(.0512,-.03089){7}{\line(1,0){.0512}}
\multiput(149.12,36.58)(.06168,-.0326){6}{\line(1,0){.06168}}
\multiput(149.49,36.38)(.06343,-.02905){6}{\line(1,0){.06343}}
\multiput(149.87,36.21)(.07797,-.03048){5}{\line(1,0){.07797}}
\multiput(150.26,36.05)(.0995,-.0325){4}{\line(1,0){.0995}}
\put(150.66,35.92){\line(1,0){.405}}
\put(151.06,35.82){\line(1,0){.41}}
\put(151.47,35.73){\line(1,0){.414}}
\put(151.89,35.67){\line(1,0){.417}}
\put(152.3,35.63){\line(1,0){.837}}
\put(153.14,35.63){\line(1,0){.417}}
\put(153.56,35.66){\line(1,0){.415}}
\put(153.97,35.72){\line(1,0){.411}}
\put(154.38,35.8){\line(1,0){.405}}
\multiput(154.79,35.91)(.0997,.0317){4}{\line(1,0){.0997}}
\multiput(155.19,36.03)(.07821,.02986){5}{\line(1,0){.07821}}
\multiput(155.58,36.18)(.06366,.02854){6}{\line(1,0){.06366}}
\multiput(155.96,36.35)(.06194,.03211){6}{\line(1,0){.06194}}
\multiput(156.33,36.55)(.05145,.03049){7}{\line(1,0){.05145}}
\multiput(156.69,36.76)(.04963,.03335){7}{\line(1,0){.04963}}
\multiput(157.04,36.99)(.04171,.0316){8}{\line(1,0){.04171}}
\multiput(157.37,37.25)(.03542,.03015){9}{\line(1,0){.03542}}
\multiput(157.69,37.52)(.03365,.0321){9}{\line(1,0){.03365}}
\multiput(157.99,37.81)(.03178,.03396){9}{\line(0,1){.03396}}
\multiput(158.28,38.11)(.03353,.04017){8}{\line(0,1){.04017}}
\multiput(158.55,38.43)(.0312,.04201){8}{\line(0,1){.04201}}
\multiput(158.8,38.77)(.03287,.04995){7}{\line(0,1){.04995}}
\multiput(159.03,39.12)(.02999,.05174){7}{\line(0,1){.05174}}
\multiput(159.24,39.48)(.03151,.06225){6}{\line(0,1){.06225}}
\multiput(159.43,39.86)(.03351,.07672){5}{\line(0,1){.07672}}
\multiput(159.59,40.24)(.02911,.0785){5}{\line(0,1){.0785}}
\multiput(159.74,40.63)(.0308,.1){4}{\line(0,1){.1}}
\put(159.86,41.03){\line(0,1){.406}}
\put(159.96,41.44){\line(0,1){.411}}
\put(160.04,41.85){\line(0,1){.415}}
\put(160.09,42.26){\line(0,1){.736}}
\put(60.13,43){\line(0,1){.418}} \put(60.12,43.42){\line(0,1){.417}}
\put(60.08,43.84){\line(0,1){.414}}
\put(60.02,44.25){\line(0,1){.41}}
\put(59.94,44.66){\line(0,1){.405}}
\multiput(59.84,45.07)(-.0321,.0996){4}{\line(0,1){.0996}}
\multiput(59.71,45.46)(-.03018,.07809){5}{\line(0,1){.07809}}
\multiput(59.56,45.85)(-.0288,.06355){6}{\line(0,1){.06355}}
\multiput(59.38,46.24)(-.03235,.06181){6}{\line(0,1){.06181}}
\multiput(59.19,46.61)(-.03069,.05132){7}{\line(0,1){.05132}}
\multiput(58.97,46.97)(-.03355,.0495){7}{\line(0,1){.0495}}
\multiput(58.74,47.31)(-.03176,.04158){8}{\line(0,1){.04158}}
\multiput(58.49,47.64)(-.03029,.0353){9}{\line(0,1){.0353}}
\multiput(58.21,47.96)(-.03224,.03353){9}{\line(0,1){.03353}}
\multiput(57.92,48.26)(-.03409,.03164){9}{\line(-1,0){.03409}}
\multiput(57.62,48.55)(-.0403,.03337){8}{\line(-1,0){.0403}}
\multiput(57.29,48.82)(-.04213,.03103){8}{\line(-1,0){.04213}}
\multiput(56.96,49.06)(-.05008,.03267){7}{\line(-1,0){.05008}}
\multiput(56.61,49.29)(-.05186,.02978){7}{\line(-1,0){.05186}}
\multiput(56.24,49.5)(-.06237,.03126){6}{\line(-1,0){.06237}}
\multiput(55.87,49.69)(-.07685,.03321){5}{\line(-1,0){.07685}}
\multiput(55.48,49.85)(-.07861,.0288){5}{\line(-1,0){.07861}}
\multiput(55.09,50)(-.1001,.0304){4}{\line(-1,0){.1001}}
\put(54.69,50.12){\line(-1,0){.407}}
\put(54.28,50.22){\line(-1,0){.412}}
\put(53.87,50.29){\line(-1,0){.415}}
\put(53.46,50.35){\line(-1,0){.418}}
\put(53.04,50.37){\line(-1,0){.419}}
\put(52.62,50.38){\line(-1,0){.418}}
\put(52.2,50.36){\line(-1,0){.416}}
\put(51.79,50.32){\line(-1,0){.413}}
\put(51.37,50.25){\line(-1,0){.409}}
\put(50.96,50.16){\line(-1,0){.403}}
\multiput(50.56,50.05)(-.07921,-.0271){5}{\line(-1,0){.07921}}
\multiput(50.16,49.91)(-.07755,-.03154){5}{\line(-1,0){.07755}}
\multiput(49.78,49.75)(-.06303,-.02991){6}{\line(-1,0){.06303}}
\multiput(49.4,49.58)(-.06123,-.03343){6}{\line(-1,0){.06123}}
\multiput(49.03,49.37)(-.05078,-.03159){7}{\line(-1,0){.05078}}
\multiput(48.68,49.15)(-.04279,-.03011){8}{\line(-1,0){.04279}}
\multiput(48.33,48.91)(-.04101,-.03249){8}{\line(-1,0){.04101}}
\multiput(48.01,48.65)(-.03476,-.0309){9}{\line(-1,0){.03476}}
\multiput(47.69,48.37)(-.03295,-.03282){9}{\line(-1,0){.03295}}
\multiput(47.4,48.08)(-.03104,-.03464){9}{\line(0,-1){.03464}}
\multiput(47.12,47.77)(-.03265,-.04089){8}{\line(0,-1){.04089}}
\multiput(46.86,47.44)(-.03028,-.04267){8}{\line(0,-1){.04267}}
\multiput(46.61,47.1)(-.03179,-.05065){7}{\line(0,-1){.05065}}
\multiput(46.39,46.74)(-.03368,-.0611){6}{\line(0,-1){.0611}}
\multiput(46.19,46.38)(-.03016,-.06291){6}{\line(0,-1){.06291}}
\multiput(46.01,46)(-.03185,-.07742){5}{\line(0,-1){.07742}}
\multiput(45.85,45.61)(-.02741,-.0791){5}{\line(0,-1){.0791}}
\put(45.71,45.22){\line(0,-1){.403}}
\put(45.6,44.81){\line(0,-1){.408}}
\put(45.51,44.41){\line(0,-1){.413}}
\put(45.44,43.99){\line(0,-1){.416}}
\put(45.39,43.58){\line(0,-1){1.67}}
\put(45.45,41.91){\line(0,-1){.412}}
\put(45.53,41.5){\line(0,-1){.407}}
\put(45.62,41.09){\line(0,-1){.401}}
\multiput(45.74,40.69)(.02848,-.07872){5}{\line(0,-1){.07872}}
\multiput(45.88,40.29)(.0329,-.07698){5}{\line(0,-1){.07698}}
\multiput(46.05,39.91)(.03101,-.06249){6}{\line(0,-1){.06249}}
\multiput(46.23,39.53)(.02958,-.05197){7}{\line(0,-1){.05197}}
\multiput(46.44,39.17)(.03248,-.05021){7}{\line(0,-1){.05021}}
\multiput(46.67,38.82)(.03086,-.04225){8}{\line(0,-1){.04225}}
\multiput(46.92,38.48)(.03321,-.04044){8}{\line(0,-1){.04044}}
\multiput(47.18,38.16)(.03151,-.03421){9}{\line(0,-1){.03421}}
\multiput(47.47,37.85)(.0334,-.03237){9}{\line(1,0){.0334}}
\multiput(47.77,37.56)(.03518,-.03043){9}{\line(1,0){.03518}}
\multiput(48.08,37.28)(.04145,-.03193){8}{\line(1,0){.04145}}
\multiput(48.41,37.03)(.0432,-.02953){8}{\line(1,0){.0432}}
\multiput(48.76,36.79)(.0512,-.03089){7}{\line(1,0){.0512}}
\multiput(49.12,36.58)(.06168,-.0326){6}{\line(1,0){.06168}}
\multiput(49.49,36.38)(.06343,-.02905){6}{\line(1,0){.06343}}
\multiput(49.87,36.21)(.07797,-.03048){5}{\line(1,0){.07797}}
\multiput(50.26,36.05)(.0995,-.0325){4}{\line(1,0){.0995}}
\put(50.66,35.92){\line(1,0){.405}}
\put(51.06,35.82){\line(1,0){.41}}
\put(51.47,35.73){\line(1,0){.414}}
\put(51.89,35.67){\line(1,0){.417}}
\put(52.3,35.63){\line(1,0){.837}}
\put(53.14,35.63){\line(1,0){.417}}
\put(53.56,35.66){\line(1,0){.415}}
\put(53.97,35.72){\line(1,0){.411}}
\put(54.38,35.8){\line(1,0){.405}}
\multiput(54.79,35.91)(.0997,.0317){4}{\line(1,0){.0997}}
\multiput(55.19,36.03)(.07821,.02986){5}{\line(1,0){.07821}}
\multiput(55.58,36.18)(.06366,.02854){6}{\line(1,0){.06366}}
\multiput(55.96,36.35)(.06194,.03211){6}{\line(1,0){.06194}}
\multiput(56.33,36.55)(.05145,.03049){7}{\line(1,0){.05145}}
\multiput(56.69,36.76)(.04963,.03335){7}{\line(1,0){.04963}}
\multiput(57.04,36.99)(.04171,.0316){8}{\line(1,0){.04171}}
\multiput(57.37,37.25)(.03542,.03015){9}{\line(1,0){.03542}}
\multiput(57.69,37.52)(.03365,.0321){9}{\line(1,0){.03365}}
\multiput(57.99,37.81)(.03178,.03396){9}{\line(0,1){.03396}}
\multiput(58.28,38.11)(.03353,.04017){8}{\line(0,1){.04017}}
\multiput(58.55,38.43)(.0312,.04201){8}{\line(0,1){.04201}}
\multiput(58.8,38.77)(.03287,.04995){7}{\line(0,1){.04995}}
\multiput(59.03,39.12)(.02999,.05174){7}{\line(0,1){.05174}}
\multiput(59.24,39.48)(.03151,.06225){6}{\line(0,1){.06225}}
\multiput(59.43,39.86)(.03351,.07672){5}{\line(0,1){.07672}}
\multiput(59.59,40.24)(.02911,.0785){5}{\line(0,1){.0785}}
\multiput(59.74,40.63)(.0308,.1){4}{\line(0,1){.1}}
\put(59.86,41.03){\line(0,1){.406}}
\put(59.96,41.44){\line(0,1){.411}}
\put(60.04,41.85){\line(0,1){.415}}
\put(60.09,42.26){\line(0,1){.736}}
\put(84.13,43){\line(0,1){.418}} \put(84.12,43.42){\line(0,1){.417}}
\put(84.08,43.84){\line(0,1){.414}}
\put(84.02,44.25){\line(0,1){.41}}
\put(83.94,44.66){\line(0,1){.405}}
\multiput(83.84,45.07)(-.0321,.0996){4}{\line(0,1){.0996}}
\multiput(83.71,45.46)(-.03018,.07809){5}{\line(0,1){.07809}}
\multiput(83.56,45.85)(-.0288,.06355){6}{\line(0,1){.06355}}
\multiput(83.38,46.24)(-.03235,.06181){6}{\line(0,1){.06181}}
\multiput(83.19,46.61)(-.03069,.05132){7}{\line(0,1){.05132}}
\multiput(82.97,46.97)(-.03355,.0495){7}{\line(0,1){.0495}}
\multiput(82.74,47.31)(-.03176,.04158){8}{\line(0,1){.04158}}
\multiput(82.49,47.64)(-.03029,.0353){9}{\line(0,1){.0353}}
\multiput(82.21,47.96)(-.03224,.03353){9}{\line(0,1){.03353}}
\multiput(81.92,48.26)(-.03409,.03164){9}{\line(-1,0){.03409}}
\multiput(81.62,48.55)(-.0403,.03337){8}{\line(-1,0){.0403}}
\multiput(81.29,48.82)(-.04213,.03103){8}{\line(-1,0){.04213}}
\multiput(80.96,49.06)(-.05008,.03267){7}{\line(-1,0){.05008}}
\multiput(80.61,49.29)(-.05186,.02978){7}{\line(-1,0){.05186}}
\multiput(80.24,49.5)(-.06237,.03126){6}{\line(-1,0){.06237}}
\multiput(79.87,49.69)(-.07685,.03321){5}{\line(-1,0){.07685}}
\multiput(79.48,49.85)(-.07861,.0288){5}{\line(-1,0){.07861}}
\multiput(79.09,50)(-.1001,.0304){4}{\line(-1,0){.1001}}
\put(78.69,50.12){\line(-1,0){.407}}
\put(78.28,50.22){\line(-1,0){.412}}
\put(77.87,50.29){\line(-1,0){.415}}
\put(77.46,50.35){\line(-1,0){.418}}
\put(77.04,50.37){\line(-1,0){.419}}
\put(76.62,50.38){\line(-1,0){.418}}
\put(76.2,50.36){\line(-1,0){.416}}
\put(75.79,50.32){\line(-1,0){.413}}
\put(75.37,50.25){\line(-1,0){.409}}
\put(74.96,50.16){\line(-1,0){.403}}
\multiput(74.56,50.05)(-.07921,-.0271){5}{\line(-1,0){.07921}}
\multiput(74.16,49.91)(-.07755,-.03154){5}{\line(-1,0){.07755}}
\multiput(73.78,49.75)(-.06303,-.02991){6}{\line(-1,0){.06303}}
\multiput(73.4,49.58)(-.06123,-.03343){6}{\line(-1,0){.06123}}
\multiput(73.03,49.37)(-.05078,-.03159){7}{\line(-1,0){.05078}}
\multiput(72.68,49.15)(-.04279,-.03011){8}{\line(-1,0){.04279}}
\multiput(72.33,48.91)(-.04101,-.03249){8}{\line(-1,0){.04101}}
\multiput(72.01,48.65)(-.03476,-.0309){9}{\line(-1,0){.03476}}
\multiput(71.69,48.37)(-.03295,-.03282){9}{\line(-1,0){.03295}}
\multiput(71.4,48.08)(-.03104,-.03464){9}{\line(0,-1){.03464}}
\multiput(71.12,47.77)(-.03265,-.04089){8}{\line(0,-1){.04089}}
\multiput(70.86,47.44)(-.03028,-.04267){8}{\line(0,-1){.04267}}
\multiput(70.61,47.1)(-.03179,-.05065){7}{\line(0,-1){.05065}}
\multiput(70.39,46.74)(-.03368,-.0611){6}{\line(0,-1){.0611}}
\multiput(70.19,46.38)(-.03016,-.06291){6}{\line(0,-1){.06291}}
\multiput(70.01,46)(-.03185,-.07742){5}{\line(0,-1){.07742}}
\multiput(69.85,45.61)(-.02741,-.0791){5}{\line(0,-1){.0791}}
\put(69.71,45.22){\line(0,-1){.403}}
\put(69.6,44.81){\line(0,-1){.408}}
\put(69.51,44.41){\line(0,-1){.413}}
\put(69.44,43.99){\line(0,-1){.416}}
\put(69.39,43.58){\line(0,-1){1.67}}
\put(69.45,41.91){\line(0,-1){.412}}
\put(69.53,41.5){\line(0,-1){.407}}
\put(69.62,41.09){\line(0,-1){.401}}
\multiput(69.74,40.69)(.02848,-.07872){5}{\line(0,-1){.07872}}
\multiput(69.88,40.29)(.0329,-.07698){5}{\line(0,-1){.07698}}
\multiput(70.05,39.91)(.03101,-.06249){6}{\line(0,-1){.06249}}
\multiput(70.23,39.53)(.02958,-.05197){7}{\line(0,-1){.05197}}
\multiput(70.44,39.17)(.03248,-.05021){7}{\line(0,-1){.05021}}
\multiput(70.67,38.82)(.03086,-.04225){8}{\line(0,-1){.04225}}
\multiput(70.92,38.48)(.03321,-.04044){8}{\line(0,-1){.04044}}
\multiput(71.18,38.16)(.03151,-.03421){9}{\line(0,-1){.03421}}
\multiput(71.47,37.85)(.0334,-.03237){9}{\line(1,0){.0334}}
\multiput(71.77,37.56)(.03518,-.03043){9}{\line(1,0){.03518}}
\multiput(72.08,37.28)(.04145,-.03193){8}{\line(1,0){.04145}}
\multiput(72.41,37.03)(.0432,-.02953){8}{\line(1,0){.0432}}
\multiput(72.76,36.79)(.0512,-.03089){7}{\line(1,0){.0512}}
\multiput(73.12,36.58)(.06168,-.0326){6}{\line(1,0){.06168}}
\multiput(73.49,36.38)(.06343,-.02905){6}{\line(1,0){.06343}}
\multiput(73.87,36.21)(.07797,-.03048){5}{\line(1,0){.07797}}
\multiput(74.26,36.05)(.0995,-.0325){4}{\line(1,0){.0995}}
\put(74.66,35.92){\line(1,0){.405}}
\put(75.06,35.82){\line(1,0){.41}}
\put(75.47,35.73){\line(1,0){.414}}
\put(75.89,35.67){\line(1,0){.417}}
\put(76.3,35.63){\line(1,0){.837}}
\put(77.14,35.63){\line(1,0){.417}}
\put(77.56,35.66){\line(1,0){.415}}
\put(77.97,35.72){\line(1,0){.411}}
\put(78.38,35.8){\line(1,0){.405}}
\multiput(78.79,35.91)(.0997,.0317){4}{\line(1,0){.0997}}
\multiput(79.19,36.03)(.07821,.02986){5}{\line(1,0){.07821}}
\multiput(79.58,36.18)(.06366,.02854){6}{\line(1,0){.06366}}
\multiput(79.96,36.35)(.06194,.03211){6}{\line(1,0){.06194}}
\multiput(80.33,36.55)(.05145,.03049){7}{\line(1,0){.05145}}
\multiput(80.69,36.76)(.04963,.03335){7}{\line(1,0){.04963}}
\multiput(81.04,36.99)(.04171,.0316){8}{\line(1,0){.04171}}
\multiput(81.37,37.25)(.03542,.03015){9}{\line(1,0){.03542}}
\multiput(81.69,37.52)(.03365,.0321){9}{\line(1,0){.03365}}
\multiput(81.99,37.81)(.03178,.03396){9}{\line(0,1){.03396}}
\multiput(82.28,38.11)(.03353,.04017){8}{\line(0,1){.04017}}
\multiput(82.55,38.43)(.0312,.04201){8}{\line(0,1){.04201}}
\multiput(82.8,38.77)(.03287,.04995){7}{\line(0,1){.04995}}
\multiput(83.03,39.12)(.02999,.05174){7}{\line(0,1){.05174}}
\multiput(83.24,39.48)(.03151,.06225){6}{\line(0,1){.06225}}
\multiput(83.43,39.86)(.03351,.07672){5}{\line(0,1){.07672}}
\multiput(83.59,40.24)(.02911,.0785){5}{\line(0,1){.0785}}
\multiput(83.74,40.63)(.0308,.1){4}{\line(0,1){.1}}
\put(83.86,41.03){\line(0,1){.406}}
\put(83.96,41.44){\line(0,1){.411}}
\put(84.04,41.85){\line(0,1){.415}}
\put(84.09,42.26){\line(0,1){.736}}
\put(108.88,43){\line(0,1){.418}}
\put(108.87,43.42){\line(0,1){.417}}
\put(108.83,43.84){\line(0,1){.414}}
\put(108.77,44.25){\line(0,1){.41}}
\put(108.69,44.66){\line(0,1){.405}}
\multiput(108.59,45.07)(-.0321,.0996){4}{\line(0,1){.0996}}
\multiput(108.46,45.46)(-.03018,.07809){5}{\line(0,1){.07809}}
\multiput(108.31,45.85)(-.0288,.06355){6}{\line(0,1){.06355}}
\multiput(108.13,46.24)(-.03235,.06181){6}{\line(0,1){.06181}}
\multiput(107.94,46.61)(-.03069,.05132){7}{\line(0,1){.05132}}
\multiput(107.72,46.97)(-.03355,.0495){7}{\line(0,1){.0495}}
\multiput(107.49,47.31)(-.03176,.04158){8}{\line(0,1){.04158}}
\multiput(107.24,47.64)(-.03029,.0353){9}{\line(0,1){.0353}}
\multiput(106.96,47.96)(-.03224,.03353){9}{\line(0,1){.03353}}
\multiput(106.67,48.26)(-.03409,.03164){9}{\line(-1,0){.03409}}
\multiput(106.37,48.55)(-.0403,.03337){8}{\line(-1,0){.0403}}
\multiput(106.04,48.82)(-.04213,.03103){8}{\line(-1,0){.04213}}
\multiput(105.71,49.06)(-.05008,.03267){7}{\line(-1,0){.05008}}
\multiput(105.36,49.29)(-.05186,.02978){7}{\line(-1,0){.05186}}
\multiput(104.99,49.5)(-.06237,.03126){6}{\line(-1,0){.06237}}
\multiput(104.62,49.69)(-.07685,.03321){5}{\line(-1,0){.07685}}
\multiput(104.23,49.85)(-.07861,.0288){5}{\line(-1,0){.07861}}
\multiput(103.84,50)(-.1001,.0304){4}{\line(-1,0){.1001}}
\put(103.44,50.12){\line(-1,0){.407}}
\put(103.03,50.22){\line(-1,0){.412}}
\put(102.62,50.29){\line(-1,0){.415}}
\put(102.21,50.35){\line(-1,0){.418}}
\put(101.79,50.37){\line(-1,0){.419}}
\put(101.37,50.38){\line(-1,0){.418}}
\put(100.95,50.36){\line(-1,0){.416}}
\put(100.54,50.32){\line(-1,0){.413}}
\put(100.12,50.25){\line(-1,0){.409}}
\put(99.71,50.16){\line(-1,0){.403}}
\multiput(99.31,50.05)(-.07921,-.0271){5}{\line(-1,0){.07921}}
\multiput(98.91,49.91)(-.07755,-.03154){5}{\line(-1,0){.07755}}
\multiput(98.53,49.75)(-.06303,-.02991){6}{\line(-1,0){.06303}}
\multiput(98.15,49.58)(-.06123,-.03343){6}{\line(-1,0){.06123}}
\multiput(97.78,49.37)(-.05078,-.03159){7}{\line(-1,0){.05078}}
\multiput(97.43,49.15)(-.04279,-.03011){8}{\line(-1,0){.04279}}
\multiput(97.08,48.91)(-.04101,-.03249){8}{\line(-1,0){.04101}}
\multiput(96.76,48.65)(-.03476,-.0309){9}{\line(-1,0){.03476}}
\multiput(96.44,48.37)(-.03295,-.03282){9}{\line(-1,0){.03295}}
\multiput(96.15,48.08)(-.03104,-.03464){9}{\line(0,-1){.03464}}
\multiput(95.87,47.77)(-.03265,-.04089){8}{\line(0,-1){.04089}}
\multiput(95.61,47.44)(-.03028,-.04267){8}{\line(0,-1){.04267}}
\multiput(95.36,47.1)(-.03179,-.05065){7}{\line(0,-1){.05065}}
\multiput(95.14,46.74)(-.03368,-.0611){6}{\line(0,-1){.0611}}
\multiput(94.94,46.38)(-.03016,-.06291){6}{\line(0,-1){.06291}}
\multiput(94.76,46)(-.03185,-.07742){5}{\line(0,-1){.07742}}
\multiput(94.6,45.61)(-.02741,-.0791){5}{\line(0,-1){.0791}}
\put(94.46,45.22){\line(0,-1){.403}}
\put(94.35,44.81){\line(0,-1){.408}}
\put(94.26,44.41){\line(0,-1){.413}}
\put(94.19,43.99){\line(0,-1){.416}}
\put(94.14,43.58){\line(0,-1){1.67}}
\put(94.2,41.91){\line(0,-1){.412}}
\put(94.28,41.5){\line(0,-1){.407}}
\put(94.37,41.09){\line(0,-1){.401}}
\multiput(94.49,40.69)(.02848,-.07872){5}{\line(0,-1){.07872}}
\multiput(94.63,40.29)(.0329,-.07698){5}{\line(0,-1){.07698}}
\multiput(94.8,39.91)(.03101,-.06249){6}{\line(0,-1){.06249}}
\multiput(94.98,39.53)(.02958,-.05197){7}{\line(0,-1){.05197}}
\multiput(95.19,39.17)(.03248,-.05021){7}{\line(0,-1){.05021}}
\multiput(95.42,38.82)(.03086,-.04225){8}{\line(0,-1){.04225}}
\multiput(95.67,38.48)(.03321,-.04044){8}{\line(0,-1){.04044}}
\multiput(95.93,38.16)(.03151,-.03421){9}{\line(0,-1){.03421}}
\multiput(96.22,37.85)(.0334,-.03237){9}{\line(1,0){.0334}}
\multiput(96.52,37.56)(.03518,-.03043){9}{\line(1,0){.03518}}
\multiput(96.83,37.28)(.04145,-.03193){8}{\line(1,0){.04145}}
\multiput(97.16,37.03)(.0432,-.02953){8}{\line(1,0){.0432}}
\multiput(97.51,36.79)(.0512,-.03089){7}{\line(1,0){.0512}}
\multiput(97.87,36.58)(.06168,-.0326){6}{\line(1,0){.06168}}
\multiput(98.24,36.38)(.06343,-.02905){6}{\line(1,0){.06343}}
\multiput(98.62,36.21)(.07797,-.03048){5}{\line(1,0){.07797}}
\multiput(99.01,36.05)(.0995,-.0325){4}{\line(1,0){.0995}}
\put(99.41,35.92){\line(1,0){.405}}
\put(99.81,35.82){\line(1,0){.41}}
\put(100.22,35.73){\line(1,0){.414}}
\put(100.64,35.67){\line(1,0){.417}}
\put(101.05,35.63){\line(1,0){.837}}
\put(101.89,35.63){\line(1,0){.417}}
\put(102.31,35.66){\line(1,0){.415}}
\put(102.72,35.72){\line(1,0){.411}}
\put(103.13,35.8){\line(1,0){.405}}
\multiput(103.54,35.91)(.0997,.0317){4}{\line(1,0){.0997}}
\multiput(103.94,36.03)(.07821,.02986){5}{\line(1,0){.07821}}
\multiput(104.33,36.18)(.06366,.02854){6}{\line(1,0){.06366}}
\multiput(104.71,36.35)(.06194,.03211){6}{\line(1,0){.06194}}
\multiput(105.08,36.55)(.05145,.03049){7}{\line(1,0){.05145}}
\multiput(105.44,36.76)(.04963,.03335){7}{\line(1,0){.04963}}
\multiput(105.79,36.99)(.04171,.0316){8}{\line(1,0){.04171}}
\multiput(106.12,37.25)(.03542,.03015){9}{\line(1,0){.03542}}
\multiput(106.44,37.52)(.03365,.0321){9}{\line(1,0){.03365}}
\multiput(106.74,37.81)(.03178,.03396){9}{\line(0,1){.03396}}
\multiput(107.03,38.11)(.03353,.04017){8}{\line(0,1){.04017}}
\multiput(107.3,38.43)(.0312,.04201){8}{\line(0,1){.04201}}
\multiput(107.55,38.77)(.03287,.04995){7}{\line(0,1){.04995}}
\multiput(107.78,39.12)(.02999,.05174){7}{\line(0,1){.05174}}
\multiput(107.99,39.48)(.03151,.06225){6}{\line(0,1){.06225}}
\multiput(108.18,39.86)(.03351,.07672){5}{\line(0,1){.07672}}
\multiput(108.34,40.24)(.02911,.0785){5}{\line(0,1){.0785}}
\multiput(108.49,40.63)(.0308,.1){4}{\line(0,1){.1}}
\put(108.61,41.03){\line(0,1){.406}}
\put(108.71,41.44){\line(0,1){.411}}
\put(108.79,41.85){\line(0,1){.415}}
\put(108.84,42.26){\line(0,1){.736}}
\put(104.63,45.13){\oval(149.75,44.75)[]}
\put(105,45.5){\oval(165.5,63)[]}
\qbezier(57.25,49)(117.75,75.75)(118.25,44.5)
\qbezier(118.25,44.5)(118.88,21.63)(56,36.25)
\put(52,43.5){\makebox(0,0)[cc]{2}}
\put(102,43){\makebox(0,0)[cc]{3}}
\put(50.75,65.25){\makebox(0,0)[cc]{6}}
\put(56.5,80.5){\makebox(0,0)[cc]{7}}
\put(103.25,59){\makebox(0,0)[cc]{0}}
\put(67,74){\vector(0,-1){.07}}\put(67,79){\line(0,-1){5}}
\put(66,64.75){\vector(0,-1){.07}}\put(66,70.25){\line(0,-1){5.5}}
\put(76.75,42.5){\makebox(0,0)[cc]{6}}
\put(60.25,43.25){\line(1,0){9.25}}
\put(84.25,43){\line(1,0){9.75}}
\put(64,46.75){\makebox(0,0)[cc]{4}}
\put(88.75,47){\makebox(0,0)[cc]{5}}
\put(86.25,43){\line(0,1){.53}} \put(86.24,43.53){\line(0,1){.529}}
\put(86.19,44.06){\line(0,1){.525}}
\put(86.12,44.58){\line(0,1){.521}}
\multiput(86.02,45.1)(-.0321,.1286){4}{\line(0,1){.1286}}
\multiput(85.89,45.62)(-.03126,.10135){5}{\line(0,1){.10135}}
\multiput(85.74,46.13)(-.03061,.08292){6}{\line(0,1){.08292}}
\multiput(85.55,46.62)(-.03006,.06954){7}{\line(0,1){.06954}}
\multiput(85.34,47.11)(-.02957,.05933){8}{\line(0,1){.05933}}
\multiput(85.1,47.59)(-.03275,.05763){8}{\line(0,1){.05763}}
\multiput(84.84,48.05)(-.03186,.04957){9}{\line(0,1){.04957}}
\multiput(84.56,48.49)(-.03106,.04299){10}{\line(0,1){.04299}}
\multiput(84.25,48.92)(-.03335,.04124){10}{\line(0,1){.04124}}
\multiput(83.91,49.33)(-.03231,.03578){11}{\line(0,1){.03578}}
\multiput(83.56,49.73)(-.03136,.03114){12}{\line(-1,0){.03136}}
\multiput(83.18,50.1)(-.036,.03206){11}{\line(-1,0){.036}}
\multiput(82.78,50.45)(-.04146,.03307){10}{\line(-1,0){.04146}}
\multiput(82.37,50.79)(-.0432,.03076){10}{\line(-1,0){.0432}}
\multiput(81.94,51.09)(-.04979,.03152){9}{\line(-1,0){.04979}}
\multiput(81.49,51.38)(-.05786,.03236){8}{\line(-1,0){.05786}}
\multiput(81.03,51.64)(-.06803,.03334){7}{\line(-1,0){.06803}}
\multiput(80.55,51.87)(-.06975,.02959){7}{\line(-1,0){.06975}}
\multiput(80.06,52.08)(-.08313,.03004){6}{\line(-1,0){.08313}}
\multiput(79.56,52.26)(-.10156,.03057){5}{\line(-1,0){.10156}}
\multiput(79.06,52.41)(-.1288,.0313){4}{\line(-1,0){.1288}}
\put(78.54,52.53){\line(-1,0){.521}}
\put(78.02,52.63){\line(-1,0){.526}}
\put(77.49,52.7){\line(-1,0){.529}}
\put(76.96,52.74){\line(-1,0){1.06}}
\put(75.9,52.73){\line(-1,0){.528}}
\put(75.38,52.68){\line(-1,0){.525}}
\put(74.85,52.61){\line(-1,0){.52}}
\multiput(74.33,52.51)(-.1284,-.033){4}{\line(-1,0){.1284}}
\multiput(73.82,52.37)(-.10114,-.03195){5}{\line(-1,0){.10114}}
\multiput(73.31,52.21)(-.08271,-.03117){6}{\line(-1,0){.08271}}
\multiput(72.81,52.03)(-.06934,-.03053){7}{\line(-1,0){.06934}}
\multiput(72.33,51.81)(-.05913,-.02998){8}{\line(-1,0){.05913}}
\multiput(71.86,51.57)(-.05741,-.03315){8}{\line(-1,0){.05741}}
\multiput(71.4,51.31)(-.04935,-.03219){9}{\line(-1,0){.04935}}
\multiput(70.95,51.02)(-.04278,-.03135){10}{\line(-1,0){.04278}}
\multiput(70.53,50.7)(-.04101,-.03363){10}{\line(-1,0){.04101}}
\multiput(70.12,50.37)(-.03556,-.03255){11}{\line(-1,0){.03556}}
\multiput(69.72,50.01)(-.03093,-.03157){12}{\line(0,-1){.03157}}
\multiput(69.35,49.63)(-.03182,-.03622){11}{\line(0,-1){.03622}}
\multiput(69,49.23)(-.03278,-.04169){10}{\line(0,-1){.04169}}
\multiput(68.67,48.82)(-.03047,-.04341){10}{\line(0,-1){.04341}}
\multiput(68.37,48.38)(-.03118,-.05){9}{\line(0,-1){.05}}
\multiput(68.09,47.93)(-.03197,-.05807){8}{\line(0,-1){.05807}}
\multiput(67.83,47.47)(-.03287,-.06826){7}{\line(0,-1){.06826}}
\multiput(67.6,46.99)(-.02911,-.06995){7}{\line(0,-1){.06995}}
\multiput(67.4,46.5)(-.02948,-.08333){6}{\line(0,-1){.08333}}
\multiput(67.22,46)(-.02988,-.10177){5}{\line(0,-1){.10177}}
\multiput(67.07,45.49)(-.0304,-.1291){4}{\line(0,-1){.1291}}
\put(66.95,44.98){\line(0,-1){.522}}
\put(66.86,44.45){\line(0,-1){.526}}
\put(66.79,43.93){\line(0,-1){1.059}}
\put(66.75,42.87){\line(0,-1){.53}}
\put(66.77,42.34){\line(0,-1){.528}}
\put(66.82,41.81){\line(0,-1){.524}}
\put(66.9,41.29){\line(0,-1){.519}}
\multiput(67.01,40.77)(.0271,-.10255){5}{\line(0,-1){.10255}}
\multiput(67.14,40.25)(.03264,-.10092){5}{\line(0,-1){.10092}}
\multiput(67.31,39.75)(.03173,-.0825){6}{\line(0,-1){.0825}}
\multiput(67.5,39.25)(.031,-.06913){7}{\line(0,-1){.06913}}
\multiput(67.72,38.77)(.03038,-.05892){8}{\line(0,-1){.05892}}
\multiput(67.96,38.3)(.03354,-.05718){8}{\line(0,-1){.05718}}
\multiput(68.23,37.84)(.03253,-.04913){9}{\line(0,-1){.04913}}
\multiput(68.52,37.4)(.03164,-.04256){10}{\line(0,-1){.04256}}
\multiput(68.84,36.97)(.03082,-.03707){11}{\line(0,-1){.03707}}
\multiput(69.18,36.57)(.03279,-.03534){11}{\line(0,-1){.03534}}
\multiput(69.54,36.18)(.03467,-.03351){11}{\line(1,0){.03467}}
\multiput(69.92,35.81)(.03644,-.03157){11}{\line(1,0){.03644}}
\multiput(70.32,35.46)(.04191,-.0325){10}{\line(1,0){.04191}}
\multiput(70.74,35.14)(.04846,-.03352){9}{\line(1,0){.04846}}
\multiput(71.17,34.83)(.05021,-.03084){9}{\line(1,0){.05021}}
\multiput(71.63,34.56)(.05829,-.03157){8}{\line(1,0){.05829}}
\multiput(72.09,34.3)(.06848,-.03241){7}{\line(1,0){.06848}}
\multiput(72.57,34.08)(.08183,-.03341){6}{\line(1,0){.08183}}
\multiput(73.06,33.88)(.08353,-.02891){6}{\line(1,0){.08353}}
\multiput(73.56,33.7)(.10197,-.02919){5}{\line(1,0){.10197}}
\put(74.07,33.56){\line(1,0){.517}}
\put(74.59,33.44){\line(1,0){.523}}
\put(75.11,33.35){\line(1,0){.527}}
\put(75.64,33.29){\line(1,0){.529}}
\put(76.17,33.26){\line(1,0){.53}} \put(76.7,33.25){\line(1,0){.53}}
\put(77.23,33.28){\line(1,0){.528}}
\put(77.76,33.33){\line(1,0){.524}}
\put(78.28,33.41){\line(1,0){.519}}
\multiput(78.8,33.52)(.10236,.0278){5}{\line(1,0){.10236}}
\multiput(79.31,33.66)(.1007,.03332){5}{\line(1,0){.1007}}
\multiput(79.81,33.83)(.08228,.03229){6}{\line(1,0){.08228}}
\multiput(80.31,34.02)(.06892,.03147){7}{\line(1,0){.06892}}
\multiput(80.79,34.24)(.05871,.03078){8}{\line(1,0){.05871}}
\multiput(81.26,34.49)(.05063,.03015){9}{\line(1,0){.05063}}
\multiput(81.72,34.76)(.04891,.03286){9}{\line(1,0){.04891}}
\multiput(82.16,35.06)(.04235,.03193){10}{\line(1,0){.04235}}
\multiput(82.58,35.38)(.03686,.03107){11}{\line(1,0){.03686}}
\multiput(82.98,35.72)(.03512,.03303){11}{\line(1,0){.03512}}
\multiput(83.37,36.08)(.03327,.03489){11}{\line(0,1){.03489}}
\multiput(83.74,36.47)(.03132,.03665){11}{\line(0,1){.03665}}
\multiput(84.08,36.87)(.03221,.04213){10}{\line(0,1){.04213}}
\multiput(84.4,37.29)(.03319,.04869){9}{\line(0,1){.04869}}
\multiput(84.7,37.73)(.0305,.05042){9}{\line(0,1){.05042}}
\multiput(84.98,38.18)(.03118,.0585){8}{\line(0,1){.0585}}
\multiput(85.23,38.65)(.03194,.0687){7}{\line(0,1){.0687}}
\multiput(85.45,39.13)(.03285,.08206){6}{\line(0,1){.08206}}
\multiput(85.65,39.62)(.02834,.08372){6}{\line(0,1){.08372}}
\multiput(85.82,40.13)(.0285,.10217){5}{\line(0,1){.10217}}
\put(85.96,40.64){\line(0,1){.518}}
\put(86.07,41.15){\line(0,1){.523}}
\put(86.16,41.68){\line(0,1){.527}}
\put(86.22,42.21){\line(0,1){.795}}
\put(86,37.5){\makebox(0,0)[cc]{7}}
\put(152.5,42.75){\makebox(0,0)[cc]{3}}
\put(109,43){\line(1,0){36.5}}
\put(160.25,43.5){\line(1,0){27.25}}
\multiput(22.75,43.5)(-.033333,.05){15}{\line(0,1){.05}}
\put(45.25,43.5){\line(-1,0){23.25}}
\put(34.25,46.25){\makebox(0,0)[cc]{4}}
\put(127.5,45.5){\makebox(0,0)[cc]{1}}
\put(170.25,47){\makebox(0,0)[cc]{5}}
\put(38.25,40){\vector(0,-1){.07}}\put(38.25,46.5){\line(0,-1){6.5}}
\put(63.75,46.25){\vector(0,1){.07}}\put(63.75,40.5){\line(0,1){5.75}}
\put(90.75,38){\vector(0,-1){.07}}\put(90.75,45.75){\line(0,-1){7.75}}
\put(78,48.5){\vector(-1,-3){.07}}\multiput(80.25,55)(-.0335821,-.0970149){67}{\line(0,-1){.0970149}}
\put(165.75,46.5){\vector(0,1){.07}}\put(165.75,41){\line(0,1){5.5}}
\put(102.25,32.75){\vector(0,-1){.07}}\put(102.25,37.75){\line(0,-1){5}}
\put(153.25,38.5){\vector(0,1){.07}}\put(153.25,33.5){\line(0,1){5}}
\end{picture}
\end{minipage}
\end{center}

\proof Let $X=\emptyset$. The order on $P_2(8)$ is $2\succ
1\succ0\succ3\succ4\succ5\succ6\succ7$. The partial orders of
anti-links:
\begin{itemize}
\item 0: $3\succ_04\succ_05\succ_06$;
\item 1: $4\succ_17\succ_15$;
\item 2: $1\succ_23\succ_25\succ_27$;
\item 3: $2\succ_34\succ_36$;
\item 4: $1\succ_43\succ_45$, $1\succ_4 0$;
\item 5: $1\succ_5 0$.
\end{itemize}

We need to show that the conditions of the theorem hold for our
order. Let $p=(v_1,...,v_k)$ be any anti-path in $P_2(8)$. Then by
the definition of $\Theta(p)$, it must contain $0$. Assume that
$\Theta(p)$ does not satisfy (*). Since $2\succ i$ for every $i\ne
2$, we have $2\in \Theta(p)$.

We consider eight different possibilities for $v_k$.

\medskip

{\bf 0. }Let { ${v_k=0}$}. Then $\Theta(p)\supseteq \{0,2, 3, 4, 5,
6\}$. This set satisfies (*).

{\bf 1.} Let $v_k=1$. Then $\Theta(p)$ contains $\{1, 2, 4, 5, 7\}$
that satisfies (*).

{\bf 2.} Let $v_k=2$. Then $\Theta(p)$ contains $\{1,2,3,5,7\}$ that
satisfies (*).

{\bf 3.} Let $v_k=3$. Then $\Theta(p)$ contains $\{0,2,3,4,6\}$ that
satisfies (*).

{\bf 4.} Let $v_k=4$. Then $\Theta(p)$ contains $\{0,1,2,3,4,5\}$
which satisfies (*)

{\bf 5.} Let $v_k=5$. Then $\Theta(p)$ contains $\{0,1,2,4,5\}$ (the
anti-star of 5). The set $\{0,1,2,3,4,5\}$ satisfies (*). So
$3\not\in \Theta(p)$.

Then $k\ne 1$ since $3\succ5$, $v_{k-1}\in\{0, 1, 2, 4\}$.

Since $3\succ_0 5$, $3\succ_2 5$ and $3\succ_4 5$, we have
$v_{k-1}=1$. Since $7\succ_1 5$, we have that $\Theta(p)$ contains
$\{0,1,2,4,5,7\}$ which satisfies (*).

{\bf 6.} Let $v_k=6$. Then $\Theta(p)$ contains $\{0,2,3,6\}$. We
can assume that $k\ne 1$ since otherwise $\Theta(p)$ would contain
0, 1, 2, 3, 4 , 5, 6 and would satisfy (*). The vertex $v_{k-1}$ is
in $\{0, 3\}$.

Suppose that $v_{k-1}=0$. Then $\Theta(p)$ contains
$\{0,2,3,4,5,6\}$ since $4, 5\succ_0 6$. Since $1\succ 0$, and
$\{0,1,2,3,4,5,6\}$ satisfies (*), we have $k-1\ne 1$ and
$v_{k-2}\in \{3,4,5\}$. Since $1\succ_50$, $1\succ_40$, $v_{k-2}=3$.
Since $1\succ 3$, $k-2\ne 1$, $v_{k-3}\in \{2,4\}$. But $1\succ_23$,
$1\succ_43$, a contradiction.

If $v_{k-1}=3$ then $\Theta(p)$ contains $\{0, 2, 3, 4, 6\}$ since
$2, 4\succ_3 6$.

Since $1\succ 3$, and $\{0,1,2,3,4,6\}$ satisfies (*), $k-1\ne 1$,
$v_{k-2}\in \{0,2,4\}$. Since $1\succ_2 3$, $1\succ_4 3$, we have
$v_{k-2}=0$. Since $1\succ 0$, we have $k-2\ne 1$, $v_{k-3}\in
\{4,5\}$. Since $1\succ_4 0, 1\succ_5 0$, we get a contradiction.

{\bf 7.} Let $v_k=7$. Then $\Theta(p)$ contains $\{0,1,2,7\}$. We
have $k\ne 1$. The vertex $v_{k-1}$ can be $1$ or $2$.

If $v_{k-1}=1$ then $\Theta(p)$ contains $\{0,1,2,4,7\}$ since $4
\succ_1 7$. That set satisfies (*).

Let $v_{k-1}=2$. Then $\Theta(p)$ contains $\{0,1,2,3,5,7\}$ since
$3,5\succ_2 7$.  That set satisfies (*).\endproof

\subsubsection{The 8-vertex graph $P_3(8)$}

\begin{center}
\unitlength .75mm 
\linethickness{0.4pt}
\ifx\plotpoint\undefined\newsavebox{\plotpoint}\fi 
\begin{picture}(199,100.75)(0,0)
\put(7.5,23.75){\line(1,0){63.25}}
\multiput(70.75,23.75)(-.051835853,.058855292){463}{\line(0,1){.058855292}}
\put(46.75,51){\line(0,-1){27}}
\multiput(46.75,24)(-.05304878,-.051829268){410}{\line(-1,0){.05304878}}
\put(25,2.75){\line(0,1){21}}
\multiput(25,23.75)(.051807229,.065060241){415}{\line(0,1){.065060241}}
\multiput(7.5,23.75)(.051801802,-.061561562){333}{\line(0,-1){.061561562}}
\put(46.75,35.5){\line(2,-1){23.5}}
\qbezier(70.5,24)(42.5,106.75)(7.5,23.5)
\put(4.5,23.25){\makebox(0,0)[cc]{0}}
\put(73.25,23.5){\makebox(0,0)[cc]{1}}
\put(28.75,3.25){\makebox(0,0)[cc]{2}}
\put(22.25,14){\makebox(0,0)[cc]{4}}
\put(22.75,27.25){\makebox(0,0)[cc]{7}}
\put(46.75,21.25){\makebox(0,0)[cc]{6}}
\put(43,35.5){\makebox(0,0)[cc]{5}}
\put(47,54.25){\makebox(0,0)[cc]{3}}
\put(10.25,52.5){\makebox(0,0)[cc]{$P_3(8)$}}
\put(7.25,23.75){\vector(-2,1){.11}}\multiput(25,14)(-.09441489,.0518617){188}{\line(-1,0){.09441489}}
\put(116,57){\circle{11.85}} \put(142.75,24){\circle{11.85}}
\put(169.25,57){\circle{11.85}}
\put(121.75,57.25){\line(1,0){41.5}}
\multiput(165.25,52.75)(-.051771117,-.066076294){367}{\line(0,-1){.066076294}}
\multiput(118.5,51.75)(.051886792,-.064690027){371}{\line(0,-1){.064690027}}
\put(187.51,47.75){\line(0,1){1.546}}
\put(187.48,49.3){\line(0,1){1.544}}
\multiput(187.41,50.84)(-.0428,.5135){3}{\line(0,1){.5135}}
\multiput(187.28,52.38)(-.0449,.3839){4}{\line(0,1){.3839}}
\multiput(187.1,53.92)(-.04608,.30574){5}{\line(0,1){.30574}}
\multiput(186.87,55.44)(-.04684,.25336){6}{\line(0,1){.25336}}
\multiput(186.59,56.96)(-.04735,.21571){7}{\line(0,1){.21571}}
\multiput(186.26,58.47)(-.04768,.18727){8}{\line(0,1){.18727}}
\multiput(185.87,59.97)(-.04789,.16496){9}{\line(0,1){.16496}}
\multiput(185.44,61.46)(-.04801,.14695){10}{\line(0,1){.14695}}
\multiput(184.96,62.93)(-.04806,.13207){11}{\line(0,1){.13207}}
\multiput(184.43,64.38)(-.04806,.11953){12}{\line(0,1){.11953}}
\multiput(183.86,65.81)(-.048003,.108799){13}{\line(0,1){.108799}}
\multiput(183.23,67.23)(-.051592,.107144){13}{\line(0,1){.107144}}
\multiput(182.56,68.62)(-.051187,.097844){14}{\line(0,1){.097844}}
\multiput(181.85,69.99)(-.050783,.089682){15}{\line(0,1){.089682}}
\multiput(181.08,71.34)(-.050378,.082449){16}{\line(0,1){.082449}}
\multiput(180.28,72.66)(-.049967,.07598){17}{\line(0,1){.07598}}
\multiput(179.43,73.95)(-.04955,.070151){18}{\line(0,1){.070151}}
\multiput(178.54,75.21)(-.051854,.068465){18}{\line(0,1){.068465}}
\multiput(177.6,76.44)(-.051253,.063193){19}{\line(0,1){.063193}}
\multiput(176.63,77.64)(-.050659,.058382){20}{\line(0,1){.058382}}
\multiput(175.62,78.81)(-.050068,.053968){21}{\line(0,1){.053968}}
\multiput(174.57,79.94)(-.051834,.052274){21}{\line(0,1){.052274}}
\multiput(173.48,81.04)(-.053543,.050523){21}{\line(-1,0){.053543}}
\multiput(172.35,82.1)(-.057952,.051151){20}{\line(-1,0){.057952}}
\multiput(171.19,83.13)(-.062758,.051786){19}{\line(-1,0){.062758}}
\multiput(170,84.11)(-.064444,.049671){19}{\line(-1,0){.064444}}
\multiput(168.78,85.05)(-.069729,.050141){18}{\line(-1,0){.069729}}
\multiput(167.52,85.96)(-.075555,.050608){17}{\line(-1,0){.075555}}
\multiput(166.24,86.82)(-.08202,.051073){16}{\line(-1,0){.08202}}
\multiput(164.92,87.63)(-.08925,.05154){15}{\line(-1,0){.08925}}
\multiput(163.59,88.41)(-.090913,.048545){15}{\line(-1,0){.090913}}
\multiput(162.22,89.13)(-.099082,.048747){14}{\line(-1,0){.099082}}
\multiput(160.83,89.82)(-.10839,.048921){13}{\line(-1,0){.10839}}
\multiput(159.43,90.45)(-.11912,.04907){12}{\line(-1,0){.11912}}
\multiput(158,91.04)(-.13165,.04918){11}{\line(-1,0){.13165}}
\multiput(156.55,91.58)(-.14654,.04925){10}{\line(-1,0){.14654}}
\multiput(155.08,92.08)(-.16455,.04928){9}{\line(-1,0){.16455}}
\multiput(153.6,92.52)(-.18686,.04926){8}{\line(-1,0){.18686}}
\multiput(152.11,92.91)(-.21531,.04917){7}{\line(-1,0){.21531}}
\multiput(150.6,93.26)(-.25296,.04899){6}{\line(-1,0){.25296}}
\multiput(149.08,93.55)(-.30533,.04866){5}{\line(-1,0){.30533}}
\multiput(147.56,93.79)(-.3835,.0481){4}{\line(-1,0){.3835}}
\multiput(146.02,93.99)(-.5132,.0471){3}{\line(-1,0){.5132}}
\put(144.48,94.13){\line(-1,0){1.543}}
\put(142.94,94.22){\line(-1,0){1.545}}
\put(141.39,94.26){\line(-1,0){1.546}}
\put(139.85,94.24){\line(-1,0){1.545}}
\put(138.3,94.18){\line(-1,0){1.542}}
\multiput(136.76,94.07)(-.3842,-.0416){4}{\line(-1,0){.3842}}
\multiput(135.22,93.9)(-.30611,-.04349){5}{\line(-1,0){.30611}}
\multiput(133.69,93.68)(-.25375,-.0447){6}{\line(-1,0){.25375}}
\multiput(132.17,93.41)(-.21611,-.04552){7}{\line(-1,0){.21611}}
\multiput(130.66,93.09)(-.18766,-.04609){8}{\line(-1,0){.18766}}
\multiput(129.16,92.73)(-.16536,-.04649){9}{\line(-1,0){.16536}}
\multiput(127.67,92.31)(-.14735,-.04677){10}{\line(-1,0){.14735}}
\multiput(126.2,91.84)(-.14571,-.05164){10}{\line(-1,0){.14571}}
\multiput(124.74,91.32)(-.13083,-.05132){11}{\line(-1,0){.13083}}
\multiput(123.3,90.76)(-.1183,-.051){12}{\line(-1,0){.1183}}
\multiput(121.88,90.15)(-.107576,-.050685){13}{\line(-1,0){.107576}}
\multiput(120.48,89.49)(-.098273,-.050358){14}{\line(-1,0){.098273}}
\multiput(119.11,88.78)(-.090109,-.050023){15}{\line(-1,0){.090109}}
\multiput(117.75,88.03)(-.082872,-.049679){16}{\line(-1,0){.082872}}
\multiput(116.43,87.24)(-.0764,-.049323){17}{\line(-1,0){.0764}}
\multiput(115.13,86.4)(-.074718,-.051835){17}{\line(-1,0){.074718}}
\multiput(113.86,85.52)(-.068901,-.051273){18}{\line(-1,0){.068901}}
\multiput(112.62,84.59)(-.063624,-.050717){19}{\line(-1,0){.063624}}
\multiput(111.41,83.63)(-.058808,-.050164){20}{\line(-1,0){.058808}}
\multiput(110.23,82.63)(-.054389,-.04961){21}{\line(-1,0){.054389}}
\multiput(109.09,81.59)(-.052711,-.05139){21}{\line(-1,0){.052711}}
\multiput(107.98,80.51)(-.050973,-.053114){21}{\line(0,-1){.053114}}
\multiput(106.91,79.39)(-.051639,-.057517){20}{\line(0,-1){.057517}}
\multiput(105.88,78.24)(-.049699,-.059202){20}{\line(0,-1){.059202}}
\multiput(104.89,77.06)(-.050215,-.064022){19}{\line(0,-1){.064022}}
\multiput(103.93,75.84)(-.050729,-.069303){18}{\line(0,-1){.069303}}
\multiput(103.02,74.59)(-.051245,-.075124){17}{\line(0,-1){.075124}}
\multiput(102.15,73.32)(-.051765,-.081585){16}{\line(0,-1){.081585}}
\multiput(101.32,72.01)(-.049024,-.08326){16}{\line(0,-1){.08326}}
\multiput(100.54,70.68)(-.049312,-.0905){15}{\line(0,-1){.0905}}
\multiput(99.8,69.32)(-.049583,-.098666){14}{\line(0,-1){.098666}}
\multiput(99.1,67.94)(-.049836,-.107972){13}{\line(0,-1){.107972}}
\multiput(98.45,66.54)(-.05007,-.1187){12}{\line(0,-1){.1187}}
\multiput(97.85,65.11)(-.05029,-.13123){11}{\line(0,-1){.13123}}
\multiput(97.3,63.67)(-.05049,-.14612){10}{\line(0,-1){.14612}}
\multiput(96.8,62.21)(-.05067,-.16413){9}{\line(0,-1){.16413}}
\multiput(96.34,60.73)(-.05084,-.18643){8}{\line(0,-1){.18643}}
\multiput(95.93,59.24)(-.05099,-.21488){7}{\line(0,-1){.21488}}
\multiput(95.58,57.73)(-.05112,-.25253){6}{\line(0,-1){.25253}}
\multiput(95.27,56.22)(-.05124,-.30491){5}{\line(0,-1){.30491}}
\multiput(95.01,54.69)(-.0513,-.3831){4}{\line(0,-1){.3831}}
\multiput(94.81,53.16)(-.0515,-.5127){3}{\line(0,-1){.5127}}
\put(94.65,51.62){\line(0,-1){1.542}}
\put(94.55,50.08){\line(0,-1){6.179}}
\multiput(94.65,43.9)(.0512,-.5128){3}{\line(0,-1){.5128}}
\multiput(94.8,42.36)(.0511,-.3831){4}{\line(0,-1){.3831}}
\multiput(95.01,40.83)(.05106,-.30494){5}{\line(0,-1){.30494}}
\multiput(95.26,39.31)(.05098,-.25256){6}{\line(0,-1){.25256}}
\multiput(95.57,37.79)(.05086,-.21491){7}{\line(0,-1){.21491}}
\multiput(95.93,36.29)(.05073,-.18646){8}{\line(0,-1){.18646}}
\multiput(96.33,34.8)(.05058,-.16416){9}{\line(0,-1){.16416}}
\multiput(96.79,33.32)(.05041,-.14615){10}{\line(0,-1){.14615}}
\multiput(97.29,31.86)(.05021,-.13126){11}{\line(0,-1){.13126}}
\multiput(97.84,30.41)(.05,-.11873){12}{\line(0,-1){.11873}}
\multiput(98.44,28.99)(.049773,-.108001){13}{\line(0,-1){.108001}}
\multiput(99.09,27.58)(.049525,-.098695){14}{\line(0,-1){.098695}}
\multiput(99.78,26.2)(.04926,-.090528){15}{\line(0,-1){.090528}}
\multiput(100.52,24.84)(.048976,-.083289){16}{\line(0,-1){.083289}}
\multiput(101.31,23.51)(.051717,-.081615){16}{\line(0,-1){.081615}}
\multiput(102.13,22.21)(.051201,-.075154){17}{\line(0,-1){.075154}}
\multiput(103,20.93)(.050689,-.069332){18}{\line(0,-1){.069332}}
\multiput(103.92,19.68)(.050177,-.064051){19}{\line(0,-1){.064051}}
\multiput(104.87,18.46)(.049665,-.059231){20}{\line(0,-1){.059231}}
\multiput(105.86,17.28)(.051606,-.057547){20}{\line(0,-1){.057547}}
\multiput(106.9,16.13)(.050943,-.053143){21}{\line(0,-1){.053143}}
\multiput(107.97,15.01)(.052681,-.051421){21}{\line(1,0){.052681}}
\multiput(109.07,13.93)(.054361,-.049641){21}{\line(1,0){.054361}}
\multiput(110.21,12.89)(.058779,-.050198){20}{\line(1,0){.058779}}
\multiput(111.39,11.89)(.063595,-.050754){19}{\line(1,0){.063595}}
\multiput(112.6,10.92)(.068871,-.051313){18}{\line(1,0){.068871}}
\multiput(113.84,10)(.074688,-.051878){17}{\line(1,0){.074688}}
\multiput(115.11,9.12)(.076371,-.049367){17}{\line(1,0){.076371}}
\multiput(116.4,8.28)(.082843,-.049727){16}{\line(1,0){.082843}}
\multiput(117.73,7.48)(.09008,-.050076){15}{\line(1,0){.09008}}
\multiput(119.08,6.73)(.098244,-.050415){14}{\line(1,0){.098244}}
\multiput(120.46,6.02)(.107547,-.050747){13}{\line(1,0){.107547}}
\multiput(121.85,5.36)(.11827,-.05107){12}{\line(1,0){.11827}}
\multiput(123.27,4.75)(.1308,-.0514){11}{\line(1,0){.1308}}
\multiput(124.71,4.19)(.14568,-.05172){10}{\line(1,0){.14568}}
\multiput(126.17,3.67)(.14732,-.04685){10}{\line(1,0){.14732}}
\multiput(127.64,3.2)(.16533,-.04659){9}{\line(1,0){.16533}}
\multiput(129.13,2.78)(.18764,-.0462){8}{\line(1,0){.18764}}
\multiput(130.63,2.41)(.21608,-.04565){7}{\line(1,0){.21608}}
\multiput(132.14,2.09)(.25372,-.04485){6}{\line(1,0){.25372}}
\multiput(133.67,1.82)(.30609,-.04367){5}{\line(1,0){.30609}}
\multiput(135.2,1.6)(.3842,-.0418){4}{\line(1,0){.3842}}
\put(136.73,1.44){\line(1,0){1.542}}
\put(138.28,1.32){\line(1,0){1.545}}
\put(139.82,1.26){\line(1,0){1.546}}
\put(141.37,1.24){\line(1,0){1.545}}
\put(142.91,1.28){\line(1,0){1.543}}
\multiput(144.46,1.37)(.5132,.0468){3}{\line(1,0){.5132}}
\multiput(145.99,1.51)(.3835,.0479){4}{\line(1,0){.3835}}
\multiput(147.53,1.7)(.30536,.04848){5}{\line(1,0){.30536}}
\multiput(149.06,1.94)(.25299,.04884){6}{\line(1,0){.25299}}
\multiput(150.57,2.24)(.21533,.04905){7}{\line(1,0){.21533}}
\multiput(152.08,2.58)(.18689,.04915){8}{\line(1,0){.18689}}
\multiput(153.58,2.97)(.16458,.04919){9}{\line(1,0){.16458}}
\multiput(155.06,3.42)(.14657,.04917){10}{\line(1,0){.14657}}
\multiput(156.52,3.91)(.13168,.0491){11}{\line(1,0){.13168}}
\multiput(157.97,4.45)(.11915,.049){12}{\line(1,0){.11915}}
\multiput(159.4,5.04)(.108418,.048858){13}{\line(1,0){.108418}}
\multiput(160.81,5.67)(.09911,.048689){14}{\line(1,0){.09911}}
\multiput(162.2,6.35)(.090942,.048492){15}{\line(1,0){.090942}}
\multiput(163.56,7.08)(.08928,.051488){15}{\line(1,0){.08928}}
\multiput(164.9,7.85)(.082049,.051025){16}{\line(1,0){.082049}}
\multiput(166.21,8.67)(.075584,.050564){17}{\line(1,0){.075584}}
\multiput(167.5,9.53)(.069758,.050101){18}{\line(1,0){.069758}}
\multiput(168.75,10.43)(.064473,.049634){19}{\line(1,0){.064473}}
\multiput(169.98,11.37)(.062788,.051749){19}{\line(1,0){.062788}}
\multiput(171.17,12.36)(.057982,.051117){20}{\line(1,0){.057982}}
\multiput(172.33,13.38)(.053572,.050491){21}{\line(1,0){.053572}}
\multiput(173.46,14.44)(.051864,.052244){21}{\line(0,1){.052244}}
\multiput(174.55,15.54)(.050099,.053939){21}{\line(0,1){.053939}}
\multiput(175.6,16.67)(.050693,.058353){20}{\line(0,1){.058353}}
\multiput(176.61,17.84)(.05129,.063164){19}{\line(0,1){.063164}}
\multiput(177.59,19.04)(.051894,.068435){18}{\line(0,1){.068435}}
\multiput(178.52,20.27)(.049591,.070122){18}{\line(0,1){.070122}}
\multiput(179.41,21.53)(.050011,.075951){17}{\line(0,1){.075951}}
\multiput(180.26,22.82)(.050425,.082419){16}{\line(0,1){.082419}}
\multiput(181.07,24.14)(.050836,.089653){15}{\line(0,1){.089653}}
\multiput(181.83,25.49)(.051244,.097814){14}{\line(0,1){.097814}}
\multiput(182.55,26.85)(.051655,.107114){13}{\line(0,1){.107114}}
\multiput(183.22,28.25)(.048066,.108772){13}{\line(0,1){.108772}}
\multiput(183.85,29.66)(.04813,.1195){12}{\line(0,1){.1195}}
\multiput(184.42,31.1)(.04814,.13204){11}{\line(0,1){.13204}}
\multiput(184.95,32.55)(.0481,.14692){10}{\line(0,1){.14692}}
\multiput(185.43,34.02)(.04799,.16493){9}{\line(0,1){.16493}}
\multiput(185.87,35.5)(.04779,.18724){8}{\line(0,1){.18724}}
\multiput(186.25,37)(.04747,.21569){7}{\line(0,1){.21569}}
\multiput(186.58,38.51)(.04699,.25334){6}{\line(0,1){.25334}}
\multiput(186.86,40.03)(.04626,.30571){5}{\line(0,1){.30571}}
\multiput(187.09,41.56)(.0451,.3838){4}{\line(0,1){.3838}}
\multiput(187.27,43.09)(.0431,.5135){3}{\line(0,1){.5135}}
\put(187.4,44.63){\line(0,1){1.544}}
\put(187.48,46.18){\line(0,1){1.573}}
\multiput(111.75,61.25)(-.0638587,.05163043){184}{\line(-1,0){.0638587}}
\multiput(174,60.5)(.06428571,.05178571){140}{\line(1,0){.06428571}}
\put(142.75,18.25){\line(0,-1){17}}
\put(141.25,43){\circle{10.5}} \put(141.25,76.25){\circle{10.5}}
\put(166.75,29){\circle{10.5}} \put(114,29.25){\circle{10.5}}
\qbezier(110,32.75)(87.63,78.38)(136.75,79.5)
\qbezier(146.25,78)(199,72.88)(171.75,30.25)
\qbezier(164.75,24.5)(145,-4.13)(116.25,24.75)
\put(141.25,71){\line(0,-1){22.5}}
\multiput(118.75,31.5)(.11128049,.05182927){164}{\line(1,0){.11128049}}
\multiput(145.75,41)(.08938547,-.05167598){179}{\line(1,0){.08938547}}
\put(116,57.25){\makebox(0,0)[cc]{5}}
\put(169.5,56.75){\makebox(0,0)[cc]{5}}
\put(142.75,24.25){\makebox(0,0)[cc]{5}}
\put(165.25,90.5){\makebox(0,0)[cc]{5}}
\put(141.25,76.5){\makebox(0,0)[cc]{4}}
\put(141,42.75){\makebox(0,0)[cc]{4}}
\put(166.5,28.75){\makebox(0,0)[cc]{4}}
\put(113.5,29.25){\makebox(0,0)[cc]{4}}
\put(111.25,64.5){\makebox(0,0)[cc]{1}}
\put(125,59.75){\makebox(0,0)[cc]{3}}
\put(118.5,48.75){\makebox(0,0)[cc]{6}}
\put(174.5,63.75){\makebox(0,0)[cc]{6}}
\put(165.25,49.25){\makebox(0,0)[cc]{1}}
\put(144,16.75){\makebox(0,0)[cc]{3}}
\put(121.25,30.75){\makebox(0,0)[cc]{2}}
\put(110.5,36.75){\makebox(0,0)[cc]{0}}
\put(116,22){\makebox(0,0)[cc]{7}}
\put(165,21){\makebox(0,0)[cc]{7}}
\put(158.5,31){\makebox(0,0)[cc]{0}}
\put(175.25,31.5){\makebox(0,0)[cc]{2}}
\put(143,69){\makebox(0,0)[cc]{7}}
\put(132.75,82.25){\makebox(0,0)[cc]{0}}
\put(149,79.75){\makebox(0,0)[cc]{2}}
\put(134.5,28){\makebox(0,0)[cc]{6}}
\put(161.5,60){\makebox(0,0)[cc]{3}}
\end{picture}

\end{center}

\proof The set $X$ is $\{4,5\}$. The order on the graph $P_3(8)$:
$4\succ 5\succ 2\succ 3\succ 0\succ 1\succ 6\succ 7$, orders on the
anti-links:

\begin{itemize}
\item 2: $3\succ_2 1, 3\succ_27$,
\item 3: $2\succ_3 0, 2\succ_3 6$,
\item 4: $3\succ_41\succ_46$,
\item 5: $2\succ_50\succ_57$.
\end{itemize}

Let $p=(v_1,...,v_k)$ be an anti-path in $P_3(8)$. Consider eight
cases for $v_k$. Note that both the graph $P_3(8)$, the load
structure, and the dissection diagram admit an automorphism $0\to 1,
1\to 0, 2\to 3, 3\to 2, 4\to 5, 5\to 4, 6\to 7, 7\to 6$. Hence we
need to consider only 4 cases instead of 8.

Assume that $\Theta(p)$ does not satisfy (*) for $S[X]=S_0$ (the
seven-punctured disc).

{\bf 0.} Let $v_k=0$. Then $\Theta(p,X)\supseteq \{0,3,6\}$. It is
enough to show that either $1$ or $2$ is in $\Theta(p)$ because
$\{0,1,3,6\}$ and $\{0,2,3,6\}$ satisfy (*).

Since $2\succ 0$, we have $k\ne 1$, $v_{k-1}\in \{3,5,6\}$. Since
$2\succ_3 0, 2\succ_50$, $v_{k-1}=6$.

Since $1\succ 6$, we have $k-1\ne 1$, $v_{k-2}\in \{3,4\}$. But
$2\succ_3 6, 1\succ_4 6$, a contradiction.

{\bf 1.} Let $v_k=1$ (this case is symmetric to the previous one).
Then $\Theta(p)\supseteq \{1,2,7\}$. It is enough to show that
either $0$ or $3$ is in $\Theta(p)$ because $\{0,1,2,7\}$ and
$\{1,2,3,7\}$ satisfy (*).

Since $3\succ 1$, we have $k\ne 1$, $v_{k-1}\in \{2,4,7\}$. Since
$3\succ_2 1, 3\succ_41$, $v_{k-1}=7$.

Since $0\succ 7$, we have $k-1\ne 1$, $v_{k-2}\in \{2,5\}$. But
$3\succ_2 7, 0\succ_5 7$, a contradiction.

{\bf 2.} Let $v_k=2$. Then $\Theta(p)\supseteq \{1,2,3,7\}$ and
satisfies (*).

{\bf 3.} Let $v_k=3$. Then $\Theta(p)\supseteq \{0,2,3,6\}$ and
satisfies (*).

{\bf 4.} Let $v_k=4$ (and $k=1$). Then $\Theta(p)\supseteq
\{1,3,6\}$ and satisfies (*).

{\bf 5.} Let $v_k=5$ (and $k=1$). Then $\Theta(p)\supseteq
\{0,2,7\}$ and satisfies (*).

{\bf 6.} Let $v_k=6$. Then $\Theta(p)\supseteq \{0,3,6\}$. It is
enough to show that $1\in \Theta(p)$ or $2\in \Theta(p)$ because
$\{0,1,3,6\}$ and $\{0,2,3,6\}$ satisfy (*).

Since $1\succ 6$, $k\ne 1$, $v_{k-1}\in \{0,3,4\}$. Since $1\succ_4
6$, $2\succ_3 6$, we have $v_{k-1}=0$.

Since $2\succ 0$, $k-1\ne 1$, $v_{k-2}\in \{3,5\}$. But $2\succ_3 0,
2\succ_50$, a contradiction.

{\bf 7.}  Let $v_k=7$ (this case is symmetric to the previous one).
Then $\Theta(p)\supseteq \{1,2,7\}$. It is enough to show that $0\in
\Theta(p)$ or $3\in \Theta(p)$ because $\{0,1,2,7\}$, $\{1,2,3,7\}$
satisfy (*).

Since $0\succ 7$, $k\ne 1$, $v_{k-1}\in \{1,2,5\}$. Since $0\succ_5
7$, $3\succ_2 7$, we have $v_{k-1}=1$.

Since $3\succ 1$, $k-1\ne 1$, $v_{k-2}\in \{2,4\}$. But $3\succ_2 1,
3\succ_41$, a contradiction.
\endproof

\subsubsection{The 8-vertex graph $P_4(8)$}

\begin{center}

\unitlength .75mm 
\linethickness{0.4pt}
\ifx\plotpoint\undefined\newsavebox{\plotpoint}\fi 
\begin{picture}(199,90.75)(0,0)
\put(7.5,23.75){\line(1,0){63.25}}
\multiput(70.75,23.75)(-.051835853,.058855292){463}{\line(0,1){.058855292}}
\put(46.75,51){\line(0,-1){27}}
\multiput(46.75,24)(-.05304878,-.051829268){410}{\line(-1,0){.05304878}}
\put(25,2.75){\line(0,1){21}}
\multiput(25,23.75)(.051807229,.065060241){415}{\line(0,1){.065060241}}
\multiput(7.5,23.75)(.051801802,-.061561562){333}{\line(0,-1){.061561562}}
\put(46.75,35.5){\line(2,-1){23.5}}
\qbezier(70.5,24)(42.5,106.75)(7.5,23.5)
\put(4.5,23.25){\makebox(0,0)[cc]{0}}
\put(73.25,23.5){\makebox(0,0)[cc]{1}}
\put(28.75,3.25){\makebox(0,0)[cc]{2}}
\put(22.25,14){\makebox(0,0)[cc]{4}}
\put(22.75,27.25){\makebox(0,0)[cc]{7}}
\put(46.75,21.25){\makebox(0,0)[cc]{6}}
\put(43,35.5){\makebox(0,0)[cc]{5}}
\put(47,54.25){\makebox(0,0)[cc]{3}}
\put(10.25,52.5){\makebox(0,0)[cc]{$P_4(8)$}}
\put(116,57){\circle{11.85}} \put(142.75,24){\circle{11.85}}
\put(169.25,57){\circle{11.85}}
\put(121.75,57.25){\line(1,0){41.5}}
\multiput(165.25,52.75)(-.051771117,-.066076294){367}{\line(0,-1){.066076294}}
\multiput(118.5,51.75)(.051886792,-.064690027){371}{\line(0,-1){.064690027}}
\put(187.51,47.75){\line(0,1){1.546}}
\put(187.48,49.3){\line(0,1){1.544}}
\multiput(187.41,50.84)(-.0428,.5135){3}{\line(0,1){.5135}}
\multiput(187.28,52.38)(-.0449,.3839){4}{\line(0,1){.3839}}
\multiput(187.1,53.92)(-.04608,.30574){5}{\line(0,1){.30574}}
\multiput(186.87,55.44)(-.04684,.25336){6}{\line(0,1){.25336}}
\multiput(186.59,56.96)(-.04735,.21571){7}{\line(0,1){.21571}}
\multiput(186.26,58.47)(-.04768,.18727){8}{\line(0,1){.18727}}
\multiput(185.87,59.97)(-.04789,.16496){9}{\line(0,1){.16496}}
\multiput(185.44,61.46)(-.04801,.14695){10}{\line(0,1){.14695}}
\multiput(184.96,62.93)(-.04806,.13207){11}{\line(0,1){.13207}}
\multiput(184.43,64.38)(-.04806,.11953){12}{\line(0,1){.11953}}
\multiput(183.86,65.81)(-.048003,.108799){13}{\line(0,1){.108799}}
\multiput(183.23,67.23)(-.051592,.107144){13}{\line(0,1){.107144}}
\multiput(182.56,68.62)(-.051187,.097844){14}{\line(0,1){.097844}}
\multiput(181.85,69.99)(-.050783,.089682){15}{\line(0,1){.089682}}
\multiput(181.08,71.34)(-.050378,.082449){16}{\line(0,1){.082449}}
\multiput(180.28,72.66)(-.049967,.07598){17}{\line(0,1){.07598}}
\multiput(179.43,73.95)(-.04955,.070151){18}{\line(0,1){.070151}}
\multiput(178.54,75.21)(-.051854,.068465){18}{\line(0,1){.068465}}
\multiput(177.6,76.44)(-.051253,.063193){19}{\line(0,1){.063193}}
\multiput(176.63,77.64)(-.050659,.058382){20}{\line(0,1){.058382}}
\multiput(175.62,78.81)(-.050068,.053968){21}{\line(0,1){.053968}}
\multiput(174.57,79.94)(-.051834,.052274){21}{\line(0,1){.052274}}
\multiput(173.48,81.04)(-.053543,.050523){21}{\line(-1,0){.053543}}
\multiput(172.35,82.1)(-.057952,.051151){20}{\line(-1,0){.057952}}
\multiput(171.19,83.13)(-.062758,.051786){19}{\line(-1,0){.062758}}
\multiput(170,84.11)(-.064444,.049671){19}{\line(-1,0){.064444}}
\multiput(168.78,85.05)(-.069729,.050141){18}{\line(-1,0){.069729}}
\multiput(167.52,85.96)(-.075555,.050608){17}{\line(-1,0){.075555}}
\multiput(166.24,86.82)(-.08202,.051073){16}{\line(-1,0){.08202}}
\multiput(164.92,87.63)(-.08925,.05154){15}{\line(-1,0){.08925}}
\multiput(163.59,88.41)(-.090913,.048545){15}{\line(-1,0){.090913}}
\multiput(162.22,89.13)(-.099082,.048747){14}{\line(-1,0){.099082}}
\multiput(160.83,89.82)(-.10839,.048921){13}{\line(-1,0){.10839}}
\multiput(159.43,90.45)(-.11912,.04907){12}{\line(-1,0){.11912}}
\multiput(158,91.04)(-.13165,.04918){11}{\line(-1,0){.13165}}
\multiput(156.55,91.58)(-.14654,.04925){10}{\line(-1,0){.14654}}
\multiput(155.08,92.08)(-.16455,.04928){9}{\line(-1,0){.16455}}
\multiput(153.6,92.52)(-.18686,.04926){8}{\line(-1,0){.18686}}
\multiput(152.11,92.91)(-.21531,.04917){7}{\line(-1,0){.21531}}
\multiput(150.6,93.26)(-.25296,.04899){6}{\line(-1,0){.25296}}
\multiput(149.08,93.55)(-.30533,.04866){5}{\line(-1,0){.30533}}
\multiput(147.56,93.79)(-.3835,.0481){4}{\line(-1,0){.3835}}
\multiput(146.02,93.99)(-.5132,.0471){3}{\line(-1,0){.5132}}
\put(144.48,94.13){\line(-1,0){1.543}}
\put(142.94,94.22){\line(-1,0){1.545}}
\put(141.39,94.26){\line(-1,0){1.546}}
\put(139.85,94.24){\line(-1,0){1.545}}
\put(138.3,94.18){\line(-1,0){1.542}}
\multiput(136.76,94.07)(-.3842,-.0416){4}{\line(-1,0){.3842}}
\multiput(135.22,93.9)(-.30611,-.04349){5}{\line(-1,0){.30611}}
\multiput(133.69,93.68)(-.25375,-.0447){6}{\line(-1,0){.25375}}
\multiput(132.17,93.41)(-.21611,-.04552){7}{\line(-1,0){.21611}}
\multiput(130.66,93.09)(-.18766,-.04609){8}{\line(-1,0){.18766}}
\multiput(129.16,92.73)(-.16536,-.04649){9}{\line(-1,0){.16536}}
\multiput(127.67,92.31)(-.14735,-.04677){10}{\line(-1,0){.14735}}
\multiput(126.2,91.84)(-.14571,-.05164){10}{\line(-1,0){.14571}}
\multiput(124.74,91.32)(-.13083,-.05132){11}{\line(-1,0){.13083}}
\multiput(123.3,90.76)(-.1183,-.051){12}{\line(-1,0){.1183}}
\multiput(121.88,90.15)(-.107576,-.050685){13}{\line(-1,0){.107576}}
\multiput(120.48,89.49)(-.098273,-.050358){14}{\line(-1,0){.098273}}
\multiput(119.11,88.78)(-.090109,-.050023){15}{\line(-1,0){.090109}}
\multiput(117.75,88.03)(-.082872,-.049679){16}{\line(-1,0){.082872}}
\multiput(116.43,87.24)(-.0764,-.049323){17}{\line(-1,0){.0764}}
\multiput(115.13,86.4)(-.074718,-.051835){17}{\line(-1,0){.074718}}
\multiput(113.86,85.52)(-.068901,-.051273){18}{\line(-1,0){.068901}}
\multiput(112.62,84.59)(-.063624,-.050717){19}{\line(-1,0){.063624}}
\multiput(111.41,83.63)(-.058808,-.050164){20}{\line(-1,0){.058808}}
\multiput(110.23,82.63)(-.054389,-.04961){21}{\line(-1,0){.054389}}
\multiput(109.09,81.59)(-.052711,-.05139){21}{\line(-1,0){.052711}}
\multiput(107.98,80.51)(-.050973,-.053114){21}{\line(0,-1){.053114}}
\multiput(106.91,79.39)(-.051639,-.057517){20}{\line(0,-1){.057517}}
\multiput(105.88,78.24)(-.049699,-.059202){20}{\line(0,-1){.059202}}
\multiput(104.89,77.06)(-.050215,-.064022){19}{\line(0,-1){.064022}}
\multiput(103.93,75.84)(-.050729,-.069303){18}{\line(0,-1){.069303}}
\multiput(103.02,74.59)(-.051245,-.075124){17}{\line(0,-1){.075124}}
\multiput(102.15,73.32)(-.051765,-.081585){16}{\line(0,-1){.081585}}
\multiput(101.32,72.01)(-.049024,-.08326){16}{\line(0,-1){.08326}}
\multiput(100.54,70.68)(-.049312,-.0905){15}{\line(0,-1){.0905}}
\multiput(99.8,69.32)(-.049583,-.098666){14}{\line(0,-1){.098666}}
\multiput(99.1,67.94)(-.049836,-.107972){13}{\line(0,-1){.107972}}
\multiput(98.45,66.54)(-.05007,-.1187){12}{\line(0,-1){.1187}}
\multiput(97.85,65.11)(-.05029,-.13123){11}{\line(0,-1){.13123}}
\multiput(97.3,63.67)(-.05049,-.14612){10}{\line(0,-1){.14612}}
\multiput(96.8,62.21)(-.05067,-.16413){9}{\line(0,-1){.16413}}
\multiput(96.34,60.73)(-.05084,-.18643){8}{\line(0,-1){.18643}}
\multiput(95.93,59.24)(-.05099,-.21488){7}{\line(0,-1){.21488}}
\multiput(95.58,57.73)(-.05112,-.25253){6}{\line(0,-1){.25253}}
\multiput(95.27,56.22)(-.05124,-.30491){5}{\line(0,-1){.30491}}
\multiput(95.01,54.69)(-.0513,-.3831){4}{\line(0,-1){.3831}}
\multiput(94.81,53.16)(-.0515,-.5127){3}{\line(0,-1){.5127}}
\put(94.65,51.62){\line(0,-1){1.542}}
\put(94.55,50.08){\line(0,-1){6.179}}
\multiput(94.65,43.9)(.0512,-.5128){3}{\line(0,-1){.5128}}
\multiput(94.8,42.36)(.0511,-.3831){4}{\line(0,-1){.3831}}
\multiput(95.01,40.83)(.05106,-.30494){5}{\line(0,-1){.30494}}
\multiput(95.26,39.31)(.05098,-.25256){6}{\line(0,-1){.25256}}
\multiput(95.57,37.79)(.05086,-.21491){7}{\line(0,-1){.21491}}
\multiput(95.93,36.29)(.05073,-.18646){8}{\line(0,-1){.18646}}
\multiput(96.33,34.8)(.05058,-.16416){9}{\line(0,-1){.16416}}
\multiput(96.79,33.32)(.05041,-.14615){10}{\line(0,-1){.14615}}
\multiput(97.29,31.86)(.05021,-.13126){11}{\line(0,-1){.13126}}
\multiput(97.84,30.41)(.05,-.11873){12}{\line(0,-1){.11873}}
\multiput(98.44,28.99)(.049773,-.108001){13}{\line(0,-1){.108001}}
\multiput(99.09,27.58)(.049525,-.098695){14}{\line(0,-1){.098695}}
\multiput(99.78,26.2)(.04926,-.090528){15}{\line(0,-1){.090528}}
\multiput(100.52,24.84)(.048976,-.083289){16}{\line(0,-1){.083289}}
\multiput(101.31,23.51)(.051717,-.081615){16}{\line(0,-1){.081615}}
\multiput(102.13,22.21)(.051201,-.075154){17}{\line(0,-1){.075154}}
\multiput(103,20.93)(.050689,-.069332){18}{\line(0,-1){.069332}}
\multiput(103.92,19.68)(.050177,-.064051){19}{\line(0,-1){.064051}}
\multiput(104.87,18.46)(.049665,-.059231){20}{\line(0,-1){.059231}}
\multiput(105.86,17.28)(.051606,-.057547){20}{\line(0,-1){.057547}}
\multiput(106.9,16.13)(.050943,-.053143){21}{\line(0,-1){.053143}}
\multiput(107.97,15.01)(.052681,-.051421){21}{\line(1,0){.052681}}
\multiput(109.07,13.93)(.054361,-.049641){21}{\line(1,0){.054361}}
\multiput(110.21,12.89)(.058779,-.050198){20}{\line(1,0){.058779}}
\multiput(111.39,11.89)(.063595,-.050754){19}{\line(1,0){.063595}}
\multiput(112.6,10.92)(.068871,-.051313){18}{\line(1,0){.068871}}
\multiput(113.84,10)(.074688,-.051878){17}{\line(1,0){.074688}}
\multiput(115.11,9.12)(.076371,-.049367){17}{\line(1,0){.076371}}
\multiput(116.4,8.28)(.082843,-.049727){16}{\line(1,0){.082843}}
\multiput(117.73,7.48)(.09008,-.050076){15}{\line(1,0){.09008}}
\multiput(119.08,6.73)(.098244,-.050415){14}{\line(1,0){.098244}}
\multiput(120.46,6.02)(.107547,-.050747){13}{\line(1,0){.107547}}
\multiput(121.85,5.36)(.11827,-.05107){12}{\line(1,0){.11827}}
\multiput(123.27,4.75)(.1308,-.0514){11}{\line(1,0){.1308}}
\multiput(124.71,4.19)(.14568,-.05172){10}{\line(1,0){.14568}}
\multiput(126.17,3.67)(.14732,-.04685){10}{\line(1,0){.14732}}
\multiput(127.64,3.2)(.16533,-.04659){9}{\line(1,0){.16533}}
\multiput(129.13,2.78)(.18764,-.0462){8}{\line(1,0){.18764}}
\multiput(130.63,2.41)(.21608,-.04565){7}{\line(1,0){.21608}}
\multiput(132.14,2.09)(.25372,-.04485){6}{\line(1,0){.25372}}
\multiput(133.67,1.82)(.30609,-.04367){5}{\line(1,0){.30609}}
\multiput(135.2,1.6)(.3842,-.0418){4}{\line(1,0){.3842}}
\put(136.73,1.44){\line(1,0){1.542}}
\put(138.28,1.32){\line(1,0){1.545}}
\put(139.82,1.26){\line(1,0){1.546}}
\put(141.37,1.24){\line(1,0){1.545}}
\put(142.91,1.28){\line(1,0){1.543}}
\multiput(144.46,1.37)(.5132,.0468){3}{\line(1,0){.5132}}
\multiput(145.99,1.51)(.3835,.0479){4}{\line(1,0){.3835}}
\multiput(147.53,1.7)(.30536,.04848){5}{\line(1,0){.30536}}
\multiput(149.06,1.94)(.25299,.04884){6}{\line(1,0){.25299}}
\multiput(150.57,2.24)(.21533,.04905){7}{\line(1,0){.21533}}
\multiput(152.08,2.58)(.18689,.04915){8}{\line(1,0){.18689}}
\multiput(153.58,2.97)(.16458,.04919){9}{\line(1,0){.16458}}
\multiput(155.06,3.42)(.14657,.04917){10}{\line(1,0){.14657}}
\multiput(156.52,3.91)(.13168,.0491){11}{\line(1,0){.13168}}
\multiput(157.97,4.45)(.11915,.049){12}{\line(1,0){.11915}}
\multiput(159.4,5.04)(.108418,.048858){13}{\line(1,0){.108418}}
\multiput(160.81,5.67)(.09911,.048689){14}{\line(1,0){.09911}}
\multiput(162.2,6.35)(.090942,.048492){15}{\line(1,0){.090942}}
\multiput(163.56,7.08)(.08928,.051488){15}{\line(1,0){.08928}}
\multiput(164.9,7.85)(.082049,.051025){16}{\line(1,0){.082049}}
\multiput(166.21,8.67)(.075584,.050564){17}{\line(1,0){.075584}}
\multiput(167.5,9.53)(.069758,.050101){18}{\line(1,0){.069758}}
\multiput(168.75,10.43)(.064473,.049634){19}{\line(1,0){.064473}}
\multiput(169.98,11.37)(.062788,.051749){19}{\line(1,0){.062788}}
\multiput(171.17,12.36)(.057982,.051117){20}{\line(1,0){.057982}}
\multiput(172.33,13.38)(.053572,.050491){21}{\line(1,0){.053572}}
\multiput(173.46,14.44)(.051864,.052244){21}{\line(0,1){.052244}}
\multiput(174.55,15.54)(.050099,.053939){21}{\line(0,1){.053939}}
\multiput(175.6,16.67)(.050693,.058353){20}{\line(0,1){.058353}}
\multiput(176.61,17.84)(.05129,.063164){19}{\line(0,1){.063164}}
\multiput(177.59,19.04)(.051894,.068435){18}{\line(0,1){.068435}}
\multiput(178.52,20.27)(.049591,.070122){18}{\line(0,1){.070122}}
\multiput(179.41,21.53)(.050011,.075951){17}{\line(0,1){.075951}}
\multiput(180.26,22.82)(.050425,.082419){16}{\line(0,1){.082419}}
\multiput(181.07,24.14)(.050836,.089653){15}{\line(0,1){.089653}}
\multiput(181.83,25.49)(.051244,.097814){14}{\line(0,1){.097814}}
\multiput(182.55,26.85)(.051655,.107114){13}{\line(0,1){.107114}}
\multiput(183.22,28.25)(.048066,.108772){13}{\line(0,1){.108772}}
\multiput(183.85,29.66)(.04813,.1195){12}{\line(0,1){.1195}}
\multiput(184.42,31.1)(.04814,.13204){11}{\line(0,1){.13204}}
\multiput(184.95,32.55)(.0481,.14692){10}{\line(0,1){.14692}}
\multiput(185.43,34.02)(.04799,.16493){9}{\line(0,1){.16493}}
\multiput(185.87,35.5)(.04779,.18724){8}{\line(0,1){.18724}}
\multiput(186.25,37)(.04747,.21569){7}{\line(0,1){.21569}}
\multiput(186.58,38.51)(.04699,.25334){6}{\line(0,1){.25334}}
\multiput(186.86,40.03)(.04626,.30571){5}{\line(0,1){.30571}}
\multiput(187.09,41.56)(.0451,.3838){4}{\line(0,1){.3838}}
\multiput(187.27,43.09)(.0431,.5135){3}{\line(0,1){.5135}}
\put(187.4,44.63){\line(0,1){1.544}}
\put(187.48,46.18){\line(0,1){1.573}}
\multiput(111.75,61.25)(-.0638587,.05163043){184}{\line(-1,0){.0638587}}
\multiput(174,60.5)(.06428571,.05178571){140}{\line(1,0){.06428571}}
\put(142.75,18.25){\line(0,-1){17}}
\put(141.25,43){\circle{10.5}} \put(141.25,76.25){\circle{10.5}}
\put(166.75,29){\circle{10.5}} \put(114,29.25){\circle{10.5}}
\qbezier(110,32.75)(87.63,78.38)(136.75,79.5)
\qbezier(146.25,78)(199,72.88)(171.75,30.25)
\qbezier(164.75,24.5)(145,-4.13)(116.25,24.75)
\put(141.25,71){\line(0,-1){22.5}}
\multiput(118.75,31.5)(.11128049,.05182927){164}{\line(1,0){.11128049}}
\multiput(145.75,41)(.08938547,-.05167598){179}{\line(1,0){.08938547}}
\put(116,57.25){\makebox(0,0)[cc]{5}}
\put(169.5,56.75){\makebox(0,0)[cc]{5}}
\put(142.75,24.25){\makebox(0,0)[cc]{5}}
\put(165.25,90.5){\makebox(0,0)[cc]{5}}
\put(141.25,76.5){\makebox(0,0)[cc]{4}}
\put(141,42.75){\makebox(0,0)[cc]{4}}
\put(166.5,28.75){\makebox(0,0)[cc]{4}}
\put(113.5,29.25){\makebox(0,0)[cc]{4}}
\put(111.25,64.5){\makebox(0,0)[cc]{3}}
\put(125,59.75){\makebox(0,0)[cc]{1}}
\put(118.5,48.75){\makebox(0,0)[cc]{6}}
\put(174.5,63.75){\makebox(0,0)[cc]{6}}
\put(165.25,49.25){\makebox(0,0)[cc]{3}}
\put(144,16.75){\makebox(0,0)[cc]{1}}
\put(110.5,36.75){\makebox(0,0)[cc]{7}}
\put(116,22){\makebox(0,0)[cc]{6}}
\put(165,21){\makebox(0,0)[cc]{6}}
\put(158.5,31){\makebox(0,0)[cc]{7}}
\put(175.25,31.5){\makebox(0,0)[cc]{2}}
\put(143,69){\makebox(0,0)[cc]{6}}
\put(132.75,82.25){\makebox(0,0)[cc]{7}}
\put(149,79.75){\makebox(0,0)[cc]{2}}
\put(134.5,28){\makebox(0,0)[cc]{6}}
\put(161.5,60){\makebox(0,0)[cc]{1}}
\multiput(25,13.5)(.10984848,.05176768){198}{\line(1,0){.10984848}}
\qbezier(141.5,34.75)(132.38,35.25)(132.75,59.75)
\qbezier(132.75,59.75)(132.75,85.13)(140.75,86)
\qbezier(140.75,86)(148.25,87)(148.75,60)
\qbezier(148.75,60)(148.25,34.13)(141.75,34.75)
\put(150.5,63.75){\makebox(0,0)[cc]{0}}
\put(180.75,37.5){\makebox(0,0)[cc]{0}}
\qbezier(106,42)(92.13,34.13)(117.75,12.75)
\qbezier(117.75,12.75)(142.75,-7.13)(166.75,15.5)
\qbezier(166.75,15.5)(184.63,34.63)(177,39.25)
\qbezier(176.75,39.5)(165.5,47)(154.25,20.5)
\qbezier(154.25,20.5)(141.75,3.88)(126.25,23.75)
\qbezier(126.25,23.75)(110.75,44.38)(106.25,41.5)
\put(127,33.25){\makebox(0,0)[cc]{2}}
\put(110.25,16.5){\makebox(0,0)[cc]{0}}
\end{picture}

\end{center}

\proof Let $X=\{4,5\}$. The order on the vertices of $P_4(8)$ is
$5\succ 4\succ 3\succ 2\succ 0\succ 1\succ 6\succ 7$. The partial
orders on the links are:

\begin{itemize}
\item 2: $3\succ_2 1$, $3\succ_2 7$,
\item 3: $2\succ_3 6$, $2\succ_3 0$,
\item 4: $3\succ_4 1$,
\item 5: $2\succ_5 0\succ_5 7$.
\end{itemize}

Let $p=(v_1,...,v_k)$ be an anti-path in $K_2$  where only $v_1$ can
be equal to $5$ or $4$. We consider eight cases for $v_k$. The set
$\Theta(p)$ must contain $4,5$, so we should consider the surface
$S[4,5]$ (the seven-punctured disc) and its subsurfaces. Assume that
$\Theta(p)$ does not satisfy (*).

{\bf 0.} Let $v_k=0$. Then $\Theta(p)$ contains $\{0,3,6\}$ and
satisfies (*).

{\bf 1.} Let $v_k=1$. Them $\Theta(p)$ contains $\{1,2,7\}$. Since
$0\succ 1$ and $\{0,1,2,7\}$ satisfies (*), $k\ne 1$, $v_{k-1}\in
\{2,4,7\}$. Since $3\succ_2 1$, $3\succ_4 1$, and $\{1,2,3,4,7\}$
satisfies (*), $v_{k-1}=7$. Since $0\succ 7$, $k-1\ne 1$,
$v_{k-2}\in \{2,5\}$. Since $3\succ_2 7, 0\succ_5 7$, we get a
contradiction.

{\bf 2.} Let $v_k=2$. Then $\Theta(p)$ contains $\{1,2,3,7\}$ and
satisfies (*).

{\bf 3.} Let $v_k=3$. Then $\Theta(p)$ contains $\{0,2,3,6\}$ and
satisfies (*).

{\bf 4.} Let $v_k=4$. Then $\Theta(p)$ contains $\{0,1,3\}$ and
satisfies (*).

{\bf 5.} Let $v_k=5$. Then $\Theta(p)$ contains $\{0,2,7\}$ and
satisfies (*).

{\bf 6.} Let $v_k=6$. Then $\Theta(p)$ contains $\{0,3,6\}$. Since
$2\succ 6$ and $\{0,2,3,6\}$ satisfies (*), $k\ne 1$, $v_{k-1}\in
\{0,3\}$. Since $2\succ_3 6$, $v_{k-1}=0$. Since $2\succ 0$, $k-1\ne
1$, $v_{k-2}\in \{3,5\}$. But $2\succ_3 0, 2\succ_50$, a
contradiction.

{\bf 7.} Let $v_k=7$. Then $\Theta(p)$ contains $\{1,2,7\}$. Since
$0\succ 7$ and $\{0,1,2,7\}$ satisfies (*), $k\ne 1$, $v_{k-1}\in
\{1,2,5\}$. Since $0\succ_5 7$, $3\succ_2 7$, and $\{1,2,3,7\}$
satisfies (*), $v_{k-1}=1$. Since $3\succ 1$, $k-1\ne 1$,
$v_{k-2}\in \{2,4\}$. But $3\succ_2 1, 3\succ_4 1$, a contradiction.
\endproof

\subsection{Non-oriented surface subgroups of genus 2}
\label{nonor}
 Here we present a different way of creating dissection
diagrams on surfaces. As examples, we present dissection diagrams
for non-orientable surfaces of (non-orientable) genus 2 for graphs
$P_1(6), P_2(6)$, $P_2(8)$, $P_3(8)$, $P_4(8)$.

\subsubsection{Graphs $P_1(6), P_2(6)$}

\begin{lemma}\label{prisms}
Suppose that $K$ contains non-adjacent vertices $u,v$ and vertices
$c\in\Lk(u)\setminus\Lk(v)$, $d\in\Lk(v)\setminus\Lk(u)$, and
$e,f\in\Lk(u)\cap\Lk(v)$ (not necessarily distinct) such that $c,d$
are adjacent, while $c,e$ are non-adjacent and $d, f$ are
non-adjacent. Then $A(K)$ contains a subgroup isomorphic to the
fundamental group of a nonorientable closed surface of Euler
characteristic $-2$ (hence non-orientable genus $4$).
\end{lemma}

\Remark Note that if $K$ contains the above configuration, then the
full subcomplex spanned by $u,v,c,d,e,f$ is either the circuit of
length 5, or graph isomorphic to $P_1(6)$ or $P_2(6)$.

\unitlength .5mm 
\linethickness{0.4pt}
\ifx\plotpoint\undefined\newsavebox{\plotpoint}\fi 
\begin{picture}(128.75,108.5)(0,0)
\multiput(57.25,97.75)(10.5625,-.0625){4}{\line(1,0){10.5625}}
\multiput(58,98)(-.06741573,-.074438202){356}{\line(0,-1){.074438202}}
\put(34,71.5){\line(0,-1){38.75}}
\multiput(34,32.75)(.067471591,-.069602273){352}{\line(0,-1){.069602273}}
\put(57.75,8.25){\line(1,0){41.75}}
\multiput(99.5,8.25)(.071527778,.067361111){360}{\line(1,0){.071527778}}
\put(125.25,32.5){\line(0,1){39}}
\multiput(125.25,71.5)(-.067408377,.068062827){382}{\line(0,1){.068062827}}
\put(46.25,85.25){\vector(3,4){.14}}\multiput(41.25,79.25)(.0666667,.08){75}{\line(0,1){.08}}
\put(106.25,14.75){\vector(-1,-1){.14}}\multiput(115.5,23.25)(-.0734127,-.06746032){126}{\line(-1,0){.0734127}}
\put(34,53.5){\vector(0,1){.14}}\put(34,47.75){\line(0,1){5.75}}
\put(125.25,45.25){\vector(0,-1){.14}}\put(125.25,57){\line(0,-1){11.75}}
\put(76.75,8.25){\vector(-1,0){.14}}\put(88.75,8.25){\line(-1,0){12}}
\put(82.5,97.75){\vector(1,0){.14}}\put(75,97.75){\line(1,0){7.5}}
\put(113,84){\vector(1,-1){.14}}\multiput(107.75,89.5)(.0673077,-.0705128){78}{\line(0,-1){.0705128}}
\put(43.5,23){\vector(-1,1){.14}}\multiput(48.5,18)(-.0666667,.0666667){75}{\line(0,1){.0666667}}
\multiput(34,61.75)(.067346939,.073469388){490}{\line(0,1){.073469388}}
\multiput(90.5,97.25)(.067382813,-.067871094){512}{\line(0,-1){.067871094}}
\multiput(125.25,40.25)(-.074309979,-.067409766){471}{\line(-1,0){.074309979}}
\multiput(65.25,8)(-.06745182,.071734475){467}{\line(0,1){.071734475}}
\put(77.5,97.75){\line(0,-1){89.25}}
\put(34,50){\line(1,0){91.25}}
\put(56.25,80.25){\vector(1,-1){.14}}\multiput(51.75,85.25)(.0671642,-.0746269){67}{\line(0,-1){.0746269}}
\put(98.25,81.75){\vector(-4,-3){.14}}\multiput(105,87.5)(-.0784884,-.0668605){86}{\line(-1,0){.0784884}}
\put(56.5,24){\vector(4,3){.14}}\multiput(50.75,19.25)(.0809859,.0669014){71}{\line(1,0){.0809859}}
\put(109.25,20.75){\vector(3,-4){.14}}\multiput(104.5,27)(.0669014,-.0880282){71}{\line(0,-1){.0880282}}
\put(61,47){\vector(0,-1){.14}}\put(61,54){\line(0,-1){7}}
\put(82.5,71.25){\vector(1,0){.14}}\put(75.25,71.25){\line(1,0){7.25}}
\put(30.5,52.75){\makebox(0,0)[cc]{$v$}}
\put(128.75,51.75){\makebox(0,0)[cc]{$v$}}
\put(78.25,101.5){\makebox(0,0)[cc]{$u$}}
\put(77.25,5.75){\makebox(0,0)[cc]{$u$}}
\put(60.25,85.5){\makebox(0,0)[cc]{$e$}}
\put(107.25,76){\makebox(0,0)[cc]{$e$}}
\put(109.25,31.25){\makebox(0,0)[cc]{$f$}}
\put(51.25,28){\makebox(0,0)[cc]{$f$}}
\put(80.25,81.25){\makebox(0,0)[cc]{$c$}}
\put(55.25,53.25){\makebox(0,0)[cc]{$d$}}
\put(32.25,72.5){\makebox(0,0)[cc]{$p_3$}}
\put(56.25,100.5){\makebox(0,0)[cc]{$p_4$}}
\put(102,99.75){\makebox(0,0)[cc]{$p_5$}}
\put(128,72){\makebox(0,0)[cc]{$p_6$}}
\put(128.75,31.25){\makebox(0,0)[cc]{$p_7$}}
\put(102.25,6.5){\makebox(0,0)[cc]{$p_8$}}
\put(54.5,6){\makebox(0,0)[cc]{$p_1$}}
\put(29.5,31.75){\makebox(0,0)[cc]{$p_2$}}
\end{picture}

\begin{proof} We construct a dissected surface . Let $p_1,..,p_8$
denote the vertices of a regular octagon in cyclic order. Let $S_0$
denote the surface obtained from this by gluing oriented sides
$[p_1,p_2]$ to $[p_5,p_6]$ and $[p_3,p_4]$ to $[p_7,p_8]$. This is a
nonorientable surface of Euler characteristic $-1$ with two boundary
components, which we shall view as curves of the dissection. Label
the boundary $[p_8,p_1]\cup[p_4,p_5]$ with a $u$ and the other with
a $v$, both oriented out of the surface. Other dissection curves and
their orientation are shown on the picture.

Now, construct $S=S_0\cup S_1$ by doubling $S_0$ along the $u$- and
$v$-curves. Note that $\pi_1(S)$ has graph of groups decomposition
with two vertex groups each isomorphic to $\pi_1(S_0)=F_2$, two
edges with infinite cyclic edge groups generated by the $u$- and
$v$-curves respectively, and stable letters represented in $A(K)$ by
$u$ and $v$ respectively.

To prove that the dissection diagram is faithful, we use Corollary
\ref{corcuts} again. The stable set $X$ is $\{u, v\}$.

\end{proof}

\subsubsection{Graphs $P_2(8), P_3(8), P_4(8)$}

\begin{center}
\unitlength .6mm 
\linethickness{0.4pt}
\ifx\plotpoint\undefined\newsavebox{\plotpoint}\fi 
\begin{picture}(280.5,230.28)(0,0)
\put(7.5,23.75){\line(1,0){63.25}}
\put(7.5,102.25){\line(1,0){63.25}}
\put(7.5,182.28){\line(1,0){63.25}}
\multiput(70.75,23.75)(-.056206089,.06381733){427}{\line(0,1){.06381733}}
\multiput(70.75,102.25)(-.056206089,.06381733){427}{\line(0,1){.06381733}}
\multiput(70.75,182.28)(-.056206089,.06381733){427}{\line(0,1){.06381733}}
\put(46.75,51){\line(0,-1){27}}
\put(46.75,129.5){\line(0,-1){27}}
\put(46.75,209.53){\line(0,-1){27}}
\multiput(46.75,24)(-.057539683,-.056216931){378}{\line(-1,0){.057539683}}
\multiput(46.75,102.5)(-.057539683,-.056216931){378}{\line(-1,0){.057539683}}
\multiput(46.75,182.53)(-.057539683,-.056216931){378}{\line(-1,0){.057539683}}
\put(25,2.75){\line(0,1){21}}
\put(25,81.25){\line(0,1){21}}
\put(25,161.28){\line(0,1){21}}
\multiput(25,23.75)(.05613577,.070496084){383}{\line(0,1){.070496084}}
\multiput(25,102.25)(.05613577,.070496084){383}{\line(0,1){.070496084}}
\multiput(25,182.28)(.05613577,.070496084){383}{\line(0,1){.070496084}}
\multiput(7.5,23.75)(.056188925,-.066775244){307}{\line(0,-1){.066775244}}
\multiput(7.5,102.25)(.056188925,-.066775244){307}{\line(0,-1){.066775244}}
\multiput(7.5,182.28)(.056188925,-.066775244){307}{\line(0,-1){.066775244}}
\multiput(46.75,35.5)(.11244019,-.0562201){209}{\line(1,0){.11244019}}
\multiput(46.75,114)(.11244019,-.0562201){209}{\line(1,0){.11244019}}
\qbezier(70.5,24)(42.5,106.75)(7.5,23.5)
\qbezier(70.5,102.5)(42.5,185.25)(7.5,102)
\qbezier(70.5,182.53)(42.5,265.28)(7.5,182.03)
\put(4.5,23.25){\makebox(0,0)[cc]{0}}
\put(4.5,101.75){\makebox(0,0)[cc]{0}}
\put(4.5,181.78){\makebox(0,0)[cc]{0}}
\put(73.25,23.5){\makebox(0,0)[cc]{1}}
\put(73.25,102){\makebox(0,0)[cc]{1}}
\put(73.25,182.03){\makebox(0,0)[cc]{1}}
\put(28.75,3.25){\makebox(0,0)[cc]{2}}
\put(28.75,81.75){\makebox(0,0)[cc]{2}}
\put(28.75,161.78){\makebox(0,0)[cc]{2}}
\put(22.25,14){\makebox(0,0)[cc]{4}}
\put(22.25,92.5){\makebox(0,0)[cc]{4}}
\put(22.25,172.53){\makebox(0,0)[cc]{4}}
\put(22.75,27.25){\makebox(0,0)[cc]{7}}
\put(22.75,105.75){\makebox(0,0)[cc]{7}}
\put(22.75,185.78){\makebox(0,0)[cc]{7}}
\put(46.75,21.25){\makebox(0,0)[cc]{6}}
\put(46.75,99.75){\makebox(0,0)[cc]{6}}
\put(46.75,179.78){\makebox(0,0)[cc]{6}}
\put(43,35.5){\makebox(0,0)[cc]{5}}
\put(43,114){\makebox(0,0)[cc]{5}}
\put(43,194.03){\makebox(0,0)[cc]{5}}
\put(47,54.25){\makebox(0,0)[cc]{3}}
\put(47,132.75){\makebox(0,0)[cc]{3}}
\put(47,212.78){\makebox(0,0)[cc]{3}}
\put(10.25,52.5){\makebox(0,0)[cc]{$P_4(8)$}}
\put(10.25,131){\makebox(0,0)[cc]{$P_3(8)$}}
\put(10.25,211.03){\makebox(0,0)[cc]{$P_2(8)$}}
\multiput(24.75,14)(.12643678,.05603448){174}{\line(1,0){.12643678}}
\multiput(24.75,172.53)(.12643678,.05603448){174}{\line(1,0){.12643678}}
\put(92.75,2.5){\framebox(166.25,66.25)[cc]{}}
\put(92.75,81){\framebox(166.25,66.25)[cc]{}}
\put(92.75,161.03){\framebox(166.25,66.25)[cc]{}}
\put(160.75,22){\framebox(24.75,25.5)[cc]{}}
\put(160.75,100.5){\framebox(24.75,25.5)[cc]{}}
\put(160.75,180.53){\framebox(24.75,25.5)[cc]{}}
\put(110.75,22){\framebox(24.75,25.5)[cc]{}}
\put(110.75,100.5){\framebox(24.75,25.5)[cc]{}}
\put(110.75,180.53){\framebox(24.75,25.5)[cc]{}}
\put(212.25,22){\framebox(24.75,25.5)[cc]{}}
\put(212.25,100.5){\framebox(24.75,25.5)[cc]{}}
\put(212.25,180.53){\framebox(24.75,25.5)[cc]{}}
\put(89.75,54.75){\makebox(0,0)[cc]{4}}
\put(89.75,133.25){\makebox(0,0)[cc]{4}}
\put(89.75,213.28){\makebox(0,0)[cc]{4}}
\put(162.75,73){\makebox(0,0)[cc]{2}}
\put(162.75,151.5){\makebox(0,0)[cc]{2}}
\put(162.75,231.53){\makebox(0,0)[cc]{2}}
\put(262,55.25){\makebox(0,0)[cc]{4}}
\put(262,133.75){\makebox(0,0)[cc]{4}}
\put(262,213.78){\makebox(0,0)[cc]{4}}
\put(163.75,4.75){\makebox(0,0)[cc]{2}}
\put(163.75,163.28){\makebox(0,0)[cc]{2}}
\put(113.25,34.5){\makebox(0,0)[cc]{5}}
\put(113.25,113){\makebox(0,0)[cc]{5}}
\put(113.25,193.03){\makebox(0,0)[cc]{5}}
\put(133,35.75){\makebox(0,0)[cc]{5}}
\put(133,114.25){\makebox(0,0)[cc]{5}}
\put(133,194.28){\makebox(0,0)[cc]{5}}
\put(122,44.75){\makebox(0,0)[cc]{3}}
\put(122,123.25){\makebox(0,0)[cc]{3}}
\put(122,203.28){\makebox(0,0)[cc]{3}}
\put(122,24.5){\makebox(0,0)[cc]{3}}
\put(122,103){\makebox(0,0)[cc]{3}}
\put(122,183.03){\makebox(0,0)[cc]{3}}
\put(171.75,44.5){\makebox(0,0)[cc]{2}}
\put(171.75,123){\makebox(0,0)[cc]{2}}
\put(171.75,203.03){\makebox(0,0)[cc]{2}}
\put(163,35.25){\makebox(0,0)[cc]{4}}
\put(163,113.75){\makebox(0,0)[cc]{4}}
\put(163,193.78){\makebox(0,0)[cc]{4}}
\put(171.75,24.5){\makebox(0,0)[cc]{2}}
\put(171.75,103){\makebox(0,0)[cc]{2}}
\put(171.75,183.03){\makebox(0,0)[cc]{2}}
\put(182.75,36){\makebox(0,0)[cc]{4}}
\put(182.75,114.5){\makebox(0,0)[cc]{4}}
\put(182.75,194.53){\makebox(0,0)[cc]{4}}
\put(224.25,44.75){\makebox(0,0)[cc]{3}}
\put(224.25,123.25){\makebox(0,0)[cc]{3}}
\put(224.25,203.28){\makebox(0,0)[cc]{3}}
\put(223.75,24.75){\makebox(0,0)[cc]{3}}
\put(223.75,103.25){\makebox(0,0)[cc]{3}}
\put(223.75,183.28){\makebox(0,0)[cc]{3}}
\put(214.25,34.75){\makebox(0,0)[cc]{5}}
\put(214.25,113.25){\makebox(0,0)[cc]{5}}
\put(214.25,193.28){\makebox(0,0)[cc]{5}}
\put(233.5,35.75){\makebox(0,0)[cc]{5}}
\put(233.5,114.25){\makebox(0,0)[cc]{5}}
\put(233.5,194.28){\makebox(0,0)[cc]{5}}
\put(92.63,197.64){\line(1,0){18.031}}
\put(92.63,189.15){\line(1,0){18.031}}
\put(135.57,196.52){\line(1,0){25.568}}
\put(135.57,189.38){\line(1,0){25.271}}
\multiput(185.81,196.81)(4.40998,.04955){6}{\line(1,0){4.40998}}
\put(185.52,189.98){\line(1,0){26.757}}
\multiput(236.95,196.81)(3.66672,.04955){6}{\line(1,0){3.66672}}
\put(236.95,191.17){\line(1,0){22}}
\put(173.03,205.73){\line(0,1){21.703}}
\put(173.62,180.16){\line(0,-1){19.027}}
\qbezier(223.57,206.03)(172.44,236.5)(122.49,205.73)
\qbezier(121.89,180.46)(174.96,152.52)(223.87,180.46)
\put(144.79,212.57){\vector(0,-1){.12}}\multiput(145.08,220.6)(-.04955,-1.33786){6}{\line(0,-1){1.33786}}
\put(142.41,184.03){\vector(0,-1){.12}}\put(142.41,199.49){\line(0,-1){15.46}}
\put(143.6,176.3){\vector(0,1){.12}}\multiput(143.89,167.98)(-.04955,1.38741){6}{\line(0,1){1.38741}}
\put(189.68,201.87){\vector(0,1){.12}}\put(189.68,186.41){\line(0,1){15.46}}
\put(102.87,200.98){\vector(0,1){.12}}\put(102.87,185.22){\line(0,1){15.757}}
\put(240.81,186.41){\vector(0,-1){.12}}\put(240.81,202.17){\line(0,-1){15.757}}
\put(168.57,214.06){\vector(-1,0){.12}}\put(176.6,214.06){\line(-1,0){8.027}}
\put(168.57,172.73){\vector(-1,0){.12}}\multiput(176.6,172.44)(-1.33786,.04955){6}{\line(-1,0){1.33786}}
\put(96.62,200.98){\makebox(0,0)[cc]{6}}
\put(96.03,186.41){\makebox(0,0)[cc]{7}}
\put(149.54,200.08){\makebox(0,0)[cc]{7}}
\put(149.84,185.81){\makebox(0,0)[cc]{6}}
\put(199.19,200.68){\makebox(0,0)[cc]{6}}
\put(199.19,186.71){\makebox(0,0)[cc]{7}}
\put(248.25,200.98){\makebox(0,0)[cc]{7}}
\put(247.95,187.3){\makebox(0,0)[cc]{6}}
\put(186.41,217.03){\makebox(0,0)[cc]{1}}
\put(190.27,172.73){\makebox(0,0)[cc]{1}}
\put(175.71,175.11){\makebox(0,0)[cc]{0}}
\put(175.11,210.49){\makebox(0,0)[cc]{0}}
\qbezier(110.89,106.73)(99.3,106.73)(99.6,147.16)
\qbezier(123.38,100.49)(122.64,90.23)(92.76,90.08)
\qbezier(110.89,118.03)(88.74,165.75)(212.27,120.11)
\qbezier(135.57,107.62)(143.3,87.26)(154.6,87.7)
\qbezier(154.6,87.7)(280.5,85.62)(236.95,107.33)
\qbezier(173.33,125.76)(182.99,137.2)(212.27,106.43)
\qbezier(236.95,121.3)(253,108.37)(252.41,80.57)
\qbezier(122.49,125.76)(122.34,156.68)(160.84,106.73)
\qbezier(185.52,120.11)(224.91,78.19)(224.46,100.49)
\qbezier(225.35,125.76)(225.95,140.77)(258.65,139.73)
\qbezier(135.87,119.81)(143.89,119.81)(155.49,103.16)
\qbezier(155.49,103.16)(163.37,91.42)(171.84,100.49)
\qbezier(92.76,85.92)(227.73,89.04)(185.52,108.22)
\qbezier(160.84,119.22)(124.42,136.76)(258.65,143.6)
\put(237.54,83.84){\makebox(0,0)[cc]{2}}
\put(97.81,126.06){\makebox(0,0)[cc]{6}}
\put(110,132){\makebox(0,0)[cc]{1}}
\put(125.46,131.7){\makebox(0,0)[cc]{7}}
\put(140.62,113.87){\makebox(0,0)[cc]{6}}
\put(137.06,97.81){\makebox(0,0)[cc]{1}}
\put(124.27,95.73){\makebox(0,0)[cc]{7}}
\put(96.92,93.65){\makebox(0,0)[cc]{7}}
\put(128.43,90.38){\makebox(0,0)[cc]{0}}
\put(188.49,103.16){\makebox(0,0)[cc]{0}}
\put(193.25,115.65){\makebox(0,0)[cc]{7}}
\put(156.08,126.06){\makebox(0,0)[cc]{0}}
\put(220.6,138.25){\makebox(0,0)[cc]{0}}
\put(230.41,130.81){\makebox(0,0)[cc]{7}}
\put(226.84,96.92){\makebox(0,0)[cc]{7}}
\put(241.11,100.79){\makebox(0,0)[cc]{1}}
\put(207.22,126.06){\makebox(0,0)[cc]{1}}
\put(204.25,110.6){\makebox(0,0)[cc]{6}}
\put(243.49,118.92){\makebox(0,0)[cc]{6}}
\put(255.38,88.3){\makebox(0,0)[cc]{6}}
\put(109.41,127.54){\vector(3,4){.12}}\multiput(99.3,115.35)(.056157,.06771874){180}{\line(0,1){.06771874}}
\put(126.65,130.81){\vector(1,0){.12}}\put(120.11,130.81){\line(1,0){6.541}}
\put(154.89,109.11){\vector(-1,-4){.12}}\multiput(157.57,124.27)(-.0557441,-.3158831){48}{\line(0,-1){.3158831}}
\put(138.25,94.84){\vector(-1,-2){.12}}\multiput(147.16,115.95)(-.05609468,-.1327574){159}{\line(0,-1){.1327574}}
\put(102.87,83.24){\vector(0,-1){.12}}\put(102.87,93.95){\line(0,-1){10.703}}
\put(256.57,84.43){\vector(1,0){.12}}\put(249.14,84.43){\line(1,0){7.433}}
\put(218.52,97.22){\vector(-3,1){.12}}\multiput(226.25,94.24)(-.1458462,.0560947){53}{\line(-1,0){.1458462}}
\put(198.3,112.68){\vector(2,3){.12}}\multiput(191.17,101.08)(.05618301,.0912974){127}{\line(0,1){.0912974}}
\put(207.81,106.14){\vector(-1,-4){.12}}\multiput(210.19,124.57)(-.055312,-.4286677){43}{\line(0,-1){.4286677}}
\put(228.33,137.35){\vector(-4,1){.12}}\multiput(236.35,135.27)(-.2112407,.0547661){38}{\line(-1,0){.2112407}}
\multiput(111.19,30.03)(-.10405562,.056157){180}{\line(-1,0){.10405562}}
\multiput(135.57,40.14)(.11645461,-.05617223){217}{\line(1,0){.11645461}}
\multiput(185.52,42.22)(.12926164,-.05601338){207}{\line(1,0){.12926164}}
\multiput(236.95,40.14)(.14023669,-.05609468){159}{\line(1,0){.14023669}}
\put(173.92,46.97){\line(0,1){21.703}}
\put(174.52,21.41){\line(0,-1){19.325}}
\qbezier(122.79,47.27)(124.12,67.49)(160.54,36.57)
\qbezier(185.52,32.41)(223.42,-1.64)(223.87,22)
\qbezier(225.65,47.27)(225.65,69.57)(258.95,49.05)
\qbezier(122.19,21.7)(122.64,.3)(92.76,20.51)
\qbezier(110.89,40.14)(68.68,61.54)(173.92,61.54)
\qbezier(173.92,61.54)(179.87,61.1)(211.98,39.84)
\qbezier(135.57,30.92)(172.29,10.85)(177.49,10.41)
\qbezier(177.49,10.41)(271.29,6.09)(237.25,29.73)
\put(104.06,30.32){\makebox(0,0)[cc]{6}}
\put(102.87,58.87){\makebox(0,0)[cc]{1}}
\put(141.22,54.11){\makebox(0,0)[cc]{7}}
\put(193.84,29.73){\makebox(0,0)[cc]{7}}
\put(230.41,53.22){\makebox(0,0)[cc]{7}}
\put(123.38,13.08){\makebox(0,0)[cc]{7}}
\put(142.41,23.78){\makebox(0,0)[cc]{1}}
\put(243.49,29.14){\makebox(0,0)[cc]{1}}
\put(244.08,41.03){\makebox(0,0)[cc]{6}}
\put(191.17,44){\makebox(0,0)[cc]{6}}
\put(141.52,41.03){\makebox(0,0)[cc]{6}}
\put(177.49,16.35){\makebox(0,0)[cc]{0}}
\put(171.25,56.49){\makebox(0,0)[cc]{0}}
\put(207.52,46.97){\makebox(0,0)[cc]{1}}
\put(107.03,37.46){\vector(0,1){.12}}\multiput(107.33,29.14)(-.04955,1.38741){6}{\line(0,1){1.38741}}
\put(96.62,43.41){\vector(-4,-3){.12}}\multiput(104.06,48.76)(-.07742234,-.05574408){96}{\line(-1,0){.07742234}}
\put(130.22,52.03){\vector(1,0){.12}}\multiput(120.41,52.62)(.89191,-.05405){11}{\line(1,0){.89191}}
\put(149.54,28.24){\vector(1,3){.12}}\multiput(146.87,20.81)(.0557441,.1548447){48}{\line(0,1){.1548447}}
\put(155.79,25.27){\vector(0,-1){.12}}\put(155.79,45.78){\line(0,-1){20.514}}
\put(203.06,37.46){\vector(3,4){.12}}\multiput(192.65,22.89)(.05594388,.07832144){186}{\line(0,1){.07832144}}
\put(209,46.08){\vector(0,1){.12}}\put(209,29.73){\line(0,1){16.352}}
\put(241.41,20.22){\vector(-2,-3){.12}}\multiput(252.41,36.87)(-.05612329,-.08494337){196}{\line(0,-1){.08494337}}
\put(221.79,51.73){\vector(-1,0){.12}}\put(228.63,51.73){\line(-1,0){6.838}}
\put(170.95,5.95){\vector(-1,0){.12}}\put(176.6,5.95){\line(-1,0){5.649}}
\put(170.06,52.92){\vector(-1,0){.12}}\put(176.3,52.92){\line(-1,0){6.243}}
\put(126.65,18.43){\vector(1,0){.12}}\put(118.62,18.43){\line(1,0){8.027}}
\multiput(25.27,92.76)(-.10695613,.05619729){164}{\line(-1,0){.10695613}}
\multiput(46.68,194.44)(-.10097042,-.05609468){212}{\line(-1,0){.10097042}}
\end{picture}

\end{center}

To the right of each graph is shown a dissected planar surface $S_0$
with four boundary components. Each boundary component is shown as a
square whose sides are alternately labeled by different generators.
The dissected surface $S$ is obtained from $S_0$ as follows: each
side of a boundary square which is labeled $4$ or $5$ is glued to
\emph{the other edge in the same square with the same label while
reversing orientation}. This produces two boundary $2$-curves and
two boundary $3$-curves which may be identified in pairs (of the
same label).

The closed surface $S$ resulting from this construction has exactly
one dissection curve of each type $2,3,4$, and $5$.

The faithfulness of these dissections diagrams is proved as in
Section \ref{prev} because it is easy to see that the sets of
vertices satisfying (*) in the proofs for $P_2(8), P_3(8), P_4(8)$
in Section \ref{prev}, also satisfy (*) for these new dissection
diagrams.

\subsection{Kim's results} \label{kim}

We note here that recent results of Kim \cite{Kim1} somewhat overlap
with the results of this paper. Kim's work uses different
techniques. Instead of embedding surface groups directly, he embeds
right angled Artin groups which are known to contain hyperbolic
surface subgroups. More precisely, a if vertices $a,b$ in a graph
$K$ are non-adjacent, we can produce a new graph $K'$ by the {\em
co-contraction} of the pair $(a,b)$. It amounts replacing the
$(a,b)$ by one vertex connected to all vertices that were connected
to both $a$ and $b$. Kim proves that $A(K')$ is a subgroup of
$A(K)$. This allowed Kim to construct a series of graphs
$K_n=C_n\opp$, which do not contain induced subgraphs $C_n$, $n\ge
5$, and such that $A(K_n)$ contains hyperbolic  surface subgroups.
In particular, this implies that $A(P_1(6))$ contains $A(C_5)$ as a subgroup
where $C_n$ is a cycle of length $n$. (This way he answered a question
from \cite{GLR} by giving an example of a weakly chordal graph $K$
such that $A(K)$ contains a hyperbolic surface group.) Note that
this fact also follows from our results because $K_5$ is $C_5$,
$K_6$ is isomorphic to $P_1(6)$, and all $K_n, n\ge 7$ contain
isomorphic copies of $P_2(6)$.

It is easy to check (using Proposition \ref{sepsub}) that Kim's
method applies to only one of our exceptional graphs, $P_1(6)$:
applying co-contraction to any other graph $P_i(j)$, one cannot get
a graph $K$ with $A(K)$ containing a hyperbolic surface subgroup. It would be interesting to find out when a co-contraction of a pair of non-adjacent vertices in a graph $K$ avoiding $C_n, n\ge 5$, or $P_1(6)-P_4(8)$ produces a graph that does not avoid the ``forbidden" subgraphs.

\section{Diagram groups}
\label{dg}

In this section, we show that the right angled Artin group
$A(P_2(6))$ is a subgroup of a diagram group, and answer a question of Guba
and Sapir from \cite{GS}. One of the definitions of diagram groups
is the following (see \cite{GS}). Consider an alphabet $X$ and a set
${\cal S}$ of {\em cells}, each cell is a disc whose boundary is
subdivided into two directed paths (the top path and the bottom
path) labeled by positive words $u$ (the top path) and $v$ (the
bottom path) in the alphabet $X$. One can consider the cell as the
rewriting rule $u\to v$. Each cell $\pi$ is an {\em elementary
$(u,v)$-diagram} with top path labeled by $u$, bottom path labeled
by $v$, and two distinguished vertices $\iota$ and $\tau$: the
common starting and ending points of the top and bottom paths. For
every word $u$ in $X$, there exists also the {\em trivial}
$(u,u)$-diagram: a path labeled by $x$. Its top path and bottom path
coincide. There are four operations allowing to construct more
complicated diagrams from the elementary ones. These are defined as
follows.

\unitlength .5mm 
\linethickness{0.4pt}
\ifx\plotpoint\undefined\newsavebox{\plotpoint}\fi 
\begin{picture}(269.75,64.88)(0,0)
\qbezier(5.25,26.5)(19.5,50.88)(33.75,25.75)
\qbezier(43.25,26.5)(57.5,50.88)(71.75,25.75)
\qbezier(81.5,26.75)(95.75,51.13)(110,26)
\qbezier(110.25,26.25)(124.5,50.63)(138.75,25.5)
\qbezier(33.75,25.75)(20,5.88)(5.25,26.5)
\qbezier(71.75,25.75)(58,5.88)(43.25,26.5)
\qbezier(110,26)(96.25,6.13)(81.5,26.75)
\qbezier(138.75,25.5)(125,5.63)(110.25,26.25)
\put(37.5,25.5){\makebox(0,0)[cc]{$+$}}
\put(75,26){\makebox(0,0)[cc]{$=$}}
\qbezier(161.25,31.25)(184,64.88)(207.75,31)
\qbezier(223.25,26)(246,59.63)(269.75,25.75)
\qbezier(161.25,18.25)(184,-15.37)(207.75,18.5)
\qbezier(223.25,26)(246,-7.62)(269.75,25.75)
\multiput(207.75,31)(-11.5625,.0625){4}{\line(-1,0){11.5625}}
\multiput(269.75,25.75)(-11.5625,.0625){4}{\line(-1,0){11.5625}}
\multiput(207.75,18.5)(-11.5625,-.0625){4}{\line(-1,0){11.5625}}
\multiput(269.75,26.25)(-11.5625,-.0625){4}{\line(-1,0){11.5625}}
\put(161.5,31.25){\line(0,1){0}}
\put(182.75,23.75){\makebox(0,0)[cc]{$\times$}}
\put(215,25.5){\makebox(0,0)[cc]{$=$}}
\end{picture}

\begin{itemize}
\item The addition: $\Delta_1+\Delta_2$ is obtained by identifying the
distinguished vertex $\tau$ of $\Delta_1$ with the initial vertex
$\iota$ of $\Delta_1$. The top and the bottom paths of
$\Delta_1+\Delta_2$ are defined in a natural way.

\item The multiplication: If the label of the bottom path of
$\Delta_1$ coincides with the label of the top path $\Delta_2$, then
$\Delta_1\Delta_2$ is defined by identifying the bottom path of
$\Delta_1$ with the top path of $\Delta_2$.

\item The inversion: $\Delta\iv$ is obtained from $\Delta$ by switching
the top and the bottom paths of the
diagram.

\item Dipole cancelation: if $\pi$ is an $(u,v)$-cell, then we
identify $\pi\pi\iv$ with the trivial $(u,u)$-diagram. Thus we can
always replace a subdiagram $\pi\pi\iv$ of a diagram $\Delta$ by the
trivial $(u,u)$-subdiagram: the resulting diagram is {\em
equivalent} to $\Delta$.
\end{itemize}

For every word $u$, the set of all $(u,u)$-diagrams forms a group
under the product operation.

\begin{example}[\cite{GS}] The R. Thompson group $F$ is the diagram
group of all $(x,x)$-diagrams corresponding to the 1-letter alphabet
$\{x\}$ and one cell $(x^2, x)$-cell.

The wreath product $\Z\wr \Z$ is the diagram group of
$(ac,ac)$-digrams over the alphabet $\{a, b_1,$ $b_2, b_3, c\}$
corresponding to cells $ab_1\to a, b_1\to b_2, b_2\to b_3. b_3\to
b_1, b_1c\to c$.

The free group $F_2$ is the diagram group of $(a, a)$-diagrams over
the alphabet $\{a,a_1,a_2, a_3, a_4\}$ and cells $a\to a_1, a_1\to
a_2, a_2\to a, a\to a_3, a_3\to a_4, a_4\to a$.

The direct product $\Z\times \Z$ is the diagram group of
$(ab,ab)$-diagram over the alphabet $\{a, a_1, a_2, b, b_1, b_2\}$
and cells $a\to a_1, a_1\to a_2, a_2\to a, b\to b_1, b_1\to b_2,
b_2\to b$.

Many right angled Artin groups are diagram groups \cite{GS3}.
\end{example}

The class of diagram groups is closed under direct and free products
\cite{GS1}, each diagram group is linearly orderable \cite{GS3}. One
can view a diagram group as a 2-dimensional analog of a free group
(a free group is the group of 1-paths of a graph; the diagram groups
are groups of 2-paths on directed 2-complexes). The word problem in
any subgroup of a diagram group is very easy to decide. In many
important cases (including the Thompson group $F$), the conjugacy
problem in a diagram group has also an easy diagrammatic solution.

As often happens with other representation questions, given a group
$G$, it is not usually easy to find out if $G$ can be a subgroup of a
diagram group. The situation is easier for right angled Artin groups
because the pairs of commuting diagrams are easy to describe (the
description is a 2-dimensional analog of the well known description
of commuting elements in the free group) \cite{GS}. For example,
\cite[Theorem 30]{GS1} shows that if $C_n$ is a cycle of odd length
$n\ge 5$ then the right angled Artin group $A(C_n)$ cannot be
embedded into a diagram group. It is quite possible (but is not
proved yet) that the same is true for even $n\ge 6$. More
restrictions on the class of right angled Artin groups that are
diagram groups are provided in \cite{GS3}. Since groups $A(C_n)$
contain hyperbolic
 surface subgroups and until this paper there were no examples
of right angled Artin groups containing hyperbolic  surface
subgroups and not containing $A(C_n)$, $n\ge 5$, this served as one
motivation for the question of Guba and Sapir of whether a diagram
group can contain a hyperbolic  surface subgroup. Another motivation
is \cite[Theorem 9.14]{GS2} which says that if the system of
rewriting rules corresponding to the cells of a diagram group is
complete (i.e. confluent and terminating), then either the diagram
group is free or it contains a copy of $\Z\times \Z$. It is not
known if one can remove the completeness assumption in that
statement (\cite[Problem 9.15]{GS2}).

Right angled Artin groups appear naturally when one studies diagram
groups. It is proved in \cite{GS3} that every countable diagram
group is a subgroup of a certain concrete (finitely presented)
universal diagram group $U$ that is a split extension of a right
angled Artin group $A$ described below and the R. Thompson group
$F$. It is well known that the R.Thompson group $F$ does not contain
free non-Abelian subgroups \cite{BrSq85}. Hence every non-elementary
hyperbolic subgroup of a diagram group must intersect the group $A$
(which is a diagram group itself \cite{GS3}).

The infinite graph $K$ corresponding to the right angled Artin group
$A$ is defined as follows. For every subinterval $\alpha\subseteq
(0,1)$ with dyadic endpoints we assign a countable set of symbols
$K_\alpha$. The union of all $K_\alpha$ is the vertex set of $K$.
Two vertices $x\in K_\alpha$ and $y\in K_\beta$ are adjacent in $K$
if and only if the intervals $\alpha$ and $\beta$ are disjoint.

Consider a set $M$ of six subintervals of the unit interval $(0,1)$:
$\alpha_1, (0,\frac 14)$, $\alpha_2=(\frac14, \frac 12)$,
$\alpha_3=(\frac12, \frac34)$, $\alpha_4=(\frac18, \frac38)$,
$\alpha_5=(\frac38, \frac 58)$, $\alpha_6=(\frac58, \frac78)$. For
each $i=1,...,6$ pick one symbol $s_i$ from $K_{\alpha_i}$. The
subgraph of the graph $K$ from the previous paragraph spanned by the
vertices $s_1,...,s_6$ is isomorphic to the graph $P_2(6)$. The
isomorphism is the following: $1\to s_3, 2\to s_2, 3\to s_6, 4\to
s_5, 5\to s_1, 6\to s_4$.

The result of Section \ref{p16p26} now implies the following answer
to the question of Guba and Sapir.

\begin{theorem}\label{thdg} The diagram group $U$ contains the fundamental group of
a hyperbolic surface.
\end{theorem}

It is shown in \cite{GS2} that every integral homology group of any
diagram group is free Abelian. The question of whether the same is
true for subgroups of diagram groups remained open. Since the first
homology group a non-orientable surface has 2-torsion, Lemma
\ref{prisms} and Theorem \ref{thdg} show that the first homology
group of a subgroup of a diagram group can have torsion.

As an unexpected corollary of Theorem \ref{thdg} and \cite[Theorem
30]{GS1} we get

\begin{corollary} The group $A(P_2(6))$ does not contain subgroups
isomorphic to $A(C_n)$ for odd $n$.
\end{corollary}

Recall that by \cite{Kim1}, $A(P_1(6))$ contains $A(C_5)$. It would
be interesting to find out if groups $A(K)$ for any of the other
graphs $K$ from Section \ref{prev} are embedded into each other.

\section{A description of graphs without long holes and induced subgraphs $P_1(6)$,
$P_2(6)$}

\label{SC}



Here we give a proof of Theorem \ref{thCW}. The proof of the following lemma is a modification of the proof of Chudnovsky and Seymour of their result quoted in the introduction.

\begin{lemma}\label{thCW2} Let $K$ be a connected and co-connected graph that
does not contain holes of length $\ge 5$ and induced subgraphs $P_1(6), P_2(6)$. Then for every two vertices $u, v$ at distance $2$ in $K$ and every co-component $W$ of $C(\{u,v\})$ one of the following three conditions holds:

\begin{itemize}
\item[(1)] $\{u\}\star W$ is a separator;
\item[(2)] $\{v\}\star W$ is a separator;
\item[(3)] $W'=C(W)\setminus \{u,v\}$ is not empty and $W\star W'$ separates $u$
and $v$.
\end{itemize}
\end{lemma}

\proof Suppose that $K$ is connected, co-connected, does not contain long holes and copies of $P_1(6), P_2(6)$. If $|K|=1$ there is nothing to prove. Let $|K|>1$.

We observe that since an anti-hole of length five is also a hole of length five, $P_1(6)$
is isomorphic to the anti-hole of length six, and $P_2(6)$ is the complement of a
path of length 5, it follows that $K$ contains no anti-hole of
length at least 5.

{\bf Step 1.} Since $K$ is co-connected, it is not a clique. Hence it has two non-adjacent vertices $u,v$ at distance 2, hence $C(\{u,v\})$ is not empty. Take any two such vertices $u, v$. We need to show that these vertices satisfy one of the Conditions (1), (2), (3).

Let $W$ be any co-component of the subgraph spanned by $C(\{u,v\})$ (i.e. a connected component of the complement graph $C(\{u,v\})\opp$).

{\bf Step 2.} Suppose that $u, v$ are disconnected by $W$.

Let $U$ be the connected component of $K\setminus (W \cup \{v\})$ containing $u$,
and let $V$ be the connected component of $K \setminus (W \cup \{u\})$ containing $v$. Then
$U$ and $V$ are disjoint. Let $X=K^0 \setminus (U \cup V \cup W)$.
Then there are no edges between any pair of the sets $U,V,X$. If $U
\cup X \neq \{u\}$, then $\{u\}\star W$ separates $K$, hence $\{u,v\}$ satisfies Condition (1). Hence we can assume $U=\{u\}$, $X=\emptyset$. Similarly, $V=\{v\}$. Thus $K=W \star \{u,v\}$ and $K\opp$ is not connected, a contradiction. Hence we can assume that $u, v$ are connected in $K\setminus W$.

{\bf Step 3.} Let us prove Condition 3. By contradiction, suppose that there exists a path $p$ from $u$ to $v$ in $K\setminus W$
contains no vertex  of $W'=C(W)\setminus \{u,v\}$.

Let $p$ be a shortest path from $u$ to $v$ in $K\setminus W$.
If $p$ has length 2,
then the vertex $w$ on $p$ that is distinct from $u, v$ is in $C(\{u,v\})$.
Since $W$ is a co-component of $C(\{u,v\})$, $w$ is in $C(W)$.

So we may assume that $p$ is an induced path of length at least three.
Let $w_1, w_2$ be the second and the third vertices of $p$ (counting from $u$ to $v$). Since $w_1,w_2$ are not in $C(W)$ there exist two vertices $z_1, z_2$ in $W$ such that $z_1$ is not connected to $w_1$ and $z_2$ is not connected to $w_2$ (it could happen that $z_1=z_2$). Since the subgraph spanned by $W$ in $K\opp$ is connected, there exists an anti-path $q$ connecting $w_1, w_2$ such that $q\setminus \{w_1,w_2\}$ is in $W$.

Since $p$ is a shortest path, any $z$ in $W$ together with
any subpath of $p$ of which only the first and the last vertices are
adjacent to $z$ form a hole. Thus if some $z\in W$ is
not adjacent to both $w_1$ and $w_2$, then $z,u, w_1, w_2$ lie in a hole of length at least 5 (the smallest hole in the subgraph spanned by $z$ and $p$ containing $z,u, w_1, w_2$).

Therefore there is no induced anti-path of length 2
from $w_1$ to $w_2$, passing through $W$. In particular, $z_1\ne z_2$, $z_1$ is adjacent to $w_2$, $z_2$ is adjacent to $w_1$.

\unitlength 1mm 
\linethickness{0.4pt}
\ifx\plotpoint\undefined\newsavebox{\plotpoint}\fi 
\begin{picture}(97,41.25)(0,0)
\multiput(17,15.5)(.098214286,-.033730159){252}{\line(1,0){.098214286}}
\put(41.75,7){\line(1,0){29}}
\multiput(54.5,34)(-.069316081,-.033733826){541}{\line(-1,0){.069316081}}
\put(76,7){\makebox(0,0)[cc]{\dots}}
\multiput(94.5,8.5)(-.0532407407,.0337301587){756}{\line(-1,0){.0532407407}}
\put(69.75,2.75){\makebox(0,0)[cc]{$p$}}
\put(14.5,15.75){\makebox(0,0)[cc]{$u$}}
\put(54.25,37){\makebox(0,0)[cc]{$z_1$}}
\put(41,4.25){\makebox(0,0)[cc]{$w_1$}}
\put(69,9.5){\makebox(0,0)[cc]{$w_2$}}
\put(97,7.25){\makebox(0,0)[cc]{$v$}}
\multiput(17.25,15.75)(.112292052,.033733826){541}{\line(1,0){.112292052}}
\multiput(78,34)(.033673469,-.052040816){490}{\line(0,-1){.052040816}}
\put(80.25,35){\makebox(0,0)[cc]{$z_2$}}
\multiput(77.93,33.68)(-.046086,-.03346){33}{\line(-1,0){.046086}}
\multiput(74.89,31.47)(-.046086,-.03346){33}{\line(-1,0){.046086}}
\multiput(71.85,29.26)(-.046086,-.03346){33}{\line(-1,0){.046086}}
\multiput(68.8,27.05)(-.046086,-.03346){33}{\line(-1,0){.046086}}
\multiput(65.76,24.85)(-.046086,-.03346){33}{\line(-1,0){.046086}}
\multiput(62.72,22.64)(-.046086,-.03346){33}{\line(-1,0){.046086}}
\multiput(59.68,20.43)(-.046086,-.03346){33}{\line(-1,0){.046086}}
\multiput(56.64,18.22)(-.046086,-.03346){33}{\line(-1,0){.046086}}
\multiput(53.6,16.01)(-.046086,-.03346){33}{\line(-1,0){.046086}}
\multiput(50.55,13.8)(-.046086,-.03346){33}{\line(-1,0){.046086}}
\multiput(47.51,11.6)(-.046086,-.03346){33}{\line(-1,0){.046086}}
\multiput(44.47,9.39)(-.046086,-.03346){33}{\line(-1,0){.046086}}
\multiput(54.18,33.93)(.032961,-.053753){29}{\line(0,-1){.053753}}
\multiput(56.09,30.81)(.032961,-.053753){29}{\line(0,-1){.053753}}
\multiput(58,27.69)(.032961,-.053753){29}{\line(0,-1){.053753}}
\multiput(59.92,24.58)(.032961,-.053753){29}{\line(0,-1){.053753}}
\multiput(61.83,21.46)(.032961,-.053753){29}{\line(0,-1){.053753}}
\multiput(63.74,18.34)(.032961,-.053753){29}{\line(0,-1){.053753}}
\multiput(65.65,15.22)(.032961,-.053753){29}{\line(0,-1){.053753}}
\multiput(67.56,12.11)(.032961,-.053753){29}{\line(0,-1){.053753}}
\multiput(69.47,8.99)(.032961,-.053753){29}{\line(0,-1){.053753}}
\put(65.5,34.5){\makebox(0,0)[cc]{$\dots$}}
\put(65.25,37.25){\makebox(0,0)[cc]{$q$}}
\end{picture}

Now, if $w_2$ is non-adjacent to $v$, then the subgraph of $K$ spanned by the vertices of the anti-loop $q\cup \{(w_2,v), (v, w_1)\}$ contains an anti-hole of length at least $5$ (the smallest anti-hole containing the anti-path $z_1, w_1, v, w_2, z_2$), a contradiction. This proves that $w_2$ is adjacent to $v$. Since $z_1$ is adjacent to $w_2$ but not to $w_1$ and $z_2$ is adjacent to $w_1$ but not to $w_2$, we can find two consecutive vertices $z, z'\in W$ on $q$ such that $z$ is adjacent to $w_2$ but not to $w_1$, and $z'$ is adjacent to $w_1$ (and possibly also to $w_1$). But then the subgraph induced by $\{z,z', u, w_1, w_2, v\}$ is  the prism $P_1(6)$ or the prism with diagonal $P_2(6)$ (depending on whether there exists an edge $(z', w_2)$ or not), a contradiction. This proves that every path from $u$ to $v$ avoiding $W$ passes through $W'$. This and Step 2 show that $W'$ is not empty. Hence $u, v$ are in different connected components of $K\setminus (W\star W')$, and $u, v$ satisfies Condition (3).
This proves the lemma.
\endproof

Finally let us prove Theorem \ref{thCW}. Recall its formulation.

\begin{theorem}\label{thCW3} A graph $K$
does not contain holes of length $\ge 5$ and induced subgraphs $P_1(6)$, $P_2(6)$ if and only if for every two vertices $u,v$ at distance 2 and every co-component $W$ of $C(\{u,v\})$, the set $W \cup (C(W)\setminus \{u,v\})$ separates $u$ from $v$.
\end{theorem}

\proof The ``only if part" follows from Step 3 of the proof of Lemma \ref{thCW2}. Indeed, we can assume that two vertices $u, v$ at distance 2 in $K$ are not separated by $W$ (otherwise they would be separated by $W\cup (C(W)\setminus \{u,v\})$ as well. But then Step 3 of the proof of Lemma \ref{thCW2} gives that $u, v$ are separated by
$W\cup (C(W)\setminus \{u,v\})$.

The ``if" part follows from the fact that both long holes and $P_1(6), P_2(6)$ contains pairs of vertices $u,v$ and co-components $W$ of $C(\{u,v\}$ such that $W\cup (C(W)\setminus \{u,v\})$ does not separate $u,v$, and that the condition of the theorem is obviously hereditary for induced subgraphs.
\endproof

\section{The proof of Theorem \ref{8vertex}}
\label{8}

Here we present the description of a computer based proof of Theorem \ref{8vertex}: if $|K^0|\le 8$ and $A(K)$ contains a hyperbolic surface subgroup, then $K$ contains one of our ``forbidden" induced subgraphs $C_n, n=5,6,7,8$, $P_1(6)-P_4(8)$.

First using independently written programs written in C by the first
author and in Maple by the third author, we checked all graphs with
at most 7 vertices and found out that each of them either can be
reduced to a tree-like graph by the moves that we considered in the
previous sections (and so the corresponding right angled Artin group
does not contain a non-abelian surface subgroup), or contains a
fully embedded circuit of length at least 5, or contains one of the
graphs $P_1(6)-P_2(7)$. An earlier computation was performed in MAGMA with the help
of Marston Conder. This program formed the basis for the later C++ program.
We wish to thank Marston Conder for his enthusiastic contribution to
this project. Several of our ``forbidden" subgraphs first appeared there.

The hardest case is, of course, when $|K^0|=8$. So we describe it in some details. The (Maple) program was created by the third author. The program successively eliminated 8-vertex graphs, first removing graphs containing ``forbidden" subgraphs, and then applying more and more complicated reduction rules.

The procedure consists of several steps. At each step, we take the first of the remaining graphs and find the simplest reduction move that eliminates it. Then we eliminate all other graphs using that reduction move, etc.

\bigskip

{\bf Step 1}. It is well known that the total number of 8-vertex graphs up to isomorphism is 12,346 (see the Brendan McKay's Web site http://cs.anu.edu.au/$\sim$bdm/data/graphs.html . We first eliminate graphs containing $C_n, n\ge 5$, $P_1(6)-P_4(8)$, and also graphs $K$ that satisfy one of the following conditions:

\begin{itemize}
\item[$(R_1)$] $K$ is disconnected;
\item[$(R_2)$] $K\opp$ is disconnected;
\item[$(R_3)$] $K$ decomposes non-trivially as an almost join;
\item[$(R_4)$] $K$ contains a splitting subset $X$ such that $K=K_1\cup_{C(X)} ...\cup_{C(X)} K_m$, and $X\cap K_i$ is dense in $K_i$ relative to $\emptyset$ for $i=1,...,m$ (applying Lemma \ref{sepsub});
\item[$(R_5)$] $K$ contains a pair of adjacent vertices $x,y$ with adjacent links $\Lk(x), \Lk(y)$ (applying Lemma \ref{edgered}) and such that the graph $K\setminus \{(x,y)\}$ does not contain forbidden subgraphs and has been eliminated already.
\end{itemize}

There are 67 graphs surviving this step. From now on, we call a graph considered on Step $i$ {\em excluded} if if does not contain the ``forbidden subgraphs" and has been excluded on steps $\le i-1$.

\bigskip

{\bf Step 2.} One of the remaining graphs is the graph $K$ on the next picture:

\begin{center}
\unitlength .6mm 
\linethickness{0.4pt}
\ifx\plotpoint\undefined\newsavebox{\plotpoint}\fi 
\begin{picture}(99,97.25)(0,0)
\put(20.5,37){\line(0,1){22.75}}
\put(20.5,59.75){\line(1,0){49}}
\put(69.5,59.75){\line(0,-1){21}}
\put(20.25,37.25){\line(1,0){49.25}}
\put(69.5,37.25){\line(0,1){4.25}}
\multiput(20.75,59.5)(.07388974,-.033690658){653}{\line(1,0){.07388974}}
\multiput(69.25,59.5)(-.074272588,-.033690658){653}{\line(-1,0){.074272588}}
\multiput(80,28)(-.033653846,.101762821){312}{\line(0,1){.101762821}}
\multiput(80,28.25)(-.039325843,.033707865){267}{\line(-1,0){.039325843}}
\multiput(80,28.5)(-.223782772,.033707865){267}{\line(-1,0){.223782772}}
\put(71.75,38.75){\makebox(0,0)[cc]{1}}
\put(72.25,60.5){\makebox(0,0)[cc]{2}}
\put(16.25,60.5){\makebox(0,0)[cc]{3}}
\put(16.75,35.75){\makebox(0,0)[cc]{4}}
\put(82.75,28.25){\makebox(0,0)[cc]{5}}
\multiput(69.25,37.25)(.109126984,.033730159){252}{\line(1,0){.109126984}}
\multiput(96.75,45.75)(-.065662651,.03373494){415}{\line(-1,0){.065662651}}
\put(99,45.75){\makebox(0,0)[cc]{6}}
\put(69.5,59.75){\line(0,1){20.75}}
\multiput(69.5,80.5)(.0337009804,-.0422794118){816}{\line(0,-1){.0422794118}}
\qbezier(69.5,80.5)(-23,97.25)(20.5,37)
\put(72.25,81.75){\makebox(0,0)[cc]{7}}
\multiput(77.75,14)(-.033730159,.092261905){252}{\line(0,1){.092261905}}
\multiput(77.75,14.25)(.033712785,.055604203){571}{\line(0,1){.055604203}}
\multiput(78,14.5)(-.0843108504,.0337243402){682}{\line(-1,0){.0843108504}}
\put(78.5,12.75){\makebox(0,0)[cc]{8}}
\end{picture}
\end{center}

Note that this graph satisfies the following property

\begin{itemize}
\item[($R_6$)] There exists a set of vertices $X$ (in our case, $X=\{3,5,7,8\}$) that is stable, and a vertex $x$ (say, $x=5$) in $X$ such  that the commutator $[\Lk(x),X']$ of the link of $x$ and the complement $X'$ of $X$ is dense in $K$ relative to $\emptyset$ (in fact for this graph, the commutator $\{4,6\}$ is even nuclear in $K$ relative to $\emptyset$).
\end{itemize}

If a graph $K$ satisfies $(R_6)$, then there could not be a faithful $K$-dissection diagram on a hyperbolic surface $S$ containing $x$-curves for all $x\in K^0$ and so the graph can be reduced (some of the vertices could be removed). Indeed, if such a faithful dissection diagram exists, we can cut the surface along $X$-curves (which are disjoint since $X$ is stable), producing a surface $(S', \partial)$ with non-Abelian fundamental group and content in $X'$. By Lemma \ref{BasicLemma1}, we can assume that one of the boundary curves of $S'$ is an $x$-curve $\alpha$. Its content is in $\Lk(x)$. Since the fundamental group of $S'$ is non-Abelian, there exists an essential closed curve $\beta$ in $S'$ intersecting $\alpha$.  By Lemma \ref{commeff} the essential closed curve $[\alpha,\beta]$ has effective content inside $[\Lk(x), X]$ which is dense, contradicting Lemma \ref{power}.

Only $50$ of the $67$ graphs do not have property $(R_6)$.

\bigskip

{\bf Step 3.} One of these 50 graphs is presented on the following picture.

\begin{center}
\unitlength .6mm 
\linethickness{0.4pt}
\ifx\plotpoint\undefined\newsavebox{\plotpoint}\fi 
\begin{picture}(107.5,64.75)(0,0)
\put(107.5,34.25){\makebox(0,0)[cc]{1}}
\put(87.75,54){\makebox(0,0)[cc]{2}}
\put(63.5,64.75){\makebox(0,0)[cc]{3}}
\put(32.5,59.5){\makebox(0,0)[cc]{4}}
\put(14.75,38.5){\makebox(0,0)[cc]{5}}
\put(36.25,15){\makebox(0,0)[cc]{6}}
\put(67.5,9.25){\makebox(0,0)[cc]{7}}
\put(91.75,18){\makebox(0,0)[cc]{8}}
\multiput(105.5,36)(-.03515625,.033691406){512}{\line(-1,0){.03515625}}
\multiput(105.25,35.75)(-.051183432,.0337278107){845}{\line(-1,0){.051183432}}
\multiput(105,36)(-.1044721408,.0337243402){682}{\line(-1,0){.1044721408}}
\multiput(105.25,35.75)(-.119154676,-.033723022){556}{\line(-1,0){.119154676}}
\multiput(105.25,35.75)(-.0542521994,-.0337243402){682}{\line(-1,0){.0542521994}}
\multiput(105,35.5)(-.033653846,-.038461538){364}{\line(0,-1){.038461538}}
\multiput(87.25,53.5)(-.32469512,.03353659){164}{\line(-1,0){.32469512}}
\multiput(87,53.75)(-.16607565,-.033687943){423}{\line(-1,0){.16607565}}
\multiput(86.75,53.75)(-.0439001848,-.0337338262){1082}{\line(-1,0){.0439001848}}
\multiput(62.75,63.5)(-.21268657,-.03358209){134}{\line(-1,0){.21268657}}
\multiput(62.5,63.5)(-.0646306818,-.0337357955){704}{\line(-1,0){.0646306818}}
\multiput(62.5,63.75)(-.0336956522,-.0670289855){690}{\line(0,-1){.0670289855}}
\multiput(62.25,63.5)(.03362573,-.29678363){171}{\line(0,-1){.29678363}}
\multiput(62.25,63.25)(.0337171053,-.0452302632){912}{\line(0,-1){.0452302632}}
\multiput(34.25,59)(-.03371869,-.036608863){519}{\line(0,-1){.036608863}}
\multiput(34,59)(.0337361531,-.0460725076){993}{\line(0,-1){.0460725076}}
\multiput(34,59)(.0535551504,-.03372835){1097}{\line(1,0){.0535551504}}
\multiput(16.75,40.25)(.0337331334,-.0344827586){667}{\line(0,-1){.0344827586}}
\end{picture}

\end{center}

Note that this graph satisfies the following property (a generalization of Property $(R_5)$).

\begin{itemize}
\item[$(R_7)$] There are two adjacent vertices $x$ and $y$ (vertices 2 and 5 in that graph) such that the commutator $[Lk(x), Lk(y)]$ is dense in the graph relative to $\emptyset$.
\end{itemize}

In that case, we can simplify the graph by removing the edge $(x,y)$, because if a dissection diagram contains an $x$-curve intersecting a $y$-curve, then the commutator of these curves will be (Lemma \ref{commeff}) an essential closed curve whose effective content is dense in our graph relative to $\emptyset$ contradicting Lemma \ref{power}.  The resulting graph does not contain forbidden subgraphs and has been eliminated on the previous steps. Hence that graph can be eliminated as well. Fifteen other graphs among remaining 50 satisfy property $(R_7)$ as well. So only 35 graphs survive Step 3.

\bigskip

{\bf Step 4.} The following graph (denoted again by $K$) is one of the 35 remaining graphs.

\begin{center}
\unitlength .6 mm 
\linethickness{0.4pt}
\ifx\plotpoint\undefined\newsavebox{\plotpoint}\fi 
\begin{picture}(107.5,64.75)(0,0)
\put(107.5,34.25){\makebox(0,0)[cc]{1}}
\put(87.75,54){\makebox(0,0)[cc]{2}}
\put(63.5,64.75){\makebox(0,0)[cc]{3}}
\put(32.5,59.5){\makebox(0,0)[cc]{4}}
\put(14.75,38.5){\makebox(0,0)[cc]{5}}
\put(36.25,15){\makebox(0,0)[cc]{6}}
\put(67.5,9.25){\makebox(0,0)[cc]{7}}
\put(91.75,18){\makebox(0,0)[cc]{8}}
\multiput(105.25,35.75)(-.051183432,.0337278107){845}{\line(-1,0){.051183432}}
\multiput(105,36)(-.1044721408,.0337243402){682}{\line(-1,0){.1044721408}}
\multiput(105.25,35.75)(-.119154676,-.033723022){556}{\line(-1,0){.119154676}}
\multiput(105.25,35.75)(-.0542521994,-.0337243402){682}{\line(-1,0){.0542521994}}
\multiput(105,35.5)(-.033653846,-.038461538){364}{\line(0,-1){.038461538}}
\multiput(87.25,53.5)(-.32469512,.03353659){164}{\line(-1,0){.32469512}}
\multiput(87,53.75)(-.16607565,-.033687943){423}{\line(-1,0){.16607565}}
\multiput(62.75,63.5)(-.21268657,-.03358209){134}{\line(-1,0){.21268657}}
\multiput(62.5,63.5)(-.0646306818,-.0337357955){704}{\line(-1,0){.0646306818}}
\multiput(62.5,63.75)(-.0336956522,-.0670289855){690}{\line(0,-1){.0670289855}}
\multiput(34.25,59)(-.03371869,-.036608863){519}{\line(0,-1){.036608863}}
\multiput(34,59)(.0337361531,-.0460725076){993}{\line(0,-1){.0460725076}}
\multiput(34,59)(.0535551504,-.03372835){1097}{\line(1,0){.0535551504}}
\multiput(105,36)(-.78571429,.03348214){112}{\line(-1,0){.78571429}}
\multiput(86.75,53.75)(-.081649832,.033670034){297}{\line(-1,0){.081649832}}
\multiput(86.5,54)(-.033723022,-.074190647){556}{\line(0,-1){.074190647}}
\multiput(16.75,39.75)(.033690658,-.034073507){653}{\line(0,-1){.034073507}}
\multiput(37.25,19.25)(.0333333,-.0333333){60}{\line(0,-1){.0333333}}
\multiput(39.25,17.25)(.21455224,-.03358209){134}{\line(1,0){.21455224}}
\multiput(68,12.75)(.093632959,.033707865){267}{\line(1,0){.093632959}}
\end{picture}
\end{center}

Graph $K$ satisfies the following property

\begin{itemize}
\item[$(R_8)$] There exists a stable set $X$, and a vertex $x\in X$ such that $L=[\Lk(x), X']$ (here $X'$ is the complement of $X$) satisfies the following condition:
\begin{itemize}

\item[(*)] $K$ is an amalgam of $K_1,...,K_m$ over $C(L)$ and
    $L\cap K_i$ is dense in $U_i=(K_i\setminus L)\cup\Lk(K_i\setminus L)$ relative to $\emptyset$ for every $i$.
\end{itemize}
\end{itemize}

For the graph $K$, $X=\{2,6,8\}$, $x=2$, $L=[\Lk(x),X']=\{3,5,7\}$, $C(L)=\{1,2,4,6\}$, $m=2$, $K_1=\{3,5\}\cup C(L)$, $K_2=\{7,8\}\cup C(L)$.

Note that if a graph $G$ satisfies ($R_8$),  then it can be simplified. Indeed, if $\Delta$ is a faithful $G$-dissection diagram on a surface $S$ containing $y$-curves for all $y\in K^0$, then as in Step 2, there exists an non-null-homotopic closed curve $\ga$ with effective content inside $L=[\Lk(x),X']$. By Lemma \ref{BasicLemma1}, there exists a collection of sets of simple closed curves $B_1,...,B_m$ satisfying the conditions of the lemma. The curve $\ga$ cannot intersect a curve $\alpha$ from $B_i$ because $\cont(\alpha)$ is adjacent to $\effcont(\ga)$ (we use Lemma \ref{comm}). We can assume that $\ga$ is in a
connected component $S'$ in $S\setminus \cup B_j$ with non-Abelian $\pi_1(S')$.
Then the content of $S'$ cannot be inside $C(L)$ (again use Lemma \ref{comm}), so $\cont(S')$ is inside $U_i=(K_i\setminus C(L))\cup (\Lk(K_i\setminus C(L)))$ for some $i=1,...,m$. Since $L\cap U_i$ is dense in $U_i$ relative to $\emptyset$, and $C(L)$ is adjacent to $L$, the effective content of $\ga$ is dense in $U_i$ relative to $L\cap U_i\supset \cont(\partial S')$. That contradicts Lemma \ref{power}.

Of the 35 remaining graphs, 11 satisfy Property $(R_8)$ and can be eliminated. It leaves 24 graphs to consider.

\bigskip

{\bf Step 5.} One of these 24 graphs is the following (we denote it by $K$, as usual).

This graph satisfies the following condition.

\begin{itemize}
\item[$(R_9)$] There exists a stable set $X$ and a vertex $x\in X$ such that
$\Lk(x)$ is dense in $K\setminus X$ relative to the collection of sets $\{\Lk(y),y\in X\}$.
\end{itemize}

For the graph $K$, $X=\{4,8\}$, $x=8$. Property $(R_9)$ implies that the graph can be reduced. Indeed, cut the surface by $X$-curves and apply Lemma \ref{power} to the resulting surface $S'$ and its boundary $x$-curve.

Four graphs among our 24 satisfy $(R_9)$, 20 graphs remain after Step 5.


\begin{center}
\unitlength .6mm 
\linethickness{0.4pt}
\ifx\plotpoint\undefined\newsavebox{\plotpoint}\fi 
\begin{picture}(107.5,64.75)(0,0)
\put(107.5,34.25){\makebox(0,0)[cc]{1}}
\put(87.75,54){\makebox(0,0)[cc]{2}}
\put(63.5,64.75){\makebox(0,0)[cc]{3}}
\put(32.5,59.5){\makebox(0,0)[cc]{4}}
\put(14.75,38.5){\makebox(0,0)[cc]{5}}
\put(36.25,15){\makebox(0,0)[cc]{6}}
\put(67.5,9.25){\makebox(0,0)[cc]{7}}
\put(91.75,18){\makebox(0,0)[cc]{8}}
\multiput(105.25,35.75)(-.051183432,.0337278107){845}{\line(-1,0){.051183432}}
\multiput(105,36)(-.1044721408,.0337243402){682}{\line(-1,0){.1044721408}}
\multiput(105.25,35.75)(-.119154676,-.033723022){556}{\line(-1,0){.119154676}}
\multiput(105.25,35.75)(-.0542521994,-.0337243402){682}{\line(-1,0){.0542521994}}
\multiput(87.25,53.5)(-.32469512,.03353659){164}{\line(-1,0){.32469512}}
\multiput(87,53.75)(-.16607565,-.033687943){423}{\line(-1,0){.16607565}}
\multiput(62.75,63.5)(-.21268657,-.03358209){134}{\line(-1,0){.21268657}}
\multiput(62.5,63.5)(-.0646306818,-.0337357955){704}{\line(-1,0){.0646306818}}
\multiput(62.5,63.75)(-.0336956522,-.0670289855){690}{\line(0,-1){.0670289855}}
\multiput(34.25,59)(-.03371869,-.036608863){519}{\line(0,-1){.036608863}}
\multiput(105,36)(-.78571429,.03348214){112}{\line(-1,0){.78571429}}
\multiput(86.75,53.75)(-.081649832,.033670034){297}{\line(-1,0){.081649832}}
\multiput(86.5,54)(-.033723022,-.074190647){556}{\line(0,-1){.074190647}}
\multiput(39.25,17.25)(.21455224,-.03358209){134}{\line(1,0){.21455224}}
\multiput(68,12.75)(.093632959,.033707865){267}{\line(1,0){.093632959}}
\multiput(86.25,54)(-.0428899083,-.0337155963){1090}{\line(-1,0){.0428899083}}
\multiput(17.25,40)(.0625,-.0337252475){808}{\line(1,0){.0625}}
\multiput(39.25,17.25)(.38297872,.03368794){141}{\line(1,0){.38297872}}
\end{picture}

\end{center}

{\bf Step 6.} One of the 20 remaining graphs (again $K$) is on the following picture.

\begin{center}
\unitlength .6mm 
\linethickness{0.4pt}
\ifx\plotpoint\undefined\newsavebox{\plotpoint}\fi 
\begin{picture}(107.5,64.75)(0,0)
\put(107.5,34.25){\makebox(0,0)[cc]{1}}
\put(87.75,54){\makebox(0,0)[cc]{2}}
\put(63.5,64.75){\makebox(0,0)[cc]{3}}
\put(32.5,59.5){\makebox(0,0)[cc]{4}}
\put(14.75,38.5){\makebox(0,0)[cc]{5}}
\put(36.25,15){\makebox(0,0)[cc]{6}}
\put(67.5,9.25){\makebox(0,0)[cc]{7}}
\put(91.75,18){\makebox(0,0)[cc]{8}}
\multiput(105.25,35.75)(-.051183432,.0337278107){845}{\line(-1,0){.051183432}}
\multiput(105,36)(-.1044721408,.0337243402){682}{\line(-1,0){.1044721408}}
\multiput(105.25,35.75)(-.119154676,-.033723022){556}{\line(-1,0){.119154676}}
\multiput(105.25,35.75)(-.0542521994,-.0337243402){682}{\line(-1,0){.0542521994}}
\multiput(87.25,53.5)(-.32469512,.03353659){164}{\line(-1,0){.32469512}}
\multiput(87,53.75)(-.16607565,-.033687943){423}{\line(-1,0){.16607565}}
\multiput(62.75,63.5)(-.21268657,-.03358209){134}{\line(-1,0){.21268657}}
\multiput(62.5,63.5)(-.0646306818,-.0337357955){704}{\line(-1,0){.0646306818}}
\multiput(62.5,63.75)(-.0336956522,-.0670289855){690}{\line(0,-1){.0670289855}}
\multiput(34.25,59)(-.03371869,-.036608863){519}{\line(0,-1){.036608863}}
\multiput(34,59)(.0337361531,-.0460725076){993}{\line(0,-1){.0460725076}}
\multiput(34,59)(.0535551504,-.03372835){1097}{\line(1,0){.0535551504}}
\multiput(105,36)(-.78571429,.03348214){112}{\line(-1,0){.78571429}}
\multiput(86.75,53.75)(-.081649832,.033670034){297}{\line(-1,0){.081649832}}
\multiput(86.5,54)(-.033723022,-.074190647){556}{\line(0,-1){.074190647}}
\multiput(16.75,39.75)(.033690658,-.034073507){653}{\line(0,-1){.034073507}}
\multiput(39.25,17.25)(.21455224,-.03358209){134}{\line(1,0){.21455224}}
\multiput(68,12.75)(.093632959,.033707865){267}{\line(1,0){.093632959}}
\multiput(86.75,54)(.03360215,-.17204301){186}{\line(0,-1){.17204301}}
\end{picture}
\end{center}

This graph satisfies the following property.

\begin{itemize}
\item[$(R_{10})$] There exist two adjacent vertices $x,y$ such that the commutator $L=[Lk(x), Lk(y)]$ satisfies Condition (*) from Step 4. In addition, the graph $K$ with edge $(x,y)$ removed is excluded.
\end{itemize}
For the graph $K$, we can take $x=6, y=1$. Then $L=\{3,5,7\}$, $C(L)=\{1,2,4,6\}$, $K_1=\{3,5\}\cup C(L)$, $K_2=\{7,8\}\cup C(L)$.

A graph satisfying Property $(R_{10})$ can be reduced by removing the edge $(x,y)$. Indeed, if a $K$-dissection diagram has $x$-curve $\alpha$ intersecting $y$-curve $\beta$, then taking the commutator $\ga=[\alpha,\beta]$ we obtain (by Lemma \ref{commeff}) a closed essential curve with effective content $L$, and we can continue as in Step 4.

Of 20 remaining graphs 6 satisfy $(R_{10})$ and 14 graphs remain.

\bigskip

{\bf Step 7.} Two of the remaining 14 graphs, $K_1$, $K_2$, are on the following picture.

\pagebreak

\begin{center}
\begin{figure}[htbp]
\unitlength .6 mm 
\linethickness{0.4pt}
\ifx\plotpoint\undefined\newsavebox{\plotpoint}\fi 
\begin{picture}(216,64.75)(0,0)
\put(107.5,34.25){\makebox(0,0)[cc]{1}}
\put(216,34.25){\makebox(0,0)[cc]{1}}
\put(87.75,54){\makebox(0,0)[cc]{2}}
\put(196.25,54){\makebox(0,0)[cc]{2}}
\put(63.5,64.75){\makebox(0,0)[cc]{3}}
\put(172,64.75){\makebox(0,0)[cc]{3}}
\put(32.5,59.5){\makebox(0,0)[cc]{4}}
\put(141,59.5){\makebox(0,0)[cc]{4}}
\put(14.75,38.5){\makebox(0,0)[cc]{5}}
\put(123.25,38.5){\makebox(0,0)[cc]{5}}
\put(36.25,15){\makebox(0,0)[cc]{6}}
\put(144.75,15){\makebox(0,0)[cc]{6}}
\put(67.5,9.25){\makebox(0,0)[cc]{7}}
\put(176,9.25){\makebox(0,0)[cc]{7}}
\put(91.75,18){\makebox(0,0)[cc]{8}}
\put(200.25,18){\makebox(0,0)[cc]{8}}
\multiput(105.25,35.75)(-.051183432,.0337278107){845}{\line(-1,0){.051183432}}
\multiput(213.75,35.75)(-.051183432,.0337278107){845}{\line(-1,0){.051183432}}
\multiput(105,36)(-.1044721408,.0337243402){682}{\line(-1,0){.1044721408}}
\multiput(105.25,35.75)(-.119154676,-.033723022){556}{\line(-1,0){.119154676}}
\multiput(213.75,35.75)(-.119154676,-.033723022){556}{\line(-1,0){.119154676}}
\multiput(105.25,35.75)(-.0542521994,-.0337243402){682}{\line(-1,0){.0542521994}}
\multiput(87.25,53.5)(-.32469512,.03353659){164}{\line(-1,0){.32469512}}
\multiput(195.75,53.5)(-.32469512,.03353659){164}{\line(-1,0){.32469512}}
\multiput(87,53.75)(-.16607565,-.033687943){423}{\line(-1,0){.16607565}}
\multiput(195.5,53.75)(-.16607565,-.033687943){423}{\line(-1,0){.16607565}}
\multiput(62.75,63.5)(-.21268657,-.03358209){134}{\line(-1,0){.21268657}}
\multiput(171.25,63.5)(-.21268657,-.03358209){134}{\line(-1,0){.21268657}}
\multiput(62.5,63.5)(-.0646306818,-.0337357955){704}{\line(-1,0){.0646306818}}
\multiput(171,63.5)(-.0646306818,-.0337357955){704}{\line(-1,0){.0646306818}}
\multiput(62.5,63.75)(-.0336956522,-.0670289855){690}{\line(0,-1){.0670289855}}
\multiput(171,63.75)(-.0336956522,-.0670289855){690}{\line(0,-1){.0670289855}}
\multiput(34.25,59)(-.03371869,-.036608863){519}{\line(0,-1){.036608863}}
\multiput(142.75,59)(-.03371869,-.036608863){519}{\line(0,-1){.036608863}}
\multiput(105,36)(-.78571429,.03348214){112}{\line(-1,0){.78571429}}
\multiput(213.5,36)(-.78571429,.03348214){112}{\line(-1,0){.78571429}}
\multiput(86.75,53.75)(-.081649832,.033670034){297}{\line(-1,0){.081649832}}
\multiput(195.25,53.75)(-.081649832,.033670034){297}{\line(-1,0){.081649832}}
\multiput(86.5,54)(-.033723022,-.074190647){556}{\line(0,-1){.074190647}}
\multiput(195,54)(-.033723022,-.074190647){556}{\line(0,-1){.074190647}}
\multiput(39.25,17.25)(.21455224,-.03358209){134}{\line(1,0){.21455224}}
\multiput(147.75,17.25)(.21455224,-.03358209){134}{\line(1,0){.21455224}}
\multiput(86.25,54)(-.0428899083,-.0337155963){1090}{\line(-1,0){.0428899083}}
\multiput(194.75,54)(-.0428899083,-.0337155963){1090}{\line(-1,0){.0428899083}}
\multiput(17.25,40)(.0625,-.0337252475){808}{\line(1,0){.0625}}
\multiput(125.75,40)(.0625,-.0337252475){808}{\line(1,0){.0625}}
\multiput(39.25,17.25)(.38297872,.03368794){141}{\line(1,0){.38297872}}
\multiput(147.75,17.25)(.38297872,.03368794){141}{\line(1,0){.38297872}}
\multiput(104.75,35.75)(-.033653846,-.038461538){364}{\line(0,-1){.038461538}}
\multiput(213.25,35.75)(-.033653846,-.038461538){364}{\line(0,-1){.038461538}}
\multiput(63,63.25)(.0337370242,-.0478662053){867}{\line(0,-1){.0478662053}}
\multiput(171.5,63.25)(.0337370242,-.0478662053){867}{\line(0,-1){.0478662053}}
\put(53.75,4.75){\makebox(0,0)[cc]{$K_1$}}
\put(162,4.75){\makebox(0,0)[cc]{$K_2$}}
\end{picture}
\end{figure}

\end{center}

These graphs satisfy the following property.

\begin{itemize}
\item[$(R_{11})$] There exists a subset $X$ in $K^0$ with $\Lk(X)\ne K^0$ and a vertex $x\in \Lk(X)$ such that the sets $L=[\Lk(X), \Lk(x)]$,
$\Lk(x) \setminus X$ satisfy (*) and there exists a vertex $z\in K^0 \setminus (X \cup \Lk(X) \cup \{x\})$ such that
    for some $t\in lk(z)\setminus (\Lk(X)\cup \{x\})$, the graph $K\setminus \{(z,t)\}$ is excluded and the commutator $[K^0\setminus (X\cup \{x,z\}), K^0\setminus (X\cup \{x,z\})]$ satisfies (*).
\end{itemize}

For the graph $K_1$, one can take $X=\{6,7\}$, $x=8$, $z=4$, $t=5$; for the graph $K_2$, one can take $X=\{2,4\}$, $x=7$, $z=8$, $t=1$.

Suppose that a graph $K$ satisfies ($R_{11}$). Then in any faithful $K$-dissection
diagram $\Delta$, an $x$-curve cannot cross the boundary of the subsurface $S(X)$
constructed as in the proof of Lemma \ref{BasicLemma1} (taking the regular
neighborhood of the graph formed by the $X$-curves, and gluing in the disc components of the complement). Indeed, otherwise we would get a curve whose effective content
satisfies (*) and get a contradiction as in Step 4. Hence any $x$-curve must be
either inside $S(X)$, or outside it.

Since $\Lk(x) \setminus X$ satisfies (*), no $x$-curve can be outside $S(X)$. Hence the content of the surface $S\setminus S(X)$ does not intersect $X\cup \{x\}$. Since $z\not\in X \cup S(X)$, $z$-curves are outside $S(X)$. Since there must be a $z$-curve intersecting a $t$-curve, $t\not\in \Lk(X)\setminus\{x\}$ (the graph $K\setminus \{(z,t)\}$ is excluded), not all $z$-curves can be parallel to the boundary of $S(X)$. Hence one of the components $S'$ of $S\setminus (S(X)\cup S(\{z\}))$ has a non-Abelian fundamental group and content in $K^0\setminus (X\cup \{x,z\})$. The commutator of two intersecting curves $\alpha$ and $\beta$ in $S'$ has effective content in the commutator $[K^0\setminus (X\cup \{x,z\}), K^0\setminus (X\cup \{x,z\})]$ (by Lemma \ref{commeff}) which satisfies (*), a contradiction.

Of 14 remaining graphs 13 satisfy ($R_{11}$).


{\bf Step 8.} The remaining graph is on the following picture.

That graph satisfies the following condition.

\begin{itemize}
\item[$(R_{12})$] There are two pairs of adjacent vertices $(x,y)$, $(z,t)$ such that

\begin{itemize}
\item[(i)] $x$ and $z$ are adjacent and the graph $K\setminus \{(x,z)\}$ is excluded;
\item[(ii)] the commutator $[\Lk(\{x,y\}), \Lk(\{z,t\})]$ satisfies (*);
\item[(iii)] $y\not\in \{z,t\}\cup \Lk(\{z,t\})$, $t\not\in \{x,y\}\cup \Lk(\{x,y\}$;
\item[(iv)] the set $\Lk(x)\cap (\{z,t\}\cup \Lk(\{z,t\}))$ is dense in $G$;
\item[(v)] the graph $K\setminus \{(z,t)\}$ is excluded;
\item[(vi)] the commutator of the set $\{x,y\}\cup \Lk(\{x,y\}) \setminus \{z\}$ with itself satisfies (*).
\end{itemize}
\end{itemize}

Indeed, one can take $(x,y)=(6,8), (z,t)=(3,4)$.

A graph $K$ satisfying $(R_{12})$ can be excluded. Indeed, suppose that a faithful $K$-dissection diagram exists on a hyperbolic surface $S$. Consider subsurfaces $S_1=S(\{x,y\})$ and $S_2=S(\{z,t\})$ (as in the previous step).
Their boundaries cannot intersect by (ii). Since by (i), there exists an $x$-curve that intersects $z$-curve, either a connected component $S'$ of $S_1$ is inside $S_2$ or a connected component $S''$ of $S_2$ is inside $S_1$. By (iii) $S'$ must be a regular neighborhood of an $x$-curve, $S''$ must be a regular neighborhood of a $z$-curve.

\begin{center}
\unitlength .6mm 
\linethickness{0.4pt}
\ifx\plotpoint\undefined\newsavebox{\plotpoint}\fi 
\begin{picture}(107.5,64.75)(0,0)
\put(107.5,34.25){\makebox(0,0)[cc]{1}}
\put(87.75,54){\makebox(0,0)[cc]{2}}
\put(63.5,64.75){\makebox(0,0)[cc]{3}}
\put(32.5,59.5){\makebox(0,0)[cc]{4}}
\put(14.75,38.5){\makebox(0,0)[cc]{5}}
\put(36.25,15){\makebox(0,0)[cc]{6}}
\put(67.5,9.25){\makebox(0,0)[cc]{7}}
\put(91.75,18){\makebox(0,0)[cc]{8}}
\multiput(105.25,35.75)(-.051183432,.0337278107){845}{\line(-1,0){.051183432}}
\multiput(105,36)(-.1044721408,.0337243402){682}{\line(-1,0){.1044721408}}
\multiput(105.25,35.75)(-.119154676,-.033723022){556}{\line(-1,0){.119154676}}
\multiput(105.25,35.75)(-.0542521994,-.0337243402){682}{\line(-1,0){.0542521994}}
\multiput(87.25,53.5)(-.32469512,.03353659){164}{\line(-1,0){.32469512}}
\multiput(87,53.75)(-.16607565,-.033687943){423}{\line(-1,0){.16607565}}
\multiput(62.75,63.5)(-.21268657,-.03358209){134}{\line(-1,0){.21268657}}
\multiput(62.5,63.5)(-.0646306818,-.0337357955){704}{\line(-1,0){.0646306818}}
\multiput(62.5,63.75)(-.0336956522,-.0670289855){690}{\line(0,-1){.0670289855}}
\multiput(34.25,59)(-.03371869,-.036608863){519}{\line(0,-1){.036608863}}
\multiput(105,36)(-.78571429,.03348214){112}{\line(-1,0){.78571429}}
\multiput(86.75,53.75)(-.081649832,.033670034){297}{\line(-1,0){.081649832}}
\multiput(86.5,54)(-.033723022,-.074190647){556}{\line(0,-1){.074190647}}
\multiput(39.25,17.25)(.21455224,-.03358209){134}{\line(1,0){.21455224}}
\multiput(86.25,54)(-.0428899083,-.0337155963){1090}{\line(-1,0){.0428899083}}
\multiput(17.25,40)(.0625,-.0337252475){808}{\line(1,0){.0625}}
\multiput(39.25,17.25)(.38297872,.03368794){141}{\line(1,0){.38297872}}
\multiput(104.75,35.75)(-.033653846,-.038461538){364}{\line(0,-1){.038461538}}
\put(53.75,4.75){\makebox(0,0)[cc]{}}
\multiput(68,13)(.095192308,.033653846){260}{\line(1,0){.095192308}}
\end{picture}
\end{center}

Suppose that an $x$-curve is in $S_2$. Then its content must be in $\Lk(x)\cap (\{z,t\}\cup \Lk(\{z,t\}))$ which is dense in $G$, a contradiction. Therefore $S_2$ does not contain $x$-curves. By (v), there exists a component $Y$ of $S_2$ with non-Abelian fundamental group. Then $Y$ cannot be inside $S_1$, so $Y$ does not contain $x$-curves or arcs. Therefore the content of $Y$ is in $\{x,y\}\cup \Lk(\{x,y\}) \setminus \{z\}$. Taking the commutator of two intersecting curves in $Y$ and using (vi), we get a contradiction.

The proof is complete.
\endproof

\noindent John Crisp:
{\small\sc \\ I.M.B.(UMR 5584 du CNRS), Universit\'e de Bourgogne, B.P. 47 870, 21078 Dijon, France.\\}
{\it E-mail: } {\tt john.crisp@gmail.com}

\vspace{3mm}

\noindent Michah Sageev:
{\small\sc \\ Technion, Israel University of Technology, Dept. of Mathematics, Haifa 32000, Israel.\\}
{\it E-mail: } {\tt sageevm@techunix.technion.ac.il}

\vspace{3mm}

\noindent Mark V. Sapir:
{\small\sc \\ Department of Mathematics, Vanderbilt University, Nashville, TN 37240.\\}
{\it E-mail: } {\tt m.sapir@vanderbilt.edu}

\end{document}